\documentclass[10pt,final,onecolumn]{siamltex}

\usepackage{amsmath,amsfonts,amssymb,amsopn}
\usepackage{graphicx}
\usepackage{ifpdf}
\ifpdf
  \usepackage[pdftex,colorlinks=true,breaklinks=true,citecolor=blue]{hyperref}
\else
  \usepackage[ps2pdf,colorlinks=true,citecolor=blue]{hyperref}
\fi
\usepackage{xcolor}
\usepackage{cite}
\usepackage{xspace}
\usepackage{subfig} 
\usepackage{mathtools} 
\usepackage{calrsfs} 

\usepackage{relsize}

\DeclareFontFamily{OT1}{pzc}{}
\DeclareFontShape{OT1}{pzc}{m}{it}%
              {<-> s * [1.25] pzcmi7t}{}
\DeclareMathAlphabet{\mathpzc}{OT1}{pzc}%
                                 {m}{it}

\newcommand{\cpfb}{\mathpzc{X}} 
\newcommand{\cpfl}{\mathpzc{Z}} 

\newcommand{\fbc}[3]{{#1}_{#2}^{#3}}
\newcommand{\lballs}[1]{\mathcal{H}_{#1}}
\newcommand{\lballss}[2]{\mathcal{\tilde H}_{#1#2}}
\newcommand{\bll}{{\mathbb{B}^3}} 
\newcommand{\rplus}{{\mathbb{R}^{+}}}
\newcommand{\bv}[1]{\boldsymbol{#1}}
\newcommand{\dfn}{\triangleq}
\newcommand{\untsph}{{\mathbb{S}^{2}}} 

\newcommand{\shc}[2]{{#1}_{#2}}
\newcommand{\lsph}{L^2(\untsph)}
\newcommand{\lball}{L^2(\bll)}
\newcommand{\lrplus}{L^2(\rplus)}

\newcommand{\rsphere}{{{R}_{\rm s}}}

\newcommand{\rradial}{{{R}_{\rm r}}}

\newcommand{\lballr}{\mathcal{H}_R}
\newcommand{\lballp}{\mathcal{H}_{PL}}
\newcommand{\lballk}{\mathcal{H}_{L}^{K}}
\newcommand{\shcc}[2]{{#1}_{#2}}

\newcommand{\conj}[1]{\left({#1}\right)^\ast} 
\newcommand{\unit}[1]{\bv{\hat{#1}}}

\newcommand{\intball}{\int_{\bll}}

\newcommand{\figref}[1]{Fig.\,\ref{#1}}

\newcommand{\secref}[1]{Section\,\ref{#1}}

\newcommand{\ddv}[1]{\ensuremath{\textrm{d}^3 {\mu}({#1})}}
\newcommand{\dds}[1]{\ensuremath{\textrm{d}^2 {\nu}({#1})}}
\newcommand{\ddrr}[1]{\textrm{d}{\upsilon}({#1})}
\newcommand{\ddr}[1]{\textrm{d}{#1}}

\newtheorem{remark}[theorem]{Remark}

\graphicspath{{/},{figs/},{pdf/},{Figures/}}

\DeclarePairedDelimiterX\norm[1]{\lVert}{\rVert}{#1}
\DeclarePairedDelimiterX\abs[1]{|}{|}{#1}
\DeclarePairedDelimiterX\innerp[2]{\langle}{\rangle}{#1,#2}
\DeclarePairedDelimiterX\parn[1]{(}{)}{#1}

\begin{document}

\title{Slepian Spatial-Spectral Concentration on the Ball}

\author{Zubair Khalid\footnotemark[2]
\and Rodney A. Kennedy\footnotemark[2]
\and Jason D. McEwen\footnotemark[3]}

\renewcommand{\thefootnote}{\fnsymbol{footnote}}

\maketitle

\footnotetext[2]{Research School of Engineering, College of Engineering and
   Computer Science, The Australian National University, Canberra,
   Australia. {E-mail: zubair.khalid@anu.edu.au, rodney.kennedy@anu.edu.au}}
\footnotetext[3]{Mullard Space Science Laboratory, University College
  London, Surrey RH5 6NT, UK. {E-mail: jason.mcewen@ucl.ac.uk }}
\footnotetext[2]{Supported by the Australian Research Council's Discovery Projects funding scheme (project no. DP1094350).}
\footnotetext[3]{Supported in part by a Newton International Fellowship from the Royal Society and the British Academy.}

\begin{abstract}
  We formulate and solve the Slepian spatial-spectral concentration
  problem on the three-dimensional ball.  Both the standard
  Fourier-Bessel and also the Fourier-Laguerre spectral domains are
  considered since the latter exhibits a number of practical advantages
  (spectral decoupling and exact computation).  The Slepian spatial
  and spectral concentration problems are formulated as eigenvalue
  problems, the eigenfunctions of which form an orthogonal family of
  concentrated functions.  Equivalence between the spatial and
  spectral problems is shown.  The spherical Shannon number on the
  ball is derived, which acts as the analog of the space-bandwidth
  product in the Euclidean setting, giving an estimate of the number
  of concentrated eigenfunctions and thus the dimension of the space
  of functions that can be concentrated in both the spatial and
  spectral domains simultaneously.  Various symmetries of the spatial
  region are considered that reduce considerably the computational
  burden of recovering eigenfunctions, either by decoupling the problem
  into smaller subproblems or by affording analytic calculations.  The
  family of concentrated eigenfunctions forms a Slepian basis that
  can be used be represent concentrated signals efficiently.  We
  illustrate our results with numerical examples and show that the
  Slepian basis indeeds permits a sparse representation of
  concentrated signals.
\end{abstract}

\section{Introduction}

It is well-known that functions cannot have finite support in both the spatial~(or time) and spectral~(or frequency) domain
at the same time~\cite{Slepian:1960,Slepian:1983}. This fundamental
problem of finding and representing the functions that are optimally
energy concentrated in both the time and frequency domains was solved by
Slepian, Landau and Pollak in the early
1960s~\cite{Slepian:1960,Landau:1961,Landau:1962,Slepian:1965}. This
problem, herein referred to as the \emph{Slepian spatial-spectral concentration problem}, or \emph{Slepian concentration problem} for short, gives rise to the orthogonal families of functions that are optimally concentrated in the spatial (spectral) domain and exactly limited in the spectral (spatial) domain.  These families of functions and their multidimensional extensions~\cite{Slepian:1964} have been extensively used in various branches of science and engineering~(e.g., signal processing~\cite{Thomson:1982,Mathew:1985}, medical imaging~\cite{Jackson:1991}, geophysics~\cite{Thomson:1976}, climatology~\cite{Thomson:1990}, to name a few).

Although the the Slepian spatial-spectral concentration problem was
initially formulated and solved in the Euclidean domain,
generalizations for various geometries and connections to wavelet
analysis have also been well-studied~(e.g.,
\cite{Meaney:1984,Daubechies:1988,Cohen:1989,Daubechies:1990,Albertella:1999,Mortlock:2002,Fernández:2003,Wieczorek:2005,Simons:2006}). We
note that the Slepian concentration problem for functions defined on
the two-sphere $\untsph$ has been thoroughly revisited and
investigated~\cite{Albertella:1999,Simons:2006,Wieczorek:2005}. The
resulting orthogonal family of band-limited spatially concentrated
functions have been applied for localized spectral
analysis~\cite{Wieczorek:2007} and spectral
estimation~\cite{Dahlen:2008} of signals~(finite energy functions)
defined on the sphere. There are also many
applications~\cite{Leistedt:2012,Simons:2011_1,Simons:2011_2} where
signals or data are defined naturally on the three-dimensional
ball, or ball for short. For example, signals defined on the ball arise
when observations made on the sphere are augmented with radial
information, such as depth, distance or redshift. Recently, a number
of signal processing techniques have been tailored and extended to
deal with signals defined on the ball (e.g.,
\cite{Simons:2011_1,Lanusse:2012,Leistedt:2012}).

In this paper, we pose, solve and analyse the Slepian concentration
problem of simultaneous spatial and spectral localization of functions
defined on the ball. By considering Slepian's quadratic~(energy)
concentration criterion, we formulate and solve the problems to: (1)
find the band-limited functions with maximum concentration in some
spatial region; and (2) find the space-limited functions with maximum
concentration in some region of the spectral domain. Each problem is
formulated as an eigenvalue problem, the solution of which gives the
orthogonal family of functions, referred as eigenfunctions, which are
either spatially concentrated while band-limited, or spectrally
concentrated while space-limited.  These eigenfunctions serve as an
alternative basis on the ball, which we call a Slepian basis, for the
representation of a band-limited or space-limited signal. We show, and
also illustrate through an example, that the representation of
band-limited spatially concentrated or space-limited spectrally
concentrated functions is sparse in the Slepian basis, which is the
essence of the Slepian spatial-spectral concentration problem. We also
derive the spherical Shannon number as an equivalent of the Shannon
number in the one dimensional Slepian concentration
problem~\cite{Slepian:1965,Percival:1993}, which serves as an estimate
of the number of concentrated functions in the Slepian basis.

For the spectral domain characterization of functions defined on the
ball we use two basis functions: (1) spherical harmonic-Bessel
functions, which arise as a solution of Helmholtz's equation in
three-dimensional spherical coordinates, and are referred to as
\emph{Fourier-Bessel}\footnote{A more appropriate terminology would be
  spherical harmonic-Bessel basis, however we adopt the established
  convention of using the term Fourier to denote the spherical
  harmonic part.} basis functions; and (2) spherical harmonic-Laguerre
functions, which are referred to as \emph{Fourier-Laguerre} basis
functions. We consider the Fourier-Laguerre functions in addition to
the standard Fourier-Bessel functions, as the Fourier-Laguerre
functions serve as a complete basis for signals defined on the ball,
enable the decoupling of the radial and angular components of the
signal, and support the exact computation of forward and inverse
Fourier-Laguerre transforms~\cite{Leistedt:2012}. We show that the
eigenvalue problem to find the eigenfunctions or Slepian basis can be
decomposed into subproblems when the spatial region of interest is
symmetric in nature. We consider two types of symmetric regions: (1)
circularly symmetric regions; and (2) circularly symmetric and
radially independent regions.

As Slepian functions on the one-dimensional Euclidean
domain~\cite{Slepian:1960,Landau:1961,Landau:1962,Slepian:1965}, and
other
geometries~\cite{Slepian:1964,Meaney:1984,Albertella:1999,Fernández:2003,Simons:2006},
have been widely useful in a diverse variety of applications, we hope
that the proposed orthogonal family of Slepian eigenfunctions on the
ball will find similar applications in fields such as cosmology,
geophysics and planetary science, where data/signals are often
inherently defined on the ball.  For example, the band-limited
spatially concentrated eigenfunctions can be used as window functions
to develop multi-window spectral estimation
techniques~\cite{Thomson:1982,Thomson:1990,Dahlen:2008,Wieczorek:2005}
for the estimation of the signal spectrum from observations made over
the limited spatial region.

We organize the remainder of the paper as follows. The mathematical
preliminaries for functions on the ball are presented in
\secref{sec:maths}. The Slepian concentration problem is posed as an
eigenvalue problem in \secref{sec:conc_problem} and the resulting
eigenfunctions are analysed in \secref{sec:eignfunctions_analysis}.
The decomposition of the eigenvalue problem into subproblems for the
case of special, but important, symmetric spatial regions is presented
in \secref{sec:special_regions}. The representation of spatially
concentrated band-limited functions in the Slepian basis is discussed
and illustrated in \secref{sec:Applications}.  Concluding remarks are
made in \secref{sec:conclusions}.

\section{Mathematical Preliminaries}
\label{sec:maths}

We review the mathematical background of signals defined on the ball in this section.  After defining coordinate systems, measures and inner products, we then review harmonic analysis on the ball, focusing on both the Fourier-Bessel and Fourier-Laguerre settings.  We conclude this section by reviewing important subspaces and operators related to the ball.

\subsection{Signals on the Sphere and Ball}

We define the ball by \mbox{$\bll \dfn \rplus\times \untsph$}, where
$\rplus$ denotes the domain $[0,\infty)$ on the real line and $\untsph
\dfn \{\bv{y} \in \mathbb{R}^3\colon \|\bv{y}\| = 1 \}$ denotes the
unit sphere. A vector $\bv{r}\in\bll$ can be
represented in spherical coordinates as
$\bv{r}\equiv\bv{r}(r,\theta,\phi)\dfn (r\sin\theta\cos\phi,\,
r\sin\theta\sin\phi,\, r\cos \theta)^{\rm T}$, where
$(\cdot)^{\rm T}$ denotes matrix or vector transpose. Here,
$r\dfn\|\bv{r}\|\in[0,\infty)$ represents the Euclidean norm of
$\bv{r}$, $\theta\in [0, \pi]$ represents the co-latitude or
elevation measured with respect to the positive $z$-axis and
$\phi\in [0, 2\pi)$ represents the longitude or azimuth and is
measured with respect to the positive $x$-axis in the $x$-$y$ plane.
The unit norm vector $ \unit{r} \equiv\unit{r}(\theta,\phi)\dfn
\bv{r}/\|\bv{r}\|= (\sin\theta\cos\phi,\, \sin\theta\sin\phi,\, \cos
\theta)^{\rm T} \in \mathbb{R}^3$ represents a point on the unit sphere
$\untsph$.

The space of square integrable complex-valued functions defined on
$\rplus$, $\untsph$ and $\bll$  form Hilbert spaces, denoted by
$\lrplus$, $\lsph$ and $\lball$, respectively, equipped with the
inner products defined by
\begin{align}
\label{eqn:innprd}
    \innerp[\big]{f}{g}_\rplus &\dfn  \int_\rplus \ddrr{r} f(r)\,g^\ast(r),
     \\ \label{eqn:innprd_2}
    \innerp[\big]{f}{g}_\untsph &\dfn  \int_\untsph \dds{\unit{r}}    f(\unit{r})\,g^\ast(\unit{r}),
    \\ \label{eqn:innprd_3}
    \innerp[\big]{f}{g}_\bll &\dfn  \intball \ddv{\bv{r}} f(\bv{r})\,g^\ast(\bv{r})
    ,
\end{align}
where $f,\,g$ are functions respectively defined on $\rplus$, $\untsph$ and $\bll$ in \eqref{eqn:innprd}, \eqref{eqn:innprd_2} and \eqref{eqn:innprd_3},  $\ddrr{r} = r^2 {\rm d}r$, $\dds{\unit{r}} = \sin \theta\,{\rm d} \theta\,{\rm d} \phi$ and $\ddv{\bv{r}} = r^2\sin \theta\,{\rm d} r\,{\rm d} \theta\,{\rm d} \phi$ represents infinitesimal length, area and the volume element respectively, $(\cdot)^\ast$ denotes complex conjugation and the integration is carried out over the respective domain. The inner products in \eqref{eqn:innprd}, \eqref{eqn:innprd_2} and \eqref{eqn:innprd_3} induce norms $\|f\| \dfn\langle f,f \rangle^{1/2}$. Throughout this paper, the functions with finite induced norm belonging to one of these spaces are referred to as signals.

\subsection{Harmonic Analysis on the Ball}

We review harmonic analysis on the ball, starting with the spherical Bessel and Laguerre transforms on the positive real line $\rplus$ and the spherical harmonic transform on the unit sphere $\untsph$, before combining these to recover the Fourier-Bessel and Fourier-Laguerre transforms on the ball, respectively.

\subsubsection{Spherical Bessel Transform}
The spherical Bessel functions, which arise as radial solutions to the Helmholtz equation in spherical coordinates, form a basis for functions on the non-negative real line $\rplus$. In this work, we consider spherical Bessel functions of the first kind, denoted by $j_\ell$ defined on $\rplus$, where $\ell$ denotes the order. The spherical Bessel functions satisfy the closure relation~\cite{Watson:1995}
\begin{align}
\int_\rplus \ddrr{r} j_\ell(k r)  j_\ell(k' r)
= \frac{\pi}{2 k^2} \, \delta(k - k')
,
\end{align}
for $r \in \rplus$ and $k \in \rplus$, and where $\delta(k-k')$ denotes the one-dimensional Dirac delta.  Consequently, we can represent a signal $f\in\lrplus$ using the following $\ell$-th order spherical Bessel inverse and forward transform, respectively,
\begin{align}
f(r) = \sqrt{\frac{2}{\pi}} \, \int_\rplus \ddr{k} f_\ell(k) k j_\ell(k r) \quad \mbox{with} \quad f_\ell(k)  \dfn \sqrt{\frac{2}{\pi}}\,  \int_\rplus \ddrr{r} f(r) k j_\ell(k r),
\end{align}
where $f_\ell(k)$ denotes the spherical Bessel trasnform.

\subsubsection{Spherical Laguerre Transform}
The Laguerre polynomials, solutions to the Laguerre differential equation~\cite{Pollard:1947,Weniger:2008}, are well known for their various applications, notably in the quantum mechanical treatment of the hydrogen atom \cite{Dunkl:2003}, and form a basis for functions on the interval $\rplus$.  We adopt the spherical Laguerre transform and associated normalisation presented by \cite{Leistedt:2012}, defining the spherical Laguerre basis functions of
non-negative integer radial degree $p$ by
\begin{align}\label{Eq:Laguerre_dfn}
K_p(r) \dfn \sqrt{\frac{p!}{(p+2)!}}\,e^{-r/2}\,L_p^{(2)}(r),
\end{align}
where $L_p^{(2)}(r)$ represents the $p$-th generalized Laguerre
polynomial of second order, defined by
\begin{align}\label{Eq:Laguerre_two_dfn}
L_p^{(2)}(r) \dfn \sum_{j=0}^p \binom{p+2}{p-j}\frac{(-r)^j}{j!}.
\end{align}
Since we use the spherical Laguerre basis functions for the
expansion of signals defined on $\rplus$ with differential measure
$\ddrr{r} = r^2 {\rm d}r$, we have chosen the second order generalized Laguerre
polynomial. The basis functions $K_p(r)$ in \eqref{Eq:Laguerre_dfn}
are orthonormal on $\rplus$, that is,
$\innerp[\big]{K_p}{K_q}_\rplus = \delta_{pq}$,
where $\delta_{pq}$ denotes the Kronecker
delta. The spherical Laguerre polynomials defined in
\eqref{Eq:Laguerre_dfn} serve as complete basis functions on
$\rplus$, where the completeness stems from the completeness of
generalized Laguerre polynomials and therefore we can expand a
signal $f\in\lrplus$ using the following spherical Laguerre inverse
and forward transform, respectively,
\begin{align}
f(r) = \sum_{p=0}^{\infty} f_p K_p(r) \quad \mbox{with} \quad f_p \dfn
\innerp[\big]{f}{K_p}_\rplus,
\end{align}
where $f_p$ denotes the spherical Laguerre coefficient of radial degree $p$.

\subsubsection{Spherical Harmonic Transform}
The spherical harmonic functions, which arise as angular solutions to the Helmholtz equation in spherical coordinates, are denoted $Y_{\ell}^m(\unit{x}) = Y_{\ell}^m(\theta, \phi)$, for integer degree ${\ell} \geq 0$ and integer order $ |m| \leq {\ell}$ and are defined by~\cite{Colton:1998,Sakurai:1994,Kennedy-book:2013}
\begin{equation}
\label{Eq:Sph_harmonics} Y_{\ell m}(\unit{r})= Y_{\ell m}(\theta,
\phi) =
\sqrt{\frac{2{\ell}+1}{4\pi}\frac{({\ell}-m)!}{({\ell}+m)!}}\,
        P_{\ell}^{m}(\cos\theta)e^{im\phi},
\end{equation}
where $P_{\ell}^{m}$ denotes the associated Legendre function of
degree $\ell$ and order $m$ with Condon-Shortley phase
included~\cite{Kennedy-book:2013}. With the above definition, the spherical
harmonic functions, or simply the spherical harmonics, form a complete
orthonormal basis for $\lsph$ and therefore a signal $f\in\lsph$
can be expanded using the following spherical harmonic inverse and
forward transform, respectively,
\begin{align}
f(\unit{r}) = \sum_{\ell=0}^\infty \sum_{m=-\ell}^\ell f_{\ell m}
Y_{\ell m}(\unit{r}) \quad \mbox{with} \quad f_{\ell m} \dfn \innerp[\big]{f}{Y_{\ell
m}}_\untsph,
\end{align}
where $f_{\ell m}$ denotes the spherical harmonic coefficient of
angular degree $\ell$ and order $m$. We note that the spherical
harmonic functions follow the conjugate symmetry property
$Y_{\ell m}(\unit{r}) = (-1)^m Y_{\ell (-m)}^\ast(\unit{r})$. We
further note that the function $f$ is real valued if its
spherical harmonic coefficients satisfy the
conjugate symmetry property $f_{\ell m} = (-1)^m f_{\ell
(-m)}^\ast$.

\subsubsection{Fourier-Bessel Transform}

We define the Fourier-Bessel functions by~\cite{Castro:2005,Leistedt:2012b,Abramo:2010}
\begin{align}\label{Eq:FB_definition}
\fbc{X}{\ell m}{}(k,\bv{r})  \triangleq \sqrt\frac{2}{\pi}\,k\,
j_\ell(k r) Y_{\ell m}(\theta,\phi), \quad \bv{r} =
\bv{r}(r,\theta,\phi).
\end{align}
With the above definition, the Fourier-Bessel
functions form a complete, orthogonal basis for $\lball$, satisfying the orthogonality relation
\begin{align}
\intball \ddv{\bv{r}} \fbc{X}{\ell m}{}(k,\bv{r}) \fbc{X}{\ell'
m'}{\ast} (k',\bv{r}) = \delta(k-k') \delta_{\ell \ell'}\delta_{m
m'}.
\end{align}
By
the completeness of the Fourier-Bessel functions, a signal $f \in \lball$
can be decomposed in the Fourier-Bessel basis by
\begin{align}\label{Eq:FB_inverse_transform}
f(\bv{r}) &=  \int_{\rplus} \ddr{k} \sum_{\ell=0}^\infty
\sum_{m=-\ell}^\ell {f}_{\ell m}(k) X_{\ell m}(k,\bv{r}) \nonumber
\\ &= \sqrt\frac{2}{\pi}  \int_{\rplus} \ddr{k} \sum_{\ell=0}^\infty
\sum_{m=-\ell}^\ell f_{\ell m}(k) k j_\ell(k r) Y_{\ell
m}(\theta,\phi)
\end{align}
where $\fbc{f}{\ell m}{}(k)$ denotes the Fourier-Bessel coefficient, of
degree $\ell$, order $m$ and radial harmonic variable $k\in\rplus$, given by
\begin{align}
\fbc{f}{\ell m}{}(k) \triangleq \sqrt{\frac{2}{\pi}}\intball
\ddv{\bv{r}} f(\bv{r}) k j_\ell(k r) Y_{\ell m}^\ast(\theta,\phi).
\end{align}
The Fourier-Bessel coefficients constitute a spectral domain
representation of signals defined on the ball. Such a spectral domain
is to referred as the \emph{Fourier-Bessel spectral domain}.

The Fourier-Bessel transform is the natural harmonic transform on the ball since the Fourier-Bessel functions are the eigenfunctions of the spherical Laplacian and thus the Fourier-Bessel transform corresponds to the standard three-dimensional Fourier transform in spherical coordinates. However, the Fourier-Bessel transform suffers from a number of practical limitations \cite{Leistedt:2012}, motivating alternative harmonic representations of the ball, such as the Fourier-Laguerre transform.

\subsubsection{Fourier-Laguerre Transform}

In the Fourier-Bessel transform, the spherical Bessel
functions are used for the decomposition of a signal along the radial line
$\rplus$. Alternatively, we can use the spherical Laguerre basis
functions for the expansion of a signal along $\rplus$. Combining the
spherical Laguerre basis functions and spherical harmonic functions,
we define the Fourier-Laguerre basis functions for a signal as
\begin{align}
Z_{\ell m p}(\bv{r}) \dfn K_p(r)Y_\ell^m(\theta,\phi),\quad \bv{r} =
\bv{r}(r,\theta,\phi).
\end{align}
By the completeness of both Laguerre polynomials and spherical
harmonics, any signal $f\in \lball$ can be expanded as \cite{Leistedt:2012}
\begin{equation}
\label{Eq:f_expansion}
    f(\bv{r}) = \sum_{p=0}^{\infty} \sum_{{\ell}=0}^{\infty}
        \sum_{m=-{\ell}}^{\ell} \shc{f}{\ell m p}  Z_{\ell m p}(\bv{r}),%
\end{equation}
where $\shc{f}{\ell m p}$ is the Fourier-Laguerre coefficient of
radial degree $p$, angular degree $\ell$ and angular order $m$, and
is obtained by the Fourier-Laguerre transform
\begin{equation}\label{Eq:fcoeff}
    \shc{f}{\ell m p}\triangleq \langle f,Z_{\ell m p} \rangle_\bll = \intball \ddv{\bv{r}} f(\bv{r}) Z_{\ell m p}^\ast(\bv{r}).%
\end{equation}
 These Fourier coefficients constitute another spectral domain representation of signals defined on the ball, which we refer to as the \emph{Fourier-Laguerre spectral domain}.

The Fourier-Laguerre transform exhibits a number of practical advantages over the Fourier-Bessel transform, namely: (1) the angular and radial components of signals are decoupled in harmonic space; and (2) exact quadrature can be developed, leading to theoretically exact forward and inverse Fourier-Laguerre transforms \cite{Leistedt:2012}.

\subsubsection{Dirac Delta on the Ball}
The Dirac delta function on the ball is defined by
\begin{align}
\delta(\bv{r},\bv{r}') \dfn \left(r^2 \sin\theta\right)^{-1} \delta(r-r') \delta(\theta-\theta')\delta(\phi-\phi'),
\end{align}
and satisfies the sifting property
\begin{align}
\int_{\bll} \ddv{\bv{r}'} f(\bv{r}') \delta(\bv{r},\bv{r}')
\,=\,f(\bv{r}).
\end{align}
The Dirac delta has the following expansion in terms of Fourier-Bessel basis
functions
\begin{align}\label{Eq:dirac_expansion_FB}
\delta(\bv{r},\bv{r}') \,&=\,
\sum_{\ell=0}^{\infty}\,\sum_{m=-\ell}^\ell\,\int_{\rplus}\ddr{k}
\fbc{X}{\ell m}{}(k,\bv{r})\, \fbc{X}{\ell m}{\ast}(k,\bv{r}')
\nonumber
\\ &= \frac{1}{2\pi^2} \sum_{\ell=0}^{\infty}\,\left(2\ell+1\right)
P_\ell^0(\unit{r}\cdot\unit{r}') \left(\int_{\rplus} \ddr{k}
\, k^2 j_\ell(k r) j_\ell(k r)  \right),
\end{align}
and has the following expansion in terms of Fourier-Laguerre basis
functions
\begin{align}\label{Eq:dirac_expansion_FL}
\delta(\bv{r},\bv{r}') &= \sum_{p=0}^{\infty}
\sum_{\ell=0}^{\infty}\,\sum_{m=-\ell}^\ell\, Z_{\ell m p}(\bv{r})\,
Z^\ast_{\ell m p}(\bv{r}') \nonumber \\ &= \sum_{p=0}^{\infty}
K_p(r)K_p(r') \sum_{\ell=0}^{\infty} \frac{2\ell+1}{4\pi}
P_\ell^0(\unit{r}\cdot\unit{r}'),
\end{align}
where $\unit{r}\cdot\unit{r}'$ denotes the three dimensional dot product between unit vectors $\unit{r}$ and $\unit{r}'$ and we have noted the addition theorem for the spherical harmonics.

\subsection{Important Subspaces of $\lball $}

Define $\lballss{K}{L}$ as the space of band-limited
functions such that the signal is band-limited in the Fourier-Bessel
spectral domain within the spectral region $\tilde A_{KL}\dfn\{0\leq k
\leq K, 0\leq \ell \leq L-1\}$ for $L\in \mathbb{Z}^+$ and $K \in
\rplus$. Due to the continuous Fourier-Bessel spectral domain $k$, $\lballss{K}{L}$ is an
infinite dimensional subspace of $\lball$.
Define $\lballs{PL}$ as the space of band-limited functions such that
the signal is band-limited in the Fourier-Laguerre spectral domain
within the spectral region $A_{PL}\dfn\{0\leq p \leq P-1, 0\leq
\ell \leq L-1\}$ for $P,L\in \mathbb{Z}^+$. $\lballs{PL}$ is a
finite dimensional subspace of $\lball$ with size $PL^2$.
Also define $\lballs{R}$
as the space of finite energy space-limited functions confined within the region
$R\subset \bll$. $\lballs{R}$ is an infinite dimensional subspace
of $\lball$.

\subsection{Important Operators of $\lball $}

Define an operator $\rm{S}$ for signals on the ball by the
general Fredholm integral equation~\cite{Kennedy-book:2013}
\begin{align}\label{Eq:gen_operator_def}
({\rm S} f)(\bv{r})\,=\,\intball \ddv{\bv{r'}} S(\bv{r},
\bv{r'})\,f(\bv{r'}),
\end{align}
where $S(\bv{r}, \bv{r'})$ is the kernel for an operator
$\rm{S}$ defined on $\bll\times\bll$.

\begin{definition}[Spatial Selection Operator] Define the spatial selection
operator ${\rm S}_R$, which selects the function in a volume
region $R\subset \bll$, with kernel $S_R(\bv{r}, \bv{r'})$ as
\begin{align}\label{Eq:kernel_spatial_def}
S_R(\bv{r}, \bv{r'}) \,\triangleq \, I_R(\bv{r}) \delta(\bv{r},
\bv{r'}),
\end{align}
where $I_R(\bv{r})=1$ for $\bv{r}\in R$ and $I_R(\bv{r})=0$ for
$\bv{r}\in \bll\backslash R$ is an indicator function of the region
$R$. $\rm{S}_R$ projects the signal $f \in \lball$ onto the subspace $\lballs{R}$.%
\end{definition}

\begin{definition}[Fourier-Bessel Spectral Selection Operator] Define the spectral selection operator ${ \rm  \tilde S}_{KL}$, which selects the harmonic contribution of functions in the Fourier-Bessel spectral domain $\tilde A_{KL}$, by its kernel
\begin{align}\label{Eq:kernel_spectral_def_FB}
\tilde S_{KL}(\bv{r}, \bv{r'}) \,\triangleq \, \sum_{\ell=0}^{L-1}\,\sum_{m=-\ell}^\ell\,\int_{k=0}^K \ddr{k} \fbc{X}{\ell m}{}(k,\bv{r})\, \fbc{X}{\ell m}{\ast}(k,\bv{r'}).\end{align}
The operator ${ \rm \tilde S}_{KL}$ projects a signal onto the subspace of Fourier-Bessel band-limited functions $\lballss{K}{L}$. \end{definition}

\begin{definition}[Fourier-Laguerre Spectral Selection Operator] Define the spectral selection operator ${\rm S}_{PL}$, which selects the harmonic contribution of functions in the Fourier-Laguerre spectral domain
$A_{PL}$, by its kernel
\begin{align}\label{Eq:kernel_spectral_def_FL} S_{PL}(\bv{r}, \bv{r'}) \,\triangleq \, \sum_{p=0}^{P-1} \sum_{\ell=0}^{L-1}\,\sum_{m=-\ell}^\ell\, Z_{\ell m p}(\bv{r})\, Z_{\ell m p}^\ast(\bv{r'}).  \end{align}
The operator ${\rm S}_{PL}$ projects a signal onto the subspace of Fourier-Laguerre band-limited functions $\lballs{PL}$.
\end{definition}

Since both the spatial and spectral selection operators are projection operators, they are idempotent and self-adjoint in nature. By noting the expansions of the Dirac delta in \eqref{Eq:dirac_expansion_FB} and \eqref{Eq:dirac_expansion_FL}, it is evident that the kernels of the spectral selection operators in \eqref{Eq:kernel_spectral_def_FB} and \eqref{Eq:kernel_spectral_def_FL} are Dirac delta functions that are band-limited in the appropriate basis.

\section{Simultaneous Concentration in Spatial and Spectral Domain}
\label{sec:conc_problem}
By virtue of the uncertainty principle, no function can be
space-limited and band-limited simultaneously. In other words, a
signal $f\in \lball$ cannot belong to the subspace $\lballs{R}$ and
the subspaces $\lballs{PL}$ or $\lballss{K}{L}$ at the same time. In
this section we first develop a framework to determine band-limited
functions $f\in \lballs{PL}$ or $f\in\lballss{K}{L}$ that are
optimally concentrated in some spatial region $R\subset\bll$.  We then
formulate the problem to determine space-limited functions
$g\in\lballs{R}$ that are optimally concentrated within the spectral
region $\tilde A_{KL}$ or $A_{PL}$. Later, we show the equivalence
between these two problems and provide the harmonic domain
formulations of the problems.

\subsection{Spatial Concentration of Band-Limited Functions}

Let $f\in \lball$ be a band-limited signal, that is, $f\in\lballp$ or
$f\in\lballk$.  The energy concentration of the function $f$ within the spatial region $R\subset\bll$ is given by

\begin{align}\label{Eq:eigen_problem_spectral}
\lambda &= \frac{\innerp[\big]{{\rm S}_R  f}{{ \rm S}_R  f}}{\innerp[\big]{ f}{f}} \nonumber \\
        &= \frac{\mathlarger{\int}_R \ddv{\bv{r}} \mathlarger{\int}_R  \ddv{\bv{r'}} f(\bv{r}) f^{\ast}(\bv{r'})\delta(\bv{r},\bv{r'}) }{\mathlarger{\int}_{\bll} \ddv{\bv{r}} f(\bv{r}) f^\ast(\bv{r})
        } \nonumber \\
        &= \frac{\mathlarger{\int}_R  \ddv{\bv{r}} |f(\bv{r})|^2  } {\mathlarger{\int}_{\bll} \ddv{\bv{r}}
        |f(\bv{r})|^2 }.
\end{align}
It is well-known that the function $f$ that renders the Rayleigh quotient \eqref{Eq:eigen_problem_spectral} stationary is a solution of
following Fredholm integral equation~\cite[Ch. 8]{Kennedy-book:2013}:
\begin{align}\label{Eq:eigen_problem_spectral2}
\int_R \ddv{\bv{r'}} \delta(\bv{r},\bv{r'}) f(\bv{r'})  = \lambda
f(\bv{r}),
\end{align}
or equivalently (using the integral representation of the spatial
selection operator),
\begin{align}\label{Eq:eigen_problem_spectral3}
{\rm S}_{R}  f = \lambda  f,
\end{align}
where we again note that $f$ is a band-limited function. The
solution of the eigenvalue problem in \eqref{Eq:eigen_problem_spectral3}
yields band-limited eigenfunctions. Since ${\rm S}_R$ is a
projection operator, these eigenvalues are positive and bounded
above by unity. The eigenvalue associated with each eigenfunction
serves as a measure of energy concentration of the eigenfunction in the
spatial region $R$. We discuss the properties of the eigenfunctions
in \secref{sec:eignfunctions_analysis}.

\subsection{Spectral Concentration of Space-Limited Functions}

Here we consider the dual of the problem posed in the previous
subsection. Instead of seeking band-limited spatially concentrated
functions, we seek space-limited functions with optimal concentration
in some spectral region. Let $g\in\lballr$ be the space-limited
function within the spatial region $R$. We maximize the concentration
of $g$ in the spectral region $\tilde A_{KL}$ or $A_{PL}$, depending
on the basis functions chosen for the characterization of the spectral
domain (Fourier-Bessel or Fourier-Laguerre, respectively). To maximize
the concentration of the space-limited signal $g$ within the spectral
region $\tilde A_{KL}$ in the Fourier-Bessel spectral domain, we
maximize the ratio
\begin{align}\label{Eq:eigen_problem_spatialFB}
\lambda \,&=\, \frac{\innerp[\big]{{\rm \tilde  S}_{KL}
g}{{ \rm \tilde S}_{KL} g}}{\innerp[\big]{g}{g}},\quad g\in\lballs{R}.
\end{align}
Similarly, to maximize the concentration within
the spectral region $A_{PL}$ in the Fourier-Laguerre spectral domain, we maximize the ratio
\begin{align}\label{Eq:eigen_problem_spatialFL}
\lambda \,&=\, \frac{\innerp[\big]{{\rm S}_{PL} g}{{\rm S}_{PL} g}}{\innerp[\big]{g}{g}},\quad g\in\lballs{R}.
\end{align}
Following a similar approach to the spatial concentration problem
above and using the integral representations of the spectral
selection operators ${ \rm \tilde S}_{KL}$ and ${\rm S}_{PL}$,
the Fourier-Bessel concentration problem in \eqref{Eq:eigen_problem_spatialFB}
results in the following eigenvalue problem
\begin{align}\label{Eq:eigen_problem_spatialFB2}
\sum_{\ell=0}^{L-1}\,\sum_{m=-\ell}^\ell\,  \int_R \ddv{\bv{r'}}
\int_{k=0}^K \ddr{k} \fbc{X}{\ell m}{}(k,\bv{r})\, \fbc{X}{\ell
m}{\ast}(k,\bv{r'}) g(\bv{r'})    = \lambda g(\bv{r}),
\end{align}
or equivalently
\begin{align}\label{Eq:eigen_problem_spatialFB3}
{ \rm \tilde S}_{KL} g = \lambda g.
\end{align}
Similarly, the Fourier-Laguerre problem in \eqref{Eq:eigen_problem_spatialFL} gives rise to following eigenvalue problem
\begin{align}\label{Eq:eigen_problem_spatialFL2}
  \sum_{p=0}^{P-1}
\sum_{\ell=0}^{L-1}\,\sum_{m=-\ell}^\ell\, \mathlarger{\int}_R
\ddv{\bv{r'}}\, Z_{\ell m p}(\bv{r})\, Z^\ast_{\ell m p}(\bv{r'})
g(\bv{r'})    = \lambda g(\bv{r}),
\end{align}
or equivalently
\begin{align}\label{Eq:eigen_problem_spatialFL3}
{\rm S}_{PL} g = \lambda g.
\end{align}
The solution of each of the eigenvalue problems presented in
\eqref{Eq:eigen_problem_spatialFB3} and
\eqref{Eq:eigen_problem_spatialFL3} provides space-limited
eigenfunctions. Again, the eigenvalues are positive and bounded
above by unity since the both ${\rm S}_{PL}$ and
${ \rm \tilde S}_{KL}$ are projection operators. The eigenvalue
associated with each eigenfunction serves as a measure of
concentration of the eigenfunction in the spectral region
$A_{PL}$ or $\tilde A_{KL}$, as we show in \secref{sec:eignfunctions_analysis}.

So far, we have formulated two types of eigenvalue problems: (1) problem
\eqref{Eq:eigen_problem_spectral3} to find band-limited
spatially concentrated functions $f$; and (2) problems
\eqref{Eq:eigen_problem_spatialFB3} and \eqref{Eq:eigen_problem_spatialFL3} to find space-limited
functions $g$ with optimal concentration within a spectral region.
We go on to show that the eigenfunctions that arise as a
solution of both problems are the same, up to a multiplicative
constant, within the spatial region $R$ and spectral region
$A_{PL}$ or $\tilde A_{KL}$.

\subsection{Equivalence of Problems}
Here we show that the concentration problems to find the spatially
concentrated band-limited functions and spectrally concentrated
space-limited functions have equivalent solutions. Equivalence is shown explicitly only for the Fourier-Bessel spectral
domain; however, the same result also holds for the Fourier-Laguerre
spectral domain.

If the signal $f$ in
\eqref{Eq:eigen_problem_spectral3} is band-limited within the
Fourier-Bessel spectral domain, that is, $f = { \rm \tilde S}_{KL} f$, then we
can write the eigenvalue problem as
\begin{align}\label{Eq:eigen_problem_spectral4}
{\rm U}f = \lambda  f,
\end{align}
where ${\rm U} = {\rm S}_{R} \tilde{\rm S}_{KL}$ denotes the
composite operator with kernel given by
\begin{align}\label{Eq:kernel_comp1}
U(\bv{r},\bv{r'}) &= I_R(\bv{r})
\sum_{\ell=0}^{L-1}\,\sum_{m=-\ell}^\ell\, \int_{k=0}^K \ddr{k}
\fbc{X}{\ell m}{}(k,\bv{r})\,
\fbc{X}{\ell m}{\ast}(k,\bv{r'})  \nonumber \\
&= I_R(\bv{r}) \tilde S_{KL}(\bv{r},\bv{r'}).
\end{align}
Similarly, if $g$ is a space-limited signal, that is $g = {\rm S}_{R} g$, then the eigenvalue problem in \eqref{Eq:eigen_problem_spatialFB3} can be expressed as
\begin{align}\label{Eq:eigen_problem_spatialFB4}
  {\rm V} g = \lambda g,
\end{align}
where ${\rm V} = { \rm \tilde S}_{KL} {\rm S}_{R} $,
with kernel
\begin{align}\label{Eq:kernel_comp2}
V(\bv{r},\bv{r'}) &= I_R(\bv{r'})
\sum_{\ell=0}^{L-1}\,\sum_{m=-\ell}^\ell\, \int_{k=0}^K \ddr{k}
\fbc{X}{\ell m}{}(k,\bv{r})\,
\fbc{X}{\ell m}{\ast}(k,\bv{r'})  \nonumber \\
&= I_R(\bv{r'}) \tilde S_{KL}(\bv{r},\bv{r'}).
\end{align}
The composite operators ${\rm U}$ and ${\rm V}$ are not
commutative in general, that is, $U(\bv{r},\bv{r'})\neq
V(\bv{r},\bv{r'})$. However, since $I_R(\bv{r'})=1$ for $\bv{r'}\in R$, the
action of these composite operators is commutative in the spatial
region $R$, that is,
\begin{align}
U(\bv{r},\bv{r'}) = V(\bv{r},\bv{r'}),\quad \bv{r},\,\bv{r'}\in
R\subset\bll.
\end{align}
Consequently, the solution of the eigenvalue problems in
\eqref{Eq:eigen_problem_spectral4} and
\eqref{Eq:eigen_problem_spatialFB4} have the same solution and the same eigenvalue within the
region $R$, that is,
\begin{align}\label{Eq:link_f_g}
g(\bv{r}) = \left({\rm S}_{R}f\right)(\bv{r}),\quad \bv{r}\in
R\subset\bll.
\end{align}
For a signal band-limited in the Fourier-Laguerre spectral domain, the
equivalent of the eigenvalue problem in
\eqref{Eq:eigen_problem_spectral4} is
\begin{align}\label{Eq:eigen_problem_spectral5}
{\rm W}f = \lambda  f,
\end{align}
where ${\rm W} = {\rm S}_{R} {\rm S}_{PL} $, and the
analogous result holds.

\begin{remark}\label{remark:equivalence_problems}
The equivalence of the spatial and spectral concentration problems implies that we only need to solve
the eigenvalue problem presented in
\eqref{Eq:eigen_problem_spectral4} or
\eqref{Eq:eigen_problem_spectral5} to obtain the band-limited spatially
concentrated eigenfunctions and the space-limited spectrally
concentrated eigenfunctions then can be obtained using
\eqref{Eq:link_f_g}, i.e.,\ by setting the space-limited eigenfunctions to the band-limited eigenfunctions in the region $R$ and zero elsewhere.
\end{remark}

\subsection{Harmonic Domain Analysis}
So far the eigenvalue problems have been formulated in the spatial
domain. Here we present the spectral domain formulation of
the eigenvalue problems presented in \eqref{Eq:eigen_problem_spectral4}
and \eqref{Eq:eigen_problem_spectral5}. Using the kernel
representation of the composite operator ${\rm U}$ given in
\eqref{Eq:kernel_comp1}, we can write \eqref{Eq:eigen_problem_spectral4} as
\begin{align}\label{Eq:eig_harmonic_analysis1}
I_{R}(\bv{r}) \intball
 \ddv{\bv{r'}}  \int_{k'=0}^{K} \ddr{k'}\,
\sum_{\ell'=0}^{L-1}\,\sum_{m'=-\ell'}^{\ell'}\, \fbc{X}{\ell'
m'}{}(k',\bv{r})\fbc{X}{\ell' m'}{\ast}(k',\bv{r'})
f(\bv{r'})  &= \lambda
f(\bv{r}).
\end{align}
By taking the Fourier-Bessel transform of \eqref{Eq:eig_harmonic_analysis1} with respect to spatial variable $\bv{r}$, we obtain the following
formulation of the eigenvalue problem in the Fourier-Bessel spectral domain
\begin{align}\label{Eq:eig_harmonic_spectral1}
 \sum_{\ell'=0}^{L-1}\,\sum_{m'=-\ell}^\ell\,\int_{k'=0}^{K} \ddr{k'}
\cpfb_{\ell' m', \ell m}{(k',k)}\fbc{f}{\ell' m'}{}(k')  =
 \lambda \fbc{f}{\ell m}{}(k),
\end{align}
with
\begin{align}\label{Eq:eig_harmonic_spectral1_matrix}
\cpfb_{\ell m, \ell' m'}{(k,k')} = \int_{R} \ddv{\bv{r}}
\fbc{X}{\ell m}{}(k,\bv{r}) \fbc{X}{\ell' m'}{\ast}(k',\bv{r}).
\end{align}
Similarly, the eigenvalue problem \eqref{Eq:eigen_problem_spectral5} can be formulated in the Fourier-Laguerre spectral domain as
\begin{align}\label{Eq:eig_harmonic_spectral2}
\sum_{p'=0}^{P} \sum_{\ell'=0}^{L-1}\,\sum_{m'=-\ell'}^{\ell'} \,
\cpfl_{ \ell' m' p', \ell m p} \shc{f}{\ell' m' p'}  = \lambda
\shc{f}{\ell m p },
\end{align}
with
\begin{align}\label{Eq:eig_harmonic_spectral2_matrix}
\cpfl_{\ell m p, \ell' m' p'} = \int_{R} \ddv{\bv{r}} Z_{\ell m
p}(\bv{r}) Z_{\ell' m' p'}^{\ast}(\bv{r}).
\end{align}
By defining the matrix $\mathbf{\cpfl}$ of size $PL^2\times PL^2$
with entries given by \eqref{Eq:eig_harmonic_spectral2_matrix} and the vector $\mathbf{f} =
\left(\shc{f}{000},\,\shc{f}{001},\,\hdots,
\,\shc{f}{(L-1)(L-1)(P-1)}\right)^{\rm T}$ of length $PL^2$ as the spectral representation of the signal $f$,
\eqref{Eq:eig_harmonic_spectral2} can be written compactly in matrix
form as
\begin{align}\label{Eq:eig_harmonic_spectral2_matrix_formulation}
\mathbf{\cpfl}\mathbf{f} = \lambda \mathbf{f}.
\end{align}
Thus, the spectral representation $\mathbf{f}$ of the band-limited spatially concentrated signal $f$ can be obtained as a solution of an algebraic~(finite dimensional) eigenvalue problem of size \mbox{$PL^2\times PL^2$}.
Due to the continuous nature of the Fourier-Bessel harmonic space (i.e.\ $k$ is continuous), an equivalent finite-dimensional matrix formulation cannot be written for the Fourier-Bessel setting.

For the eigenvalue problems \eqref{Eq:eigen_problem_spectral4} and \eqref{Eq:eigen_problem_spectral5}, which are expressed in terms of the spatial domain representation of the signal, we have obtained here the equivalent harmonic formulations \eqref{Eq:eig_harmonic_spectral1} and \eqref{Eq:eig_harmonic_spectral2_matrix_formulation} respectively.
In the next section, we discuss the properties of band-limited and
space-limited eigenfunctions.

\section{Analysis of Eigenfunctions and Eigenvalue
Spectrum}\label{sec:eignfunctions_analysis}

We first study the properties of the both band-limited and space-limited eigenfunctions in this section, both for the Fourier-Bessel and Fourier-Laguerre scenarios.  In both of these scenarios we also study the eigenvalue spectrum and calculate the analog of the Shannon number in the one dimensional Slepian concentration problem.

\subsection{Properties of Fourier-Bessel Band-Limited Eigenfunctions}

\hfill The \newline \mbox{Fourier-Bessel} band-limited and spatially
concentrated eigenfunctions are recovered in the spectral domain by
solving the eigenvalue problem given in
\eqref{Eq:eig_harmonic_spectral1}.  That is, we obtain $\fbc{f}{\ell
  m}{}(k)$ for $0 \leq k \leq K$, $0\leq \ell \leq L-1$ and $|m| \leq
\ell$. In practice, the spectrum along $k\in\rplus$ is discretized to
solve the eigenvalue problem in \eqref{Eq:eig_harmonic_spectral1} as
an algebraic eigenvalue problem (we further elaborate this in Section~\ref{sec:special_regions}).
Since both \eqref{Eq:eigen_problem_spectral4} and
\eqref{Eq:eig_harmonic_spectral1} are equivalent and the operator
${\rm U}$ in \eqref{Eq:eigen_problem_spectral4} is a composite
projection operator, the eigenfunctions are
orthogonal  and the associated eigenvalue of each eigenfunction is
real, positive and bounded above by unity. We choose the
eigenfunctions to be orthonormal. Since the
spectral response is continuous along $k$, the number of
eigenfunctions is (theoretically) infinite and depends on the resolution of the discretization of the spectrum along $k$. We order eigenfunctions $f^1,\,f^2,\,\hdots$ and eigenvalues
$\lambda_1,\,\lambda_2,\,\hdots $ such that
$1\leq\lambda_1\leq\lambda_2 \leq \hdots \leq 0$.

The eigenfunctions $f^\alpha \in \lball $ in the spatial domain can be
recovered from their spectral representation through the inverse
Fourier-Bessel transform \eqref{Eq:FB_inverse_transform}, where
$\alpha$ is used to index the eigenfunctions.  Since the Hermitian
symmetry property $ \cpfb_{\ell m, \ell' m'}{(k,k')} =
\left(\cpfb_{\ell' m', \ell m}{(k',k)}\right)^\ast$ is satisfied, as
is directly apparent from \eqref{Eq:eig_harmonic_spectral1_matrix}, it
follows that the eigenvalues are real and the eigenfunctions
orthogonal.  The band-limited eigenfunctions are orthonormal in both
the Fourier-Bessel spectral domain and in the entire spatial domain $\bll$, that is,
\begin{align}\label{Eq:properties_eig_FB_spectral_first}
\sum_{\ell=0}^{L-1} \sum_{m=-\ell}^\ell\int_\rplus \ddr{k} \fbc{f}{\ell
m}{\alpha}(k)\left( \fbc{f}{\ell m}{\beta}(k)\right)^\ast &= \delta_{\alpha\beta},
\\
\label{Eq:properties_eig_FB_spectral_second}
\int_\bll \ddv{\bv{r}} f^\alpha(\bv{r})  \conj{f^\beta(\bv{r})} &=
\delta_{\alpha\beta}.
\end{align}
The eigenfunctions are also orthogonal (but not orthonormal) within the spatial region
$R$, with
\begin{align}\label{Eq:properties_eig_FB_spatial}
\int_R \ddv{\bv{r}} f^\alpha(\bv{r}) \conj{f^\beta(\bv{r})} =
\lambda_\alpha\delta_{\alpha\beta},
\end{align}
which is obtained using the spectral domain formulation of the eigenvalue problem in \eqref{Eq:eig_harmonic_spectral1} and the orthonormality relation in \eqref{Eq:properties_eig_FB_spectral_first}.  It is clear from \eqref{Eq:properties_eig_FB_spatial} that the eigenvalue $\lambda$
associated with each unit energy band-limited eigenfunction provides a measure of energy concentration within the spatial region $R$.

\subsection{Properties of Fourier-Laguerre Band-Limited Eigenfunctions}

The Fourier-Laguerre band-limited and spatially concentrated eigenfunctions are recovered in the spectral domain by solving the eigenvalue problem given in \eqref{Eq:eig_harmonic_spectral2_matrix_formulation}, to recover
the eigenvector $\mathbf{f}$.  Since the matrix
$\mathbf{\cpfl }$ in
\eqref{Eq:eig_harmonic_spectral2_matrix_formulation} is of size
$PL^2 \times PL^2$, the number of eigenvectors/eigenfunctions is also $PL^2$.
The matrix $\mathbf{\cpfl}$ is positive definite, therefore the
eigenvalues are positive, real and bounded by unity and the eigenvectors are orthogonal.  We chose the eigenvectors to be orthonormal. We again order
 eigenvectors $\mathbf{f}^1,\,\mathbf{f}^2,\,\hdots,\,\mathbf{f}^{PL^2}$ and eigenvalues $\lambda_1,\,\lambda_2,\,\hdots,\,\lambda_{PL^2}$ such that
$1\leq\lambda_1\leq\lambda_2 \leq \hdots \leq \lambda_{PL^2} \leq 0$.

The eigenfunctions $f^\alpha(\bv{r})$ in the spatial domain can be recovered from their spectral representation though the inverse Fourier-Laguerre transform \eqref{Eq:f_expansion}, where again $\alpha$ is used to index the eigenfunctions.  Again, the reality of eigenvalues and orthogonality of eigenvectors also follows from the fact that the matrix $\mathbf{\cpfl}$ in \eqref{Eq:eig_harmonic_spectral2_matrix},
with entries given in
\eqref{Eq:eig_harmonic_spectral2_matrix_formulation}, is Hermitian
symmetric. The band-limited eigenfunctions are orthonormal in both the
Fourier-Laguerre spectral domain and in the entire spatial domain $\bll$, that is,
\begin{align}\label{Eq:properties_eig_FL_spectral_first}
(\mathbf{f}^\alpha)^{\rm H} \mathbf{f}^\beta
=\sum_{p=0}^{P-1}\sum_{\ell=0}^{L-1} \sum_{m=-\ell}^\ell &\left(f_{\ell m
p}^\alpha\right)^\ast f_{\ell m p}^\beta = \delta_{\alpha\beta},
\\
\label{Eq:properties_eig_FL_spectral_second}
\int_\bll \ddv{\bv{r}} f^\alpha(\bv{r}) \conj{f^\beta(\bv{r})} &=
\delta_{\alpha\beta},
\end{align}
where $(\cdot)^{\rm H}$ denotes the Hermitian transpose operation, and are orthogonal (but not orthonormal) within the spatial region $R$,
\begin{align}\label{Eq:properties_eig_FL_spatial}
\int_R \ddv{\bv{r}} f^\alpha(\bv{r})  \conj{f^\beta(\bv{r})} =
\lambda_\alpha\delta_{\alpha\beta},
\end{align}
which is obtained by using \eqref{Eq:eig_harmonic_spectral2} and the orthonormality relation in \eqref{Eq:properties_eig_FL_spectral_first}.  It is clear from \eqref{Eq:properties_eig_FL_spatial} that the eigenvalue $\lambda_\alpha$ associated with each unit energy band-limited eigenfunction $f^\alpha$ provides a measure of energy concentration within the spatial region $R$.

\subsection{Properties of Space-Limited Eigenfunctions} As highlighted earlier, the space-limited spectrally concentrated eigenfunctions $g(\bv{r})$ of the eigenvalue problems \eqref{Eq:eigen_problem_spatialFB3} and \eqref{Eq:eigen_problem_spatialFL3} can be obtained from the band-limited spatially concentrated eigenfunctions of problems \eqref{Eq:eig_harmonic_spectral1} and \eqref{Eq:eig_harmonic_spectral2_matrix_formulation}, respectively.  This is achieved by employing the relationship
given in \eqref{Eq:link_f_g}, i.e., by setting the space-limited eigenfunctions to the band-limited eigenfunctions in the region $R$ and zero elsewhere.

Each band-limited eigenfunction $f^\alpha(\bv{r})$ thus provides a space-limited eigenfunction $g^\alpha(\bv{r})$ and the eigenvalue associated with each eigenfunction provides a measure of energy concentration within the Fourier-Bessel or Fourier-Laguerre spectral domain.  In order to normalize the space-limited eigenfunctions to unit energy, we update the relationship between $f^\alpha(\bv{r})$ and $g^\alpha(\bv{r})$, given in \eqref{Eq:link_f_g}, by noting \eqref{Eq:properties_eig_FB_spatial} and \eqref{Eq:properties_eig_FL_spatial}, as
\begin{align}\label{Eq:ftog_in_spatial}
g^\alpha(\bv{r}) = \frac{1}{\sqrt{\lambda_\alpha}}
\left({\rm S}_{R} f\right)(\bv{r}), \quad \lambda_\alpha \neq 0.
\end{align}
The revised space-limited eigenfunctions
$g^\alpha(\bv{r})$ then satisfy
\begin{align}
\int_\bll \ddv{\bv{r}} g^\alpha(\bv{r}) g^\beta(\bv{r}) &= \int_R
\ddv{\bv{r}} g^\alpha(\bv{r})  g^\beta(\bv{r}) =
\delta_{\alpha\beta}.
\end{align}

The relation between the band-limited eigenfunctions
$f^\alpha(\bv{r})$ and space-limited eigenfunctions
$g^\alpha(\bv{r})$ can be expressed in the spectral domain as
\begin{align}\label{Eq:relation_g_f_spectral1}
\fbc{g}{\ell m}{\alpha}(k) &= \frac{1}{\sqrt{\lambda_\alpha}}
\sum_{\ell'=0}^{L-1}\,\sum_{m'=-\ell}^\ell\,\int_{k'=0}^{K} \ddr{k'}
\cpfb_{\ell' m', \ell m}{(k',k)}\fbc{f}{\ell' m'}{\alpha}(k') ,\quad
f^\alpha(\bv{r}) \in \lballss{K}{L},
\nonumber \\
g_{\ell m p}^\alpha &= \frac{1}{\sqrt{\lambda_\alpha}}
\sum_{p'=0}^{P} \sum_{\ell'=0}^{L-1}\,\sum_{m'=-\ell'}^{\ell'} \,
\cpfl_{\ell' m' p', \ell m p} \shc{f}{\ell' m' p'}^\alpha ,\quad
f^\alpha(\bv{r}) \in \lballs{PL},
\end{align}
for $\ell,\,p \in \mathbb{Z}^+$ and $k \in \rplus$,
where we have used the kernel representation of the spatial-selection
operator ${\rm S}_R$. Note that the harmonic representations of the space-limited eigenfunctions
$\fbc{g}{\ell m}{\alpha}(k)$ and $g_{\ell m p}^\alpha$ are obviously not band-limited and the expressions given above apply over the entire harmonic domains.
Using \eqref{Eq:eig_harmonic_spectral1}
and \eqref{Eq:eig_harmonic_spectral2} these relations can also be
expressed within the finite spectral regions as
\begin{align}\label{Eq:relation_g_f_spectral2}
\fbc{g}{\ell m}{\alpha}(k) &= \sqrt{\lambda_\alpha} \fbc{f}{\ell
m}{\alpha}(k) ,\quad  \ell,\,k\in \tilde A_{KL},
\nonumber \\
g_{\ell m p}^\alpha &= \sqrt{\lambda_\alpha} f_{\ell m p}^\alpha
,\quad \, \ell,\,p \in A_{PL}.
\end{align}
Finally, using \eqref{Eq:properties_eig_FB_spectral_first} or \eqref{Eq:properties_eig_FL_spectral_first} in conjunction with \eqref{Eq:relation_g_f_spectral2},
we can obtain the orthogonality relations for space-limited eigenfunctions in the spectral domain
\begin{align}
\sum_{\ell=0}^{L-1} \sum_{m=-\ell}^\ell \int_{k=0}^K \ddr{k} \,\fbc{g}{\ell m}{\alpha}(k) \left(\fbc{g}{\ell m}{\beta}(k)\right)^\ast &= \lambda_\alpha \delta_{\alpha\beta},
\nonumber
\\
\sum_{p=0}^{P-1}\sum_{\ell=0}^{L-1} \sum_{m=-\ell}^\ell  g_{\ell m p}^\alpha \left(g_{\ell m p}^\beta\right)^\ast &= \lambda_\alpha \delta_{\alpha\beta}.
\label{Eq:properties_eig_spectral}
\end{align}
It is clear from \eqref{Eq:properties_eig_spectral} that the eigenvalue $\lambda$ associated with each unit energy space-limited eigenfunction provides a measure of energy concentration within the spectral regions $\tilde A_{KL}$ or $A_{PL}$.

\begin{remark}\label{remark:regarding_leakage}
We highlight again that the eigenvalue $\lambda_\alpha$ serves as a
measure of the concentration of the (unit energy) band-limited eigenfunctions
$f^\alpha(\bv{r})$ in the spatial domain, and the concentration of
the (unit energy) space-limited eigenfunctions $g^\alpha(\bv{r})$ in the spectral
domain. The energy $1-\lambda_\alpha$ leaked by $f^\alpha(\bv{r})$
into the spatial region $\bll\backslash R$ is equal to the energy
leaked by $g^\alpha(\bv{r})$ outside of the spectral region
$\tilde A_{KL}$ or $A_{PL}$, depending on the chosen spectral domain.
\end{remark}

\begin{remark} Since $\lballs{R}$ is an infinite dimensional subspace of $\lball$, the number of space-limited eigenfunctions arising from the solution of $\eqref{Eq:eigen_problem_spatialFB4}$ is infinite. However, when we seek the spectral concentration in the Fourier-Laguerre spectral domain, the equivalent problem in the spectral domain given in \eqref{Eq:eig_harmonic_spectral2_matrix_formulation} is finite dimensional, giving rise to a finite number ($PL^2$) of band-limited, and thus also space-limited, eigenfunctions. The remaining space-limited eigenfunctions, which do not have any energy in the spatial region $R$ nor in the spectral region $A_{PL}$, belong to the infinite dimensional null space of ${\rm V}$ and have associated eigenvalue equal to zero. \end{remark}

\subsection{Eigenvalue Spectrum}
Following Remark \ref{remark:regarding_leakage}, the band-limited
eigenfunctions and space-limited eigenfunctions, which are well
concentrated within the spatial region and spectral region
respectively, have eigenvalue near unity, whereas those which are
poorly concentrated have value near zero. If the spectrum of
eigenvalues $\lambda_1,\,\lambda_2,\,\hdots,\,$ has a narrow
transition width from values near zero to values near unity, as in
the case of one dimensional Slepian concentration problem~\cite{Slepian:1965}
and spherical concentration problem~\cite{Simons:2006}, the sum of
all of the eigenvalues well approximates the number of significant
eigenfunctions with eigenvalue near unity. We elaborate this fact
later in the paper with the help of examples. If the signal is
band-limited in the Fourier-Bessel spectral domain, the sum of
eigenvalues is given by
\begin{align}\label{Eq:sum_eigen_values_FB}
\tilde N_{KL} = \textrm{tr}({\rm U}) &= \int_{\bll} \ddv{\bv{r}} U(\bv{r},\bv{r}) = \int_{k=0}^{K} \ddr{k} \sum_{\ell=0}^{L-1}\,\sum_{m=-\ell}^\ell \cpfb_{\ell m, \ell m}{(k,k)}  \nonumber \\
&=  \sum_{\ell=0}^{L-1}\frac{2\ell+1}{2\pi^2} \int_R \ddv{\bv{r}} \left(\int_{k=0}^{K}\ddr{k} k^2 j_\ell(k r) j_\ell(k r) \right) \nonumber \\
&=  \sum_{\ell=0}^{L-1}\frac{2\ell+1}{4\pi^2}K^3 \int_R
\ddv{\bv{r}} \left( \left(j_\ell(Kr)\right)^2 -
j_{\ell-1}(Kr)j_{\ell+1}(Kr) \right),
\end{align}
with $\bv{r} =
\bv{r}(r,\theta,\phi)$. Here we use $\textrm{tr}(\cdot)$ to denote the trace of an operator or matrix.

For a signal band-limited in the Fourier-Laguerre spectral domain, the
sum of eigenvalues is given by
\begin{align}\label{Eq:sum_eigen_values_FL}
N_{PL} &= \textrm{tr}({\rm W}) =
\textrm{tr}(\mathbf{\cpfl})
= \int_{\bll} \ddv{\bv{r}} W(\bv{r},\bv{r}) \nonumber \\ &=  \sum_{p=0}^{P-1} \sum_{\ell=0}^{L-1}\,\sum_{m=-\ell}^\ell \cpfl_{\ell m p, \ell m p} \nonumber \\
&= \frac{L^2}{4\pi}  \sum_{p=0}^{P-1} \int_R \ddv{\bv{r}} K_p(r)
K_p(r) , \quad \bv{r} = \bv{r}(r,\theta,\phi).
\end{align}
We note that the sum of eigenvalues given in
\eqref{Eq:sum_eigen_values_FB} and \eqref{Eq:sum_eigen_values_FL}
involve an integral of a positive function over the spatial region
$R$. Thus, if the volume of the spatial region is relatively
smaller, the sum of eigenvalues, which also indicates the number of
significant eigenfunctions, is also smaller.

\begin{remark}
  $\tilde N_{KL}$ and $N_{PL}$ serve as an analog of the Shannon
  number in the one-dimensional Euclidean Slepian concentration
  problem~\cite{Slepian:1965,Percival:1993}, in which the Shannon
  number is given by a space-bandwidth product. Here, we refer to
  $N_{KL}$ or $N_{PL}$ as spherical Shannon numbers. Since the
  significant eigenfunctions have value near unity and the remaining
  eigenfunctions have eigenvalue near zero, the number of significant
  eigenfunctions is given approximately by $\tilde N_{KL}$ or
  $N_{PL}$. Thus, a function spatially concentrated in $R$ and
  band-limited in the Fourier-Bessel (Fourier-Laguerre) domain within
  the spectral region $\tilde A_{KL}$ ($A_{PL}$), can be represented
  approximately by $\tilde N_{KL}$ ($N_{PL}$) band-limited orthonormal
  eigenfunctions. This is the essence of the concentration problem:
  spatially concentrated band-limited functions on the ball, which
  belong to the infinite dimensional space $\lballss{K}{L}$, can be
  represented by $\tilde N_{KL}$ finite orthonormal basis functions of
  the same space. Similarly, the dimension of the space required to
  represent spatially concentrated band-limited signals belonging to
  the finite dimensional space $\lballss{K}{L}$ of size $PL^2$ can be
  reduced from $PL^2$ to $N_{PL}$. This fact is further explored in
  \secref{sec:Applications}.  \end{remark}

\section{Eigenfunctions Concentrated within Spatially Symmetric Regions}
\label{sec:special_regions}

In the preceding section we analysed the properties of eigenfunctions arising from the spatial-spectral concentration problem on the ball.  We have yet to discuss the computation of these eigenfunctions. As mentioned briefly earlier, the eigenvalue problem of \eqref{Eq:eig_harmonic_spectral1} to find eigenfunctions band-limited in the Fourier-Bessel domain can be transformed into an algebraic eigenvalue problem by discretizing the spectrum along $k\in\rplus$. If we discretize along $k$ by selecting $M$ uniform samples $k_n = nK/M,\,n=1,\,2,\,\hdots,\,M$, between $0 <k \leq K$, the eigenvalue problem in \eqref{Eq:eig_harmonic_spectral1} becomes an algebraic eigenvalue problem of size $ML^2$, for which we need to compute a Hermitian matrix of size $ML^2\times ML^2$.  For each element of this matrix we must evaluate the integral given in \eqref{Eq:eig_harmonic_spectral1_matrix} over the spatial region of interest. Similarly, the eigenvalue problem of \eqref{Eq:eig_harmonic_spectral2}, of size $PL^2\times PL^2$, must be solved to find band-limited eigenfunctions in the Fourier-Laguerre domain.  Again, for each element of the matrix we must evaluate the integral given in \eqref{Eq:eig_harmonic_spectral2_matrix} over the spatial region of interest.  Since observations in practical applications can support very high band-limits, the direct computation of eigenfunctions by solving the algebraic eigenvalue problems is computationally intensive, motivating alternative procedures.

In this section we analyze the eigenvalue problems when the spatial region of concentration is symmetric in nature. Under certain symmetries we show that the eigenvalue problems given in the spectral domain in \eqref{Eq:eig_harmonic_spectral1} and \eqref{Eq:eig_harmonic_spectral2_matrix_formulation}, originally formulated in \eqref{Eq:eigen_problem_spectral4} and \eqref{Eq:eigen_problem_spectral5}, decompose into subproblems, which reduce the computational burden. In some special cases $\cpfb_{\ell m, \ell' m'}{(k,k')}$ in \eqref{Eq:eig_harmonic_spectral1_matrix} and $\cpfl_{\ell m p, \ell' m' p'}$ in \eqref{Eq:eig_harmonic_spectral2_matrix} can be computed analytically. We consider two types of symmetric regions: (1) circularly symmetric regions only; and (2) circularly symmetric and radially independent regions.

\subsection{Circularly Symmetric Region Only}

First, we analyze the case when the spatial region is circularly
symmetric, that is, rotationally symmetric around some axis defined by its center $\bv{r}_0 \dfn \bv{r}_0(r_0,\theta_0,\phi_0) \in \bll $. For convenience, we consider the region with center $\bv{\eta}_0 \dfn \bv{\eta}_0(r_0,0,0) \in \bll $ on the $z$-axis and refer to such a region as an azimuthally symmetric region. Through a rotation of $\theta_0$ around the $y$-axis, followed by a rotation of $\phi_0$ around the $z$-axis, the azimuthally symmetric region with center at $\bv{\eta}_0$  can be transformed into the circularly symmetric region with center at $\bv{r}_0$.  When considering the azimuthally symmetric region, the azimuthal angle becomes independent
of $r$ and $\theta$ and therefore we can write the integral over
region $R$ as $\int_R\equiv \int_{(r,\theta)}\int_{\phi=0}^{2\pi}$.
Noting this decoupling of the integral and the orthonormality of the
complex exponentials, we can simplify $\cpfb_{\ell m, \ell'
m'}{(k,k')}$ given in \eqref{Eq:eig_harmonic_spectral1_matrix} as
\begin{align}\label{Eq:symmetric1_matrix_simple1}
\cpfb_{\ell m, \ell' m'}{(k,k')} &=\delta_{m m'} 2\pi\int_{(r,\theta)} r^2\sin\theta \ddr{r} \ddr{\theta}\,   \fbc{X}{\ell m}{}(k,\bv{r}) \fbc{X}{\ell' m'}{\ast}(k',\bv{r}) , \quad \bv{r} = \bv{r}(r,\theta,0) \nonumber \\
&= \delta_{mm'}\cpfb_{\ell m, \ell' m}{(k,k')},
\end{align}
which can be used to decompose the spectral domain eigenvalue problem in \eqref{Eq:eig_harmonic_spectral1} into $2L-1$ subproblems,
\begin{align}\label{Eq:symmetric1_sub_problem1}
\int_{k'=0}^{K} \ddr{k'} \sum_{\ell'=0}^{L-1}\,\, \cpfb_{\ell' m, \ell
m}{(k',k)}\fbc{f}{\ell' m}{}(k')  =  \fbc{f}{\ell m}{}(k),
\end{align}
for each $m \in \{-(L-1),\,\hdots,\,(L-1)\}$. Furthermore, it can be
easily shown that {$\cpfb_{\ell m, \ell' m}{(k,k')} = \cpfb_{\ell
(-m), \ell' (-m)}{(k,k')}$} which leaves us with $L$ subproblems for
\linebreak \mbox{$m \in \{0,\,1,\,\hdots,\,(L-1)\}$}. The subproblem in
\eqref{Eq:symmetric1_sub_problem1} gives the spectral domain
representation $\fbc{f}{\ell m}{}(k)$ of eigenfunctions for each $m$,
which can be used to obtain the spatial functions for each $m$ as
\begin{align}\label{Eq:spatial_eigen_subproblem1}
f^{(m)}(r,\theta) = \int_{k=0}^K \ddr{k} \sum_{\ell=m}^{L-1} \shc{f}{\ell m
}(k) \fbc{X}{\ell m}{}(k,\bv{r}), \quad \bv{r} = \bv{r}(r,\theta,0).
\end{align}
Furthermore, since $\cpfb_{\ell' m, \ell m}{(k',k)} = \cpfb_{\ell m, \ell' m}{(k,k')}$ is real valued, the spectral domain
representation $\fbc{f}{\ell m}{}(k)$ of the eigenfunction and consequently the spatial eigenfunction $f^{(m)}(r,\theta)$ is real valued. We note that the superscript $m$ on the left
hand side of \eqref{Eq:spatial_eigen_subproblem1} indicates angular
order and should not be confused with the rank of the eigenfunction.

Similarly, the symmetry of the region along the $z$-axis also simplifies
the elements $\cpfl_{\ell m p, \ell' m' p'}$ in
\eqref{Eq:eig_harmonic_spectral2_matrix} as
\begin{align}
\cpfl_{\ell m p, \ell' m' p'} = \delta_{mm'}\cpfl_{\ell m p, \ell' m p'} = \delta_{(-m)(-m')}\cpfl_{\ell (-m) p, \ell' (-m) p'},
\end{align}
due to which the matrix $\mathbf{\cpfl}$ becomes block diagonal. By
defining the matrix $\mathbf{\cpfl}^{(m)}$ of size $P(L-m)\times
P(L-m)$, with entries
\begin{align}\label{Eq:symmetric1_matrix_new}
\cpfl^{(m)}_{\ell p,\ell'p' } = \cpfl_{\ell m p, \ell' m p'},
\end{align}
the matrix eigenvalue problem in \eqref{Eq:eig_harmonic_spectral2_matrix_formulation}
can be decomposed into $L$ subproblems of the following form
\begin{align}\label{Eq:symmetric1_sub_problem2}
\cpfl^{(m)} \mathbf{f}^{(m)} = \lambda \mathbf{f}^{(m)},
\end{align}
for each $m \in\{0,,1,\,\hdots,\,L-1\}$.  The vector $\mathbf{f}^{(m)} $ of length $P(L-m)$ is constructed by
 $\mathbf{f}^{(m)}
=\left(\shc{f}{m m 0},\,\shc{f}{m m 1},\,\hdots,\,\shc{f}{m
m(P-1)},\hdots \shc{f}{(L-1)(m)(P-1)}\right)^{\rm T}$
and represents the spectral domain response of the
band-limited eigenfunction for given order $m$, given by
\begin{align}\label{Eq:spatial_eigen_subproblem2}
f^{(m)}(r,\theta) = \sum_{p=0}^{P-1} \sum_{\ell=m}^{L-1}
\shc{f}{\ell m p} Z_{\ell m p}(r,\theta,0).
\end{align}
Again, due to the fact that the matrix $\cpfl^{(m)}$ is real valued and symmetric, both $\mathbf{f}^{(m)}$ and $f^{(m)}(r,\theta)$ are real valued.

Using the spatial domain eigenfunction $f^{(m)}(r,\theta)$ for each
$m$ given by \eqref{Eq:spatial_eigen_subproblem1} and
\eqref{Eq:spatial_eigen_subproblem2}, the eigenfunction $f(\bv{r})$
can be obtained by scaling with the complex exponential $e^{i m \phi}$,
characterizing the variation of an eigenfunction along azimuth, giving
\begin{align}\label{Eq:rtheta_to_ball}
f(\bv{r}) =  f^{(m)}(r,\theta)e^{i m \phi} = \sum_{p=0}^{P-1} \sum_{\ell=m}^{L-1}
\shc{f}{\ell m p} Z_{\ell m p}(r,\theta,\phi), \quad  \bv{r} =
\bv{r}(r,\theta,\phi).
\end{align}

To summarise this subsection, the symmetry of the spatial region along
azimuth allows the decomposition of the large eigenvalue problem
that includes all angular orders $ -(L-1)\leq m \leq L-1$ into smaller $L$
subproblems, each for single angular order $m \in \{0, 1,2,\hdots,L-1 \}$.

\subsection{Circularly Symmetric and Radially Independent Region}

In the previous subsection we showed that the eigenvalue problem
decomposes into subproblems for circularly symmetric regions. In
addition to circular symmetry~(azimuthal symmetry as a special case), if the spatial region $R$ is also
radially independent, that is, if $r$ and $\theta$ are independent, then
the eigenvalue problem decomposes further into subproblems. For such
a symmetric region $R \dfn\{ {R}_1 \leq r \leq {R}_2,\, \theta_1
\leq \theta \leq \theta_2, \, 0 \leq \phi <2 \pi \}$, the integral
over the spatial region $R$ decouples as $\int_R\equiv
\int_{r=R_1}^{R_2} \int_{\theta=\theta_1}^{\theta_2}
\int_{\phi=0}^{2\pi} $, which can be incorporated to simplify
$\cpfb_{\ell m, \ell' m'}{(k,k')}$ given in
\eqref{Eq:symmetric1_matrix_simple1} as
\begin{align}\label{Eq:symmetric2_matrix_simple1}
\cpfb_{\ell m, \ell' m'}{(k,k')} =&\delta_{m m'}
\underbrace{\frac{2}{\pi} \int_{r={R}_1}^{{R}_2}r^2 \ddr{r}
 k k' j_\ell(kr) j_\ell'(k'r)}_{C_{\ell, \ell'}(k,k')} \nonumber
\\&
\times \underbrace{2\pi \int_{\theta=\theta_1}^{\theta_2}\sin\theta
\ddr{\theta}\, Y_{\ell m}(\theta,0)Y_{\ell'
m}^\ast(\theta,0)}_{G^m_{\ell,\ell'}}.
\end{align}
The integral $C_{\ell, \ell'}(k,k')$ can be evaluated
analytically for some special cases. When $k=k'$,
\begin{align}
C_{\ell, \ell'}(k,k) =
k^{2+\ell+\ell'}2^{(-2-\ell-\ell')}\Gamma(2+\ell+\ell')\left(
 F_1(\ell,\ell',k{R}_2) -
F_1(\ell,\ell',k{R}_1) \right),
\end{align}
where
\begin{align}
F_1(\ell,\ell',k{R}) &= {R}^{3+\ell+\ell'} {_p}F^q\bigg(\left[ \frac
{2 + \ell + \ell'}{2}, \frac{3 + \ell + \ell'}{2}, \frac {3 + \ell +
\ell'}{2}\right],\hdots \nonumber
\\& \quad  \hdots \left[ \frac{3 + 2\ell}{2},
 \frac{5 + \ell + \ell'}{2}, \frac{3 + 2\ell'}{2}, 2 + \ell + \ell'\right], -k^2 {R}^2 \bigg)
\end{align}
is the Hypergeometric generalized regularized
function.\footnote{Mathematica: HypergeometricPFQRegularized}
When $\ell=\ell'$,
\begin{align}
C_{\ell, \ell}(k,k') = \begin{cases} \frac{2\sqrt{k
k'}}{\pi(k^2-k'^2)}\bigg(
{R}_2^2 k'j_{\ell-1}(k'{R}_2) j_{\ell}(k{R}_2)
 - {R}_2^2 k j_{\ell-1}(k{R}_2) j_{\ell}(k'{R}_2)
\\
 - {R}_1^2 k'j_{\ell-1}(k'{R}_1)
j_{\ell}(k{R}_1)
 + {R}_1^2 k j_{\ell-1}(k{R}_1) j_{\ell}(k'{R}_1)
 \bigg) & k\neq k'\\
 \frac{k}{2\pi}\left( T(\ell,k,{R}_2)- T(\ell,k,{R}_2)
 \right) & k = k'
 \end{cases}
\end{align}
with $T(\ell,k,{R_1}) =  {R_1}^3\left( j_\ell^2(k{R_1}) - j_{\ell-1}(k{R_1})
j_{\ell+1}(k{R_1}) \right)$. The integral ${G^m_{\ell,\ell'}}$ can be evaluated
analytically for all $ \ell,\,\ell'\geq m$ by~\cite{Simons:2006}
\begin{align}
G^m_{\ell,\ell'} &= (-1)^m \frac{\sqrt{(2\ell+1)(2\ell'+1)}}{2}
\sum_{j=|\ell-\ell'|}^{|\ell+\ell'|}      \Bigg(
{\begin{array}{ccc}
  \ell & j & \ell' \\
  0 & 0 & 0
 \end{array} } \Bigg)
 \Bigg(
{\begin{array}{ccc}
  \ell & j & \ell' \\
  m & 0 & -m
 \end{array} } \Bigg)
\nonumber \\
& \quad \times  \left(P_{j-1}^0(\cos\theta_2) +
P_{j+1}^0(\cos\theta_1) - P_{j+1}^0(\cos\theta_2) -
P_{j-1}^0(\cos\theta_1) \right),
\end{align}
where the arrays of indices are Wigner-$3j$
symbols~\cite{Sakurai:1994}.
We note that the integral $C_{\ell,\ell'}(k,k')$ depends on only ${R}_1$ and ${R}_2$, and the integral denoted by $G^m_{\ell,\ell'}$ depends on only $\theta_1$ and $\theta_2$. For brevity, this dependence is not explicit in the notation. The azimuthal symmetry of the spatial region allows the eigenvalue problem to be decomposed into sub-problems and the independence between $r$ and $\theta$ enables the analytic computation of $\cpfb_{\ell m, \ell' m'}{(k,k')}$ when $\ell=\ell'$ or $k=k'$. For the Fourier-Bessel setting, the independence between $r$ and $\theta$ in the definition of the region does not allow further decomposition of the problem due to the coupling between the the radial and angular spectral components, characterized by, respectively, the Bessel functions $j_{\ell}(kr)$ and the spherical harmonics $Y_{\ell m}(\theta,\phi)$ (i.e., the harmonic index $\ell$ is shared).

For the Fourier-Laguerre setting, the independence between the radial component $r$ and angular colatitude $\theta$ allows us to further decompose the eigenvalue problem given in
\eqref{Eq:symmetric1_sub_problem2} since the radial and angular spectral domains are also decoupled. Considering the azimuthally and
radially symmetric region $R$, the elements of the matrix
$\mathbf{\cpfl}^m$ given in \eqref{Eq:symmetric1_matrix_new} can
be expressed as
\begin{align}
\cpfl^m_{\ell p,\ell'p' } =
\underbrace{\int_{r={R}_1}^{{R}_2}r^2 \ddr{r}
K_p(r)K_{p'}(r)}_{E_{p,p'}}
\underbrace{2\pi\int_{\theta=\theta_1}^{\theta_2}\sin\theta \ddr{\theta}\,
Y_{\ell m}(\theta,0)Y_{\ell' m}^\ast(\theta,0)}_{ G^m_{\ell,\ell'} }.
\end{align}
Since there is no dependence between the integrals along the radial
and angular spectral components, the fixed order eigenvalue problem
in \eqref{Eq:symmetric1_sub_problem2} can be decomposed into two
separate eigenvalue subproblems:
\begin{align}\label{Eq:symmetric2_FL_subproblem1}
\sum_{p'=0}^{P-1} E_{p',p} f_{p'} = \lambda^1 f_p, \quad f \in \lrplus ,
\end{align}
\begin{align}\label{Eq:symmetric2_FL_subproblem2}
\sum_{\ell'=m}^{L-1} G^m_{\ell',\ell}
 f_{\ell' m} = \lambda^2 f_{\ell m}, \quad f\in\lsph.
\end{align}
The eigenvalue problem in \eqref{Eq:symmetric2_FL_subproblem1}
maximizes the concentration of the band-limited signal defined on
$\rplus$ in the interval $r\in[R_1,R_2]$ and the eigenvalue problem
in \eqref{Eq:symmetric2_FL_subproblem2} maximizes the concentration
of the signal defined on $\untsph$ in the region characterized
by colatitude $\theta\in[\theta_1,\theta_2]$. We note that the eigenvalue problem in \eqref{Eq:symmetric2_FL_subproblem1} is
independent of angular order $m$ and is an algebraic eigenvalue
problem of size $P$, therefore its solution provides the spectral domain
representation of $P$ orthonormal eigenfunctions defined
on $\rplus$. The sum of eigenvalues is given by
\begin{equation}\label{Eq:shannon_radial_only}
N^P = \sum\limits_{p=0}^{P-1} E_{p,p},
\end{equation}
which also indicates the number
of significant eigenfunctions with eigenvalues near unity.
An analytic expression to evaluate $E_{p,p'}$ can be obtained by using the
definition of the Laguerre basis functions in \eqref{Eq:Laguerre_dfn} and \eqref{Eq:Laguerre_two_dfn}, yielding
\begin{align}
E_{p,p'} = \sqrt{\frac{p!\,p'!}{(p+2)!(p'+2)!}}& e^{-r}\, \sum_{j=0}^p \sum_{j'=0}^{p'} \frac{(-1)^{j+j'}}{j!\,j'!} \binom{p+2}{p-j} \binom{p'+2}{p'-j'}
\nonumber \\
& \times \int \ddr{r}\, e^{-r} r^{j+j'+2} ,
\end{align}
where the integral can be evaluated using the upper incomplete gamma function~\cite{Amore:2005} as
\begin{align}
\int \ddr{r}\, e^{-r} r^{j+j'+2} = (j+j'+2)! \, \sum_{a=0}^{j+j'+2}\frac{e^{-R_1}(R_1)^a - e^{-R_2}(R_2)^a  }{a!} .
\end{align}

Similarly, the solution of the eigenvalue problem in
\eqref{Eq:symmetric2_FL_subproblem2} for each $m$ gives rise to $L-m$ eigenfunctions defined on $\untsph$. The sum of eigenvalues for all orders is
given by
\begin{equation}
N_L = \sum\limits_{m=0}^{L-1}\sum\limits_{\ell=m}^{L-1}
G^m_{\ell,\ell} = \frac{L^2}{2}\left(\cos\theta_1 -
\cos\theta_2\right).
\end{equation}
Explicit expressions to determine the sum of eigenvalues for each order $m$ have
been derived in \cite{Simons:2006}, where the subproblem in
\eqref{Eq:symmetric2_FL_subproblem2} to find band-limited functions on the sphere with optimal spatial concentration in the polar cap region about the North pole has
been investigated in detail~\cite{Simons:2006}. By the decompositions due to the symmetry of the region, the
spherical Shannon number in \eqref{Eq:sum_eigen_values_FL} can be
expressed as
\begin{align}
N_{PL} = N^P N_L.
\end{align}

For the Fourier-Laguerre setting,
the band-limited eigenfunction $f^{(m)}(r,\theta)$ given by the solution of the fixed order eigenvalue problem in \eqref{Eq:symmetric1_sub_problem2} can be expressed in terms of the solution of the subproblems given in \eqref{Eq:symmetric2_FL_subproblem1} and \eqref{Eq:symmetric2_FL_subproblem2} as
\begin{align}
f^{(m)}(r,\theta) = \sum_{p=0}^{P-1} \sum_{\ell=m}^{L-1} f_p f_{\ell
m} Z_{\ell m p}(r,\theta,0),
\end{align}
where the eigenvalue $\lambda = \lambda^1 \lambda^2$ is a measure of
concentration in the spatial region $R$. The band-limited eigenfunction $f(\bv{r})$
can be obtained by scaling $f^{(m)}(r,\theta)$ with the complex exponential $e^{i m \phi}$,
as given in \eqref{Eq:rtheta_to_ball}.

\begin{remark}\label{remark:radially_independent}
If the region $R\subset\bll$ is not circularly symmetric but is radially independent, that is, $R = \rsphere \times \rradial$ with $\rsphere \subset \untsph$ and $\rradial \subset \rplus$, then due to the separability of the Fourier-Laguerre functions, the eigenvalue problem in the Fourier-Laguerre domain formulated in \eqref{Eq:eig_harmonic_spectral2_matrix_formulation} can be decomposed into subproblems to separately find band-limited functions in the region $\rradial$ as formulated in \eqref{Eq:symmetric2_FL_subproblem1} and in the region $\rsphere$ as formulated in \cite{Simons:2006}. Using \eqref{Eq:sum_eigen_values_FL} in conjunction with \eqref{Eq:shannon_radial_only}, the spherical Shannon number for a radially independent region $R = \rsphere \times \rradial$ is given by
\begin{align}\label{Eq:Shannon_symmetric_not_radial}
N_{PL} =   \frac{N^P \,L^2}{4\pi} \int_\rsphere {\emph{d}}^2\nu(\unit{r}).
\end{align}
\noindent It must be noted that no simplification is possible for the problem in the Fourier-Bessel domain, formulated in \eqref{Eq:eig_harmonic_spectral1}, due to the spectral coupling between the radial and angular components in the Fourier-Bessel functions.
\end{remark}

\begin{remark}
  Unlike the eigenvalue problem to find band-limited eigenfunctions in the Fourier-Bessel domain, the separation of the integrals over the spatial region due to the independence between $r$ and $\theta$ completely decouples the problem to find band-limited eigenfunctions in the Fourier-Laguerre domain, such that the eigenfunctions can be independently concentrated along radial component and angular component. Due to such decoupling, the eigenfunctions can be computed efficiently compared to the case when the spatial region does not exhibit symmetry and/or independence.  If one is interested in computing eigenfunctions in an arbitrary spatial region, the region can be approximated by the union of $T$ subregions $R = R_1 \cup R_2\cup\hdots\cup R_T$ such that each subregion is azimuthally symmetric and radially independent. Furthermore, as indicated earlier, an azimuthally symmetric region can be rotated to a circularly symmetric region.  The rotation of eigenfunctions can be performed in the spectral domain through the use of Wigner-$D$ functions~\cite{Kennedy-book:2013,Sakurai:1994}.
\end{remark}

%
\begin{figure}[t]
    \centering
    \subfloat[Fourier-Bessel]{
        \includegraphics[scale=0.4]{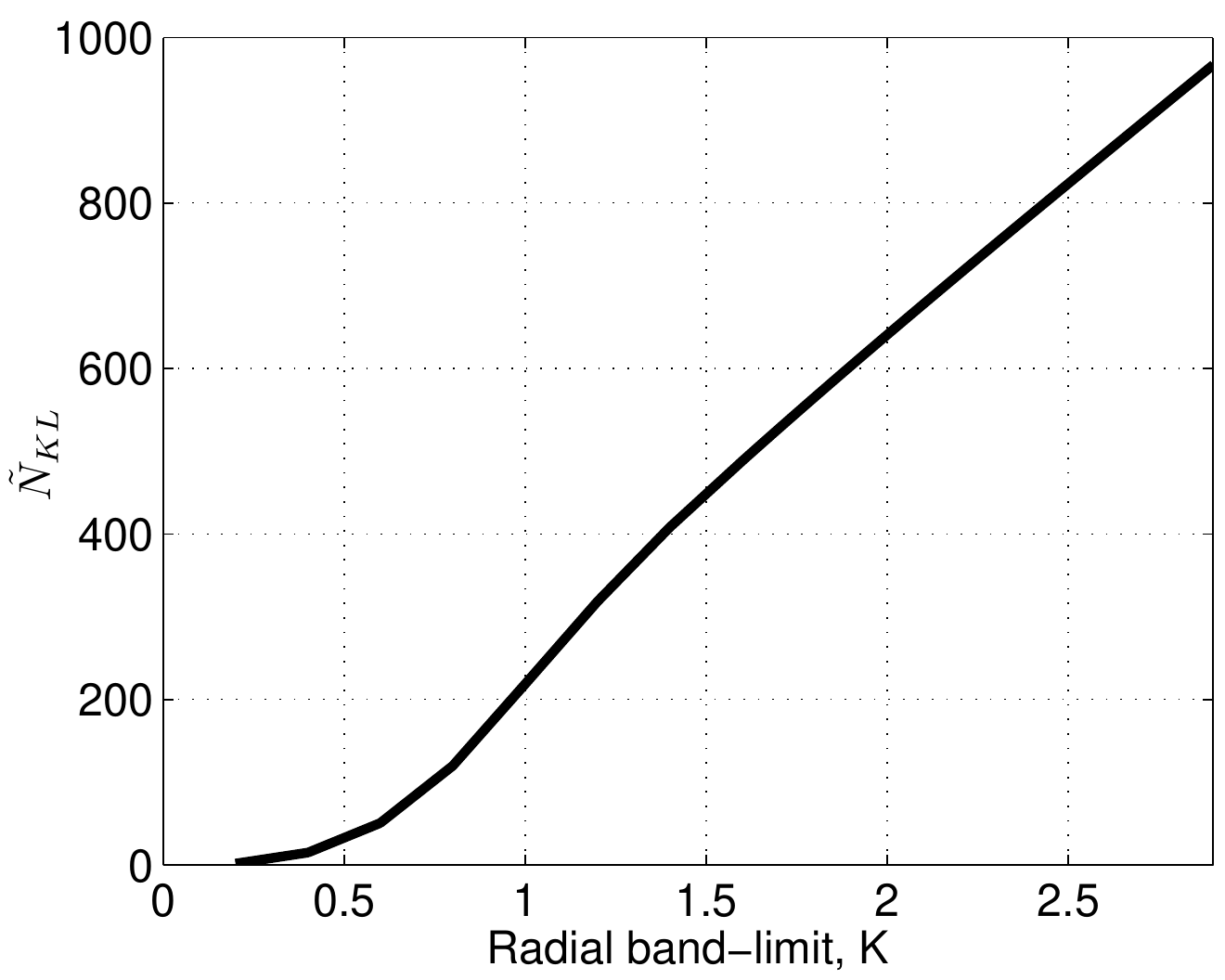}}
        \hspace{6mm}
        \subfloat[Fourier-Laguerre]{
        \includegraphics[scale=0.4]{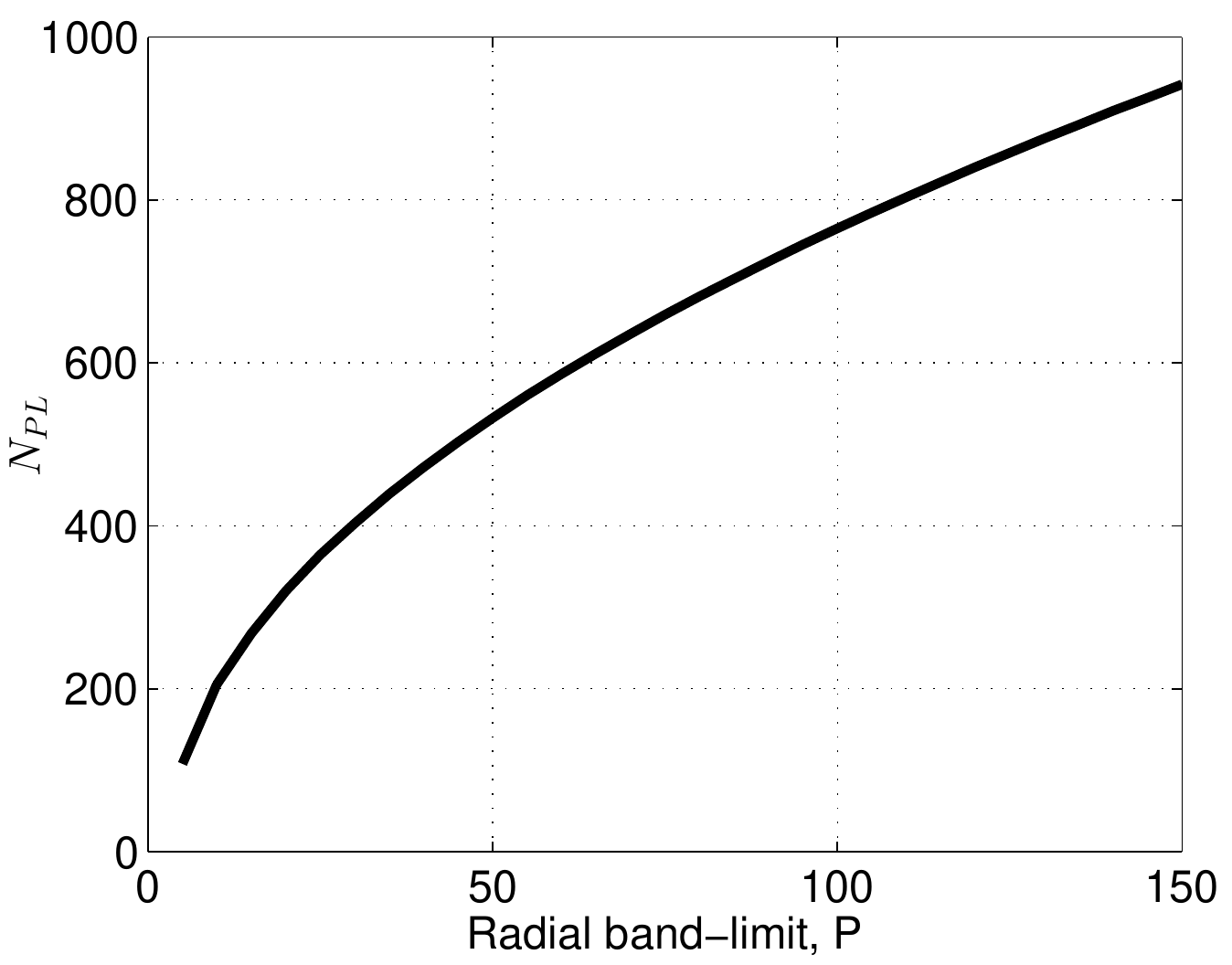}}
    \hspace{-3mm}

    \caption{\textbf{Analytical plots for the spherical Shannon
        number} for (a) $\tilde N_{KL}$ given in
      \eqref{Eq:sum_eigen_values_FB} for different values of
      band-limit $K$ and (b) $N_{PL}$ given in
      \eqref{Eq:sum_eigen_values_FL} for different values of
      band-limit $P$. The angular band-limit is $L=20$ and the spatial
      region of concentration is $R=\{ 15 \leq r \leq 25,\, \pi/8 \leq
      \theta \leq 3\pi/8, \, 0 \leq \phi <2 \pi \}$. }
    \label{fig:shannon_comparison}
\end{figure}

\begin{figure}[t]
    \centering
    \subfloat[Fourier-Bessel]{
        \includegraphics[scale=0.4]{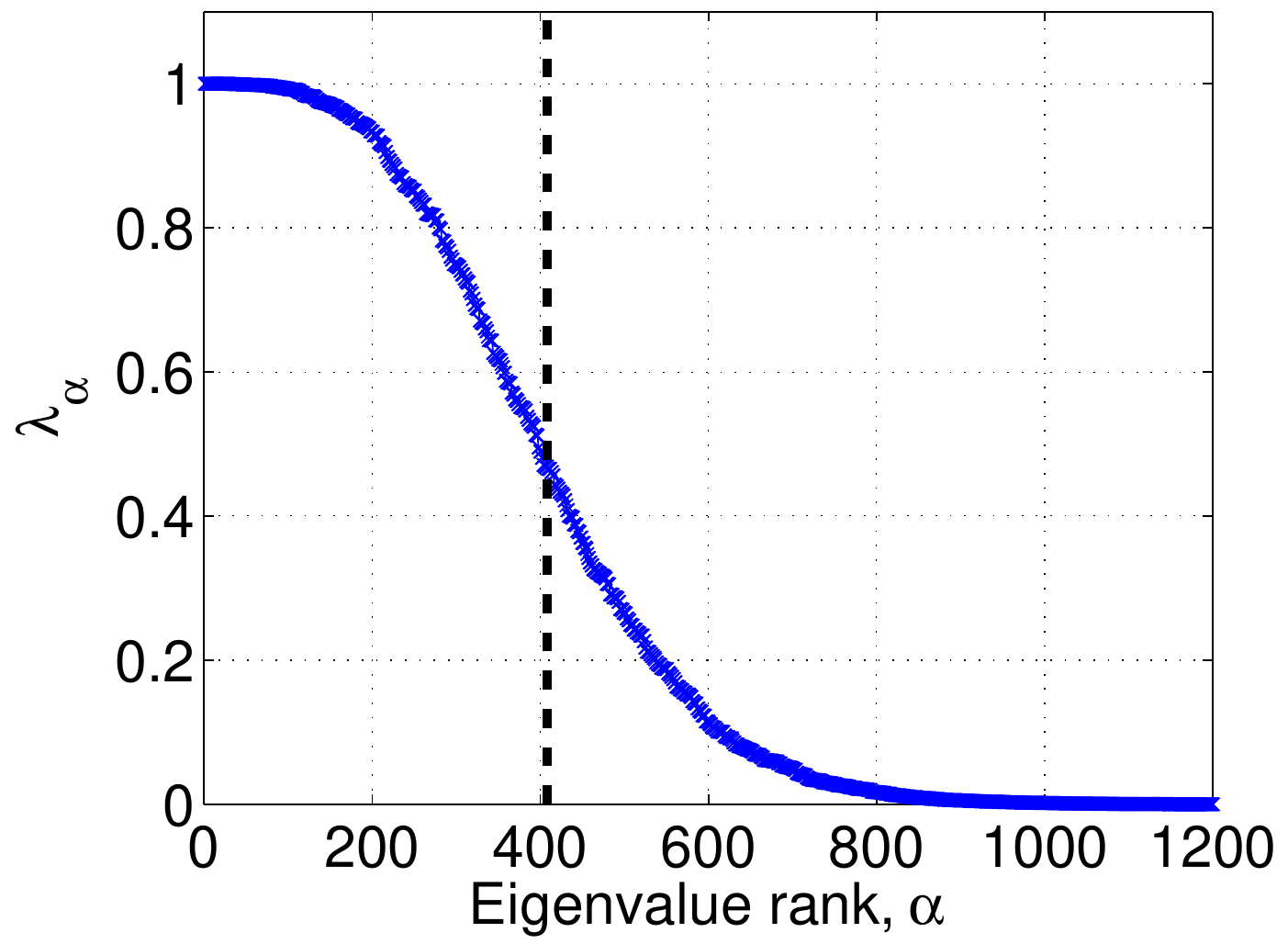}}
        \subfloat[Fourier-Laguerre]{
        \includegraphics[scale=0.4]{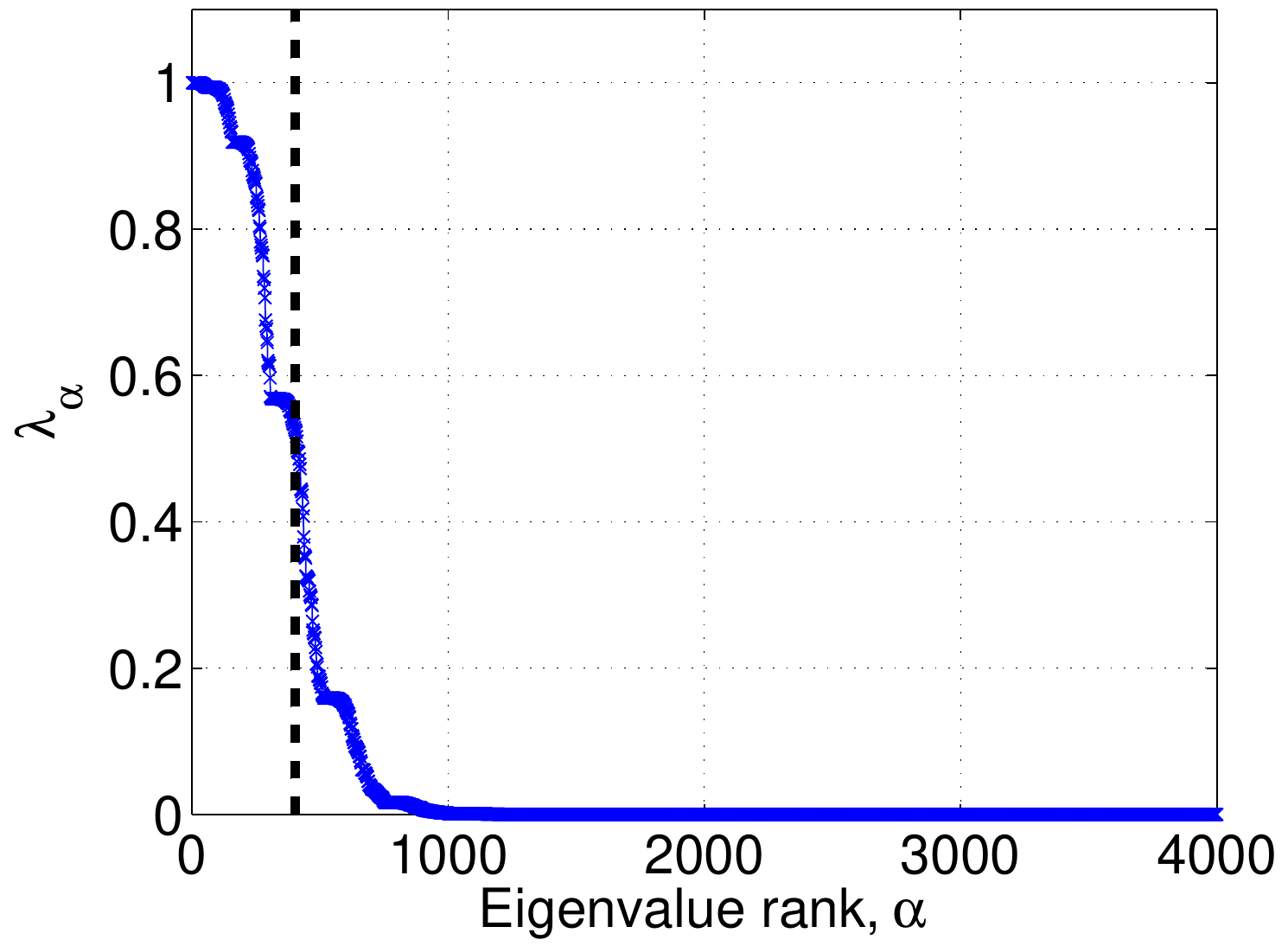}}
    \hspace{-3mm}

    \subfloat[Fourier-Laguerre radial]{
        \includegraphics[scale=0.4]{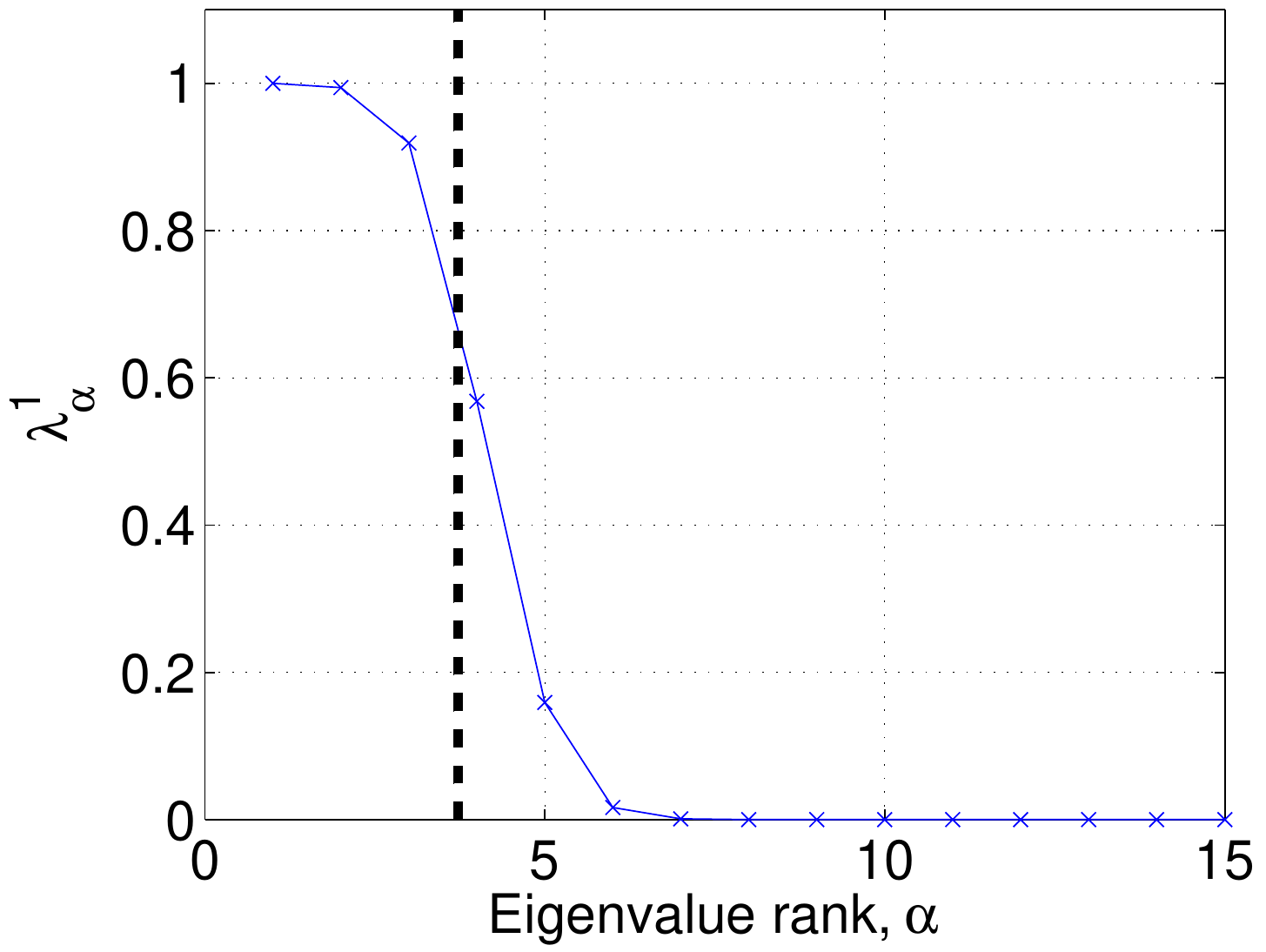}}
        \subfloat[Fourier-Laguerre angular]{
        \includegraphics[scale=0.4]{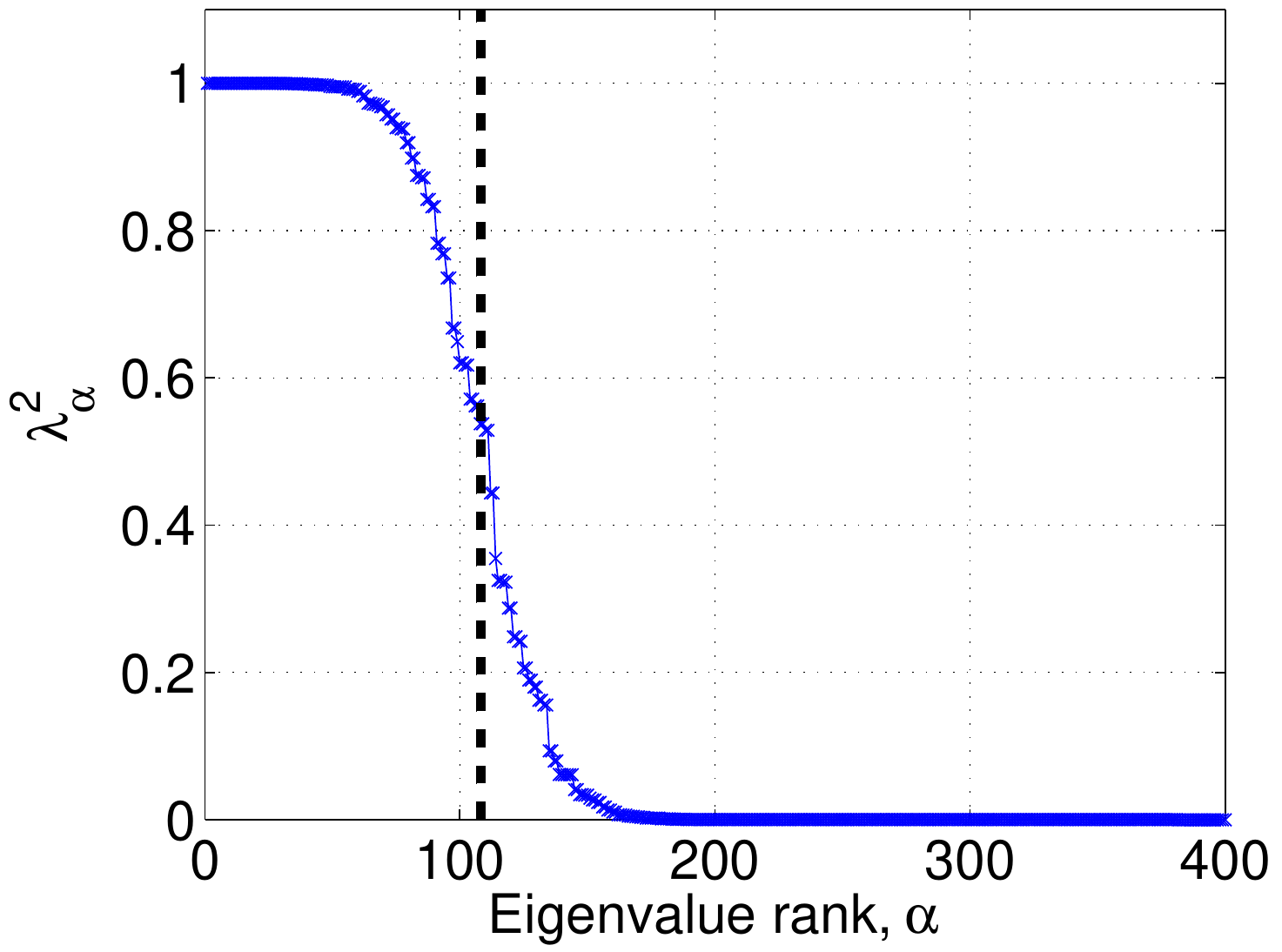}}

      \caption{\textbf{Spectrum of eigenvalues} for the eigenvalue
        problems formulated in: (a) \eqref{Eq:eig_harmonic_spectral1}
        to find spatially concentrated eigenfunctions band-limited in
        the Fourier-Bessel domain in the spectral region $\tilde
        A_{(1.4,20)}$; and (b)
        \eqref{Eq:eig_harmonic_spectral2_matrix_formulation} to find
        spatially concentrated eigenfunctions band-limited in the
        Fourier-Laguerre domain in the spectral region $A_{(30,20)}$.
        The spatial region of concentration is \mbox{$R=\{ 15 \leq r
          \leq 25,$} \mbox{$\pi/8 \leq \theta \leq 3\pi/8, \, 0 \leq
          \phi <2 \pi \}$}. Since the region is circularly symmetric
        and radially independent, the problem in the Fourier-Laguerre
        domain can be decomposed into two separate eigenvalue problems
        in radial and angular regions formulated in
        \eqref{Eq:symmetric2_FL_subproblem1} and
        \eqref{Eq:symmetric2_FL_subproblem2}, for which the spectrum
        of eigenvalues is plotted in (c) and (d) respectively. Note
        that the transition of the eigenvalue spectrum from unity to
        zero in (b) is not smooth since the radial and angular
        components are interspersed.  The sum of eigenvalues, given by
        the spherical Shannon number for each case is (a) $\tilde
        N_{KL} = 408.33$ , (b) $N_{PL} = 403.21$, (c) $N^P = 3.72$ and
        (d) $N_L= 108.24$, as indicated by the vertical dashed
        lines. The spherical Shannon numbers estimate the location of
        the transition in the eigenvalue spectrum accurately.}
    \label{fig:eigen_spectrum_lag_sphere}
\end{figure}

\begin{figure*}[!htb]
    \centering
    \vspace{-10mm}
    \hspace{-14mm}
    \subfloat{
        \includegraphics[scale=0.180]{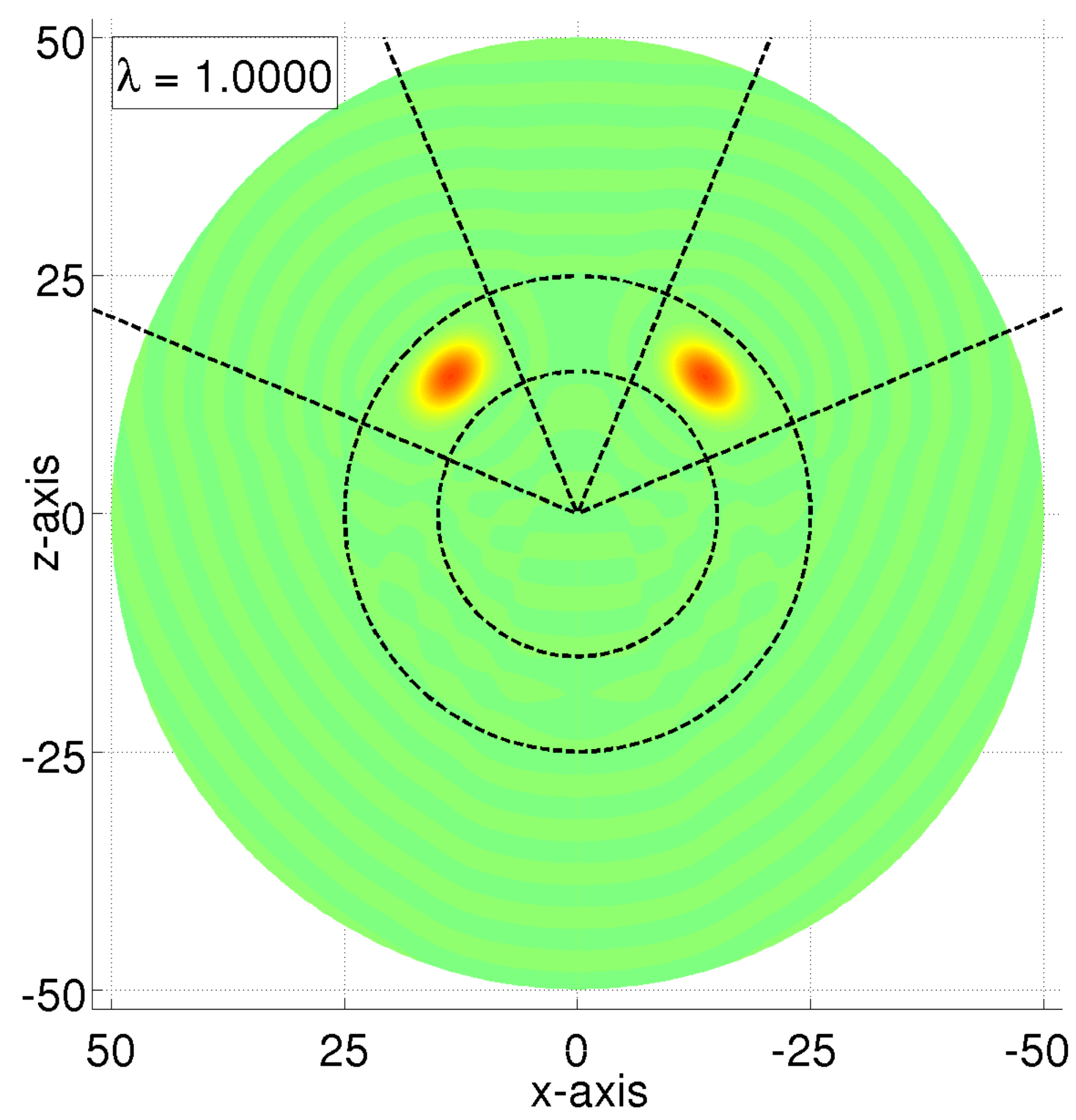}}
    \hspace{3mm}
    \subfloat{
        \includegraphics[scale=0.180]{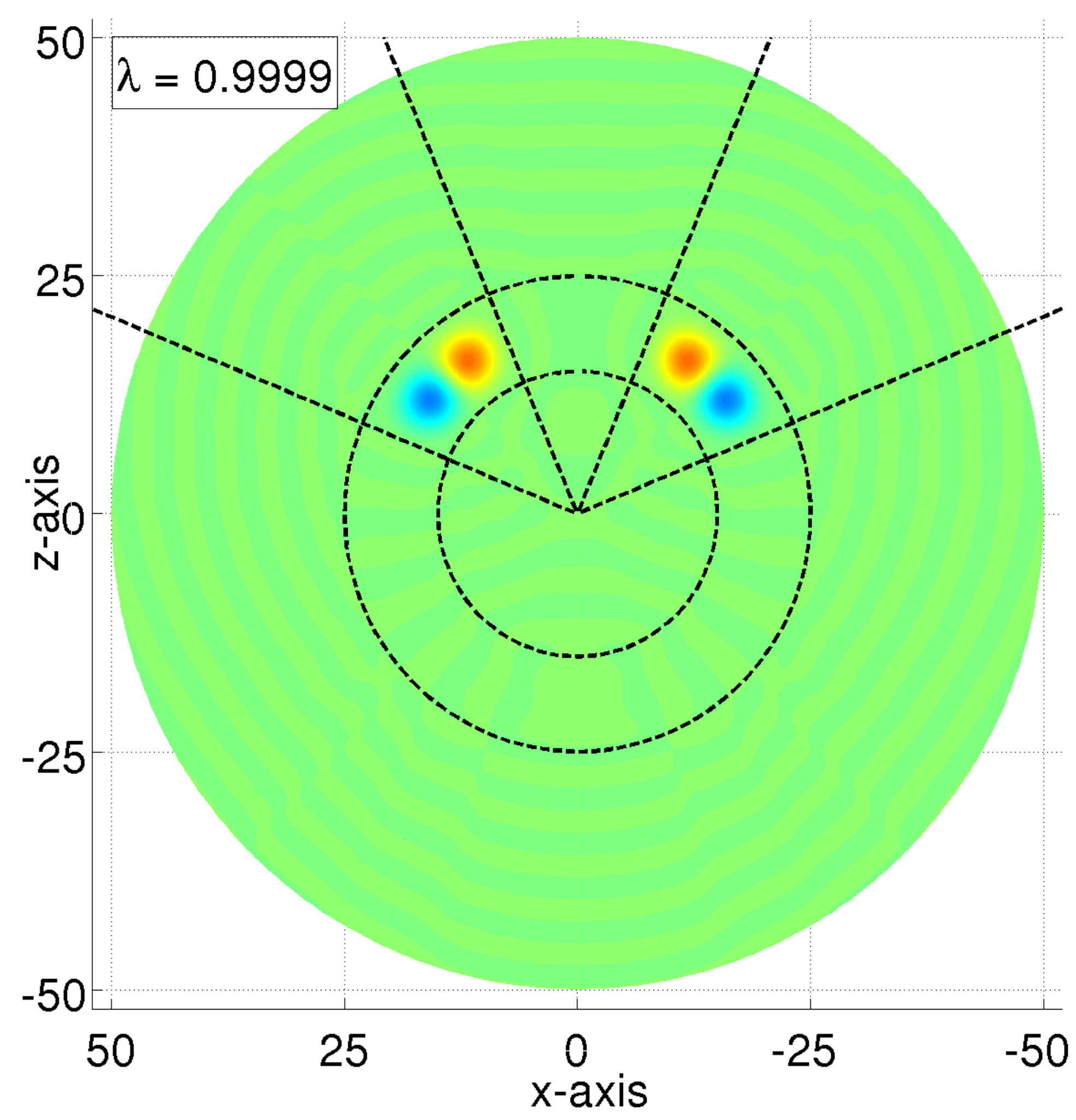}}
    \hspace{3mm}
    \subfloat{
        \includegraphics[scale=0.180]{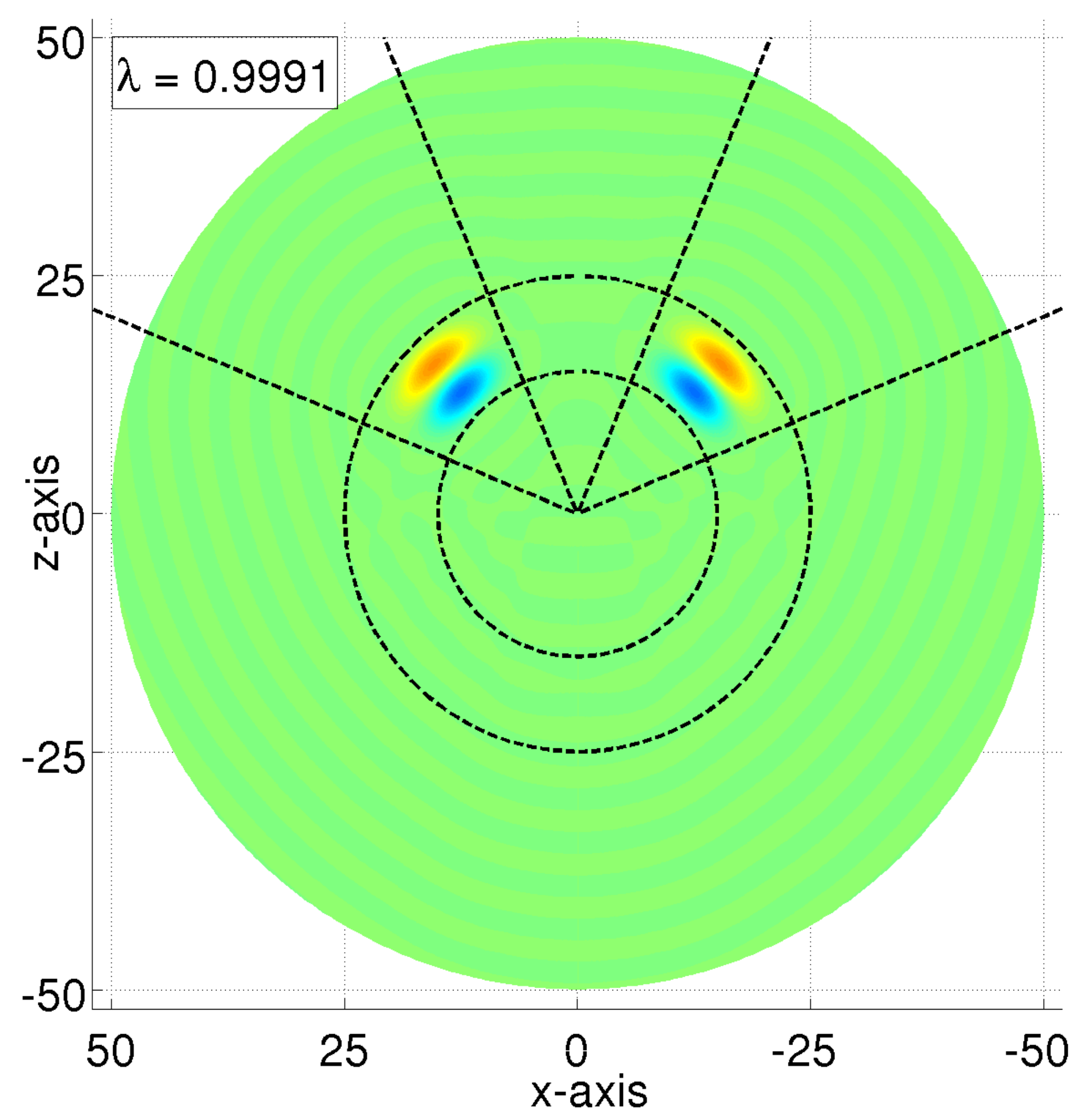}}

    \vspace{-3mm}
    \hspace{-14mm}
    \subfloat{
        \includegraphics[scale=0.180]{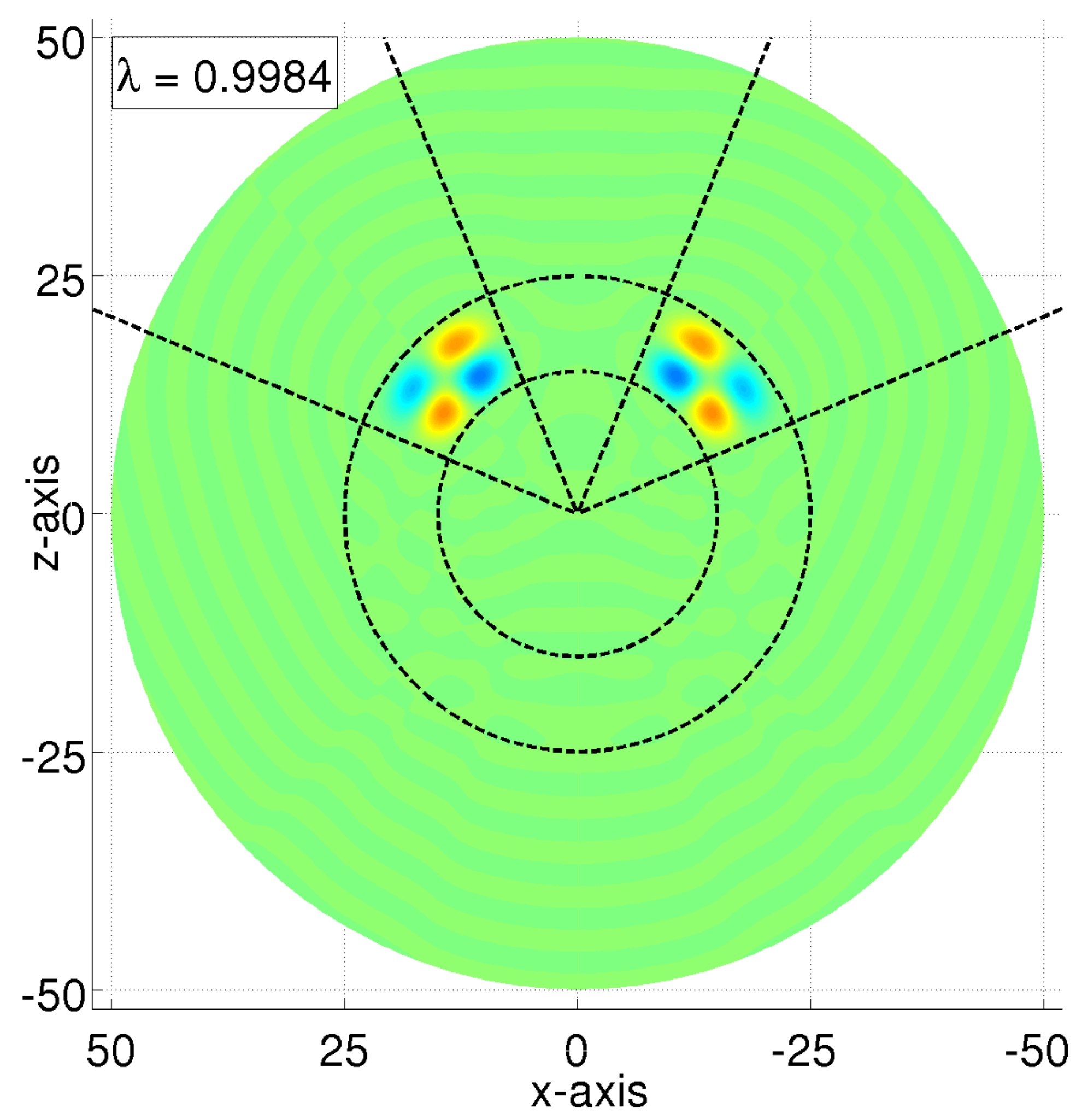}}
    \hspace{3mm}
    \subfloat{
        \includegraphics[scale=0.180]{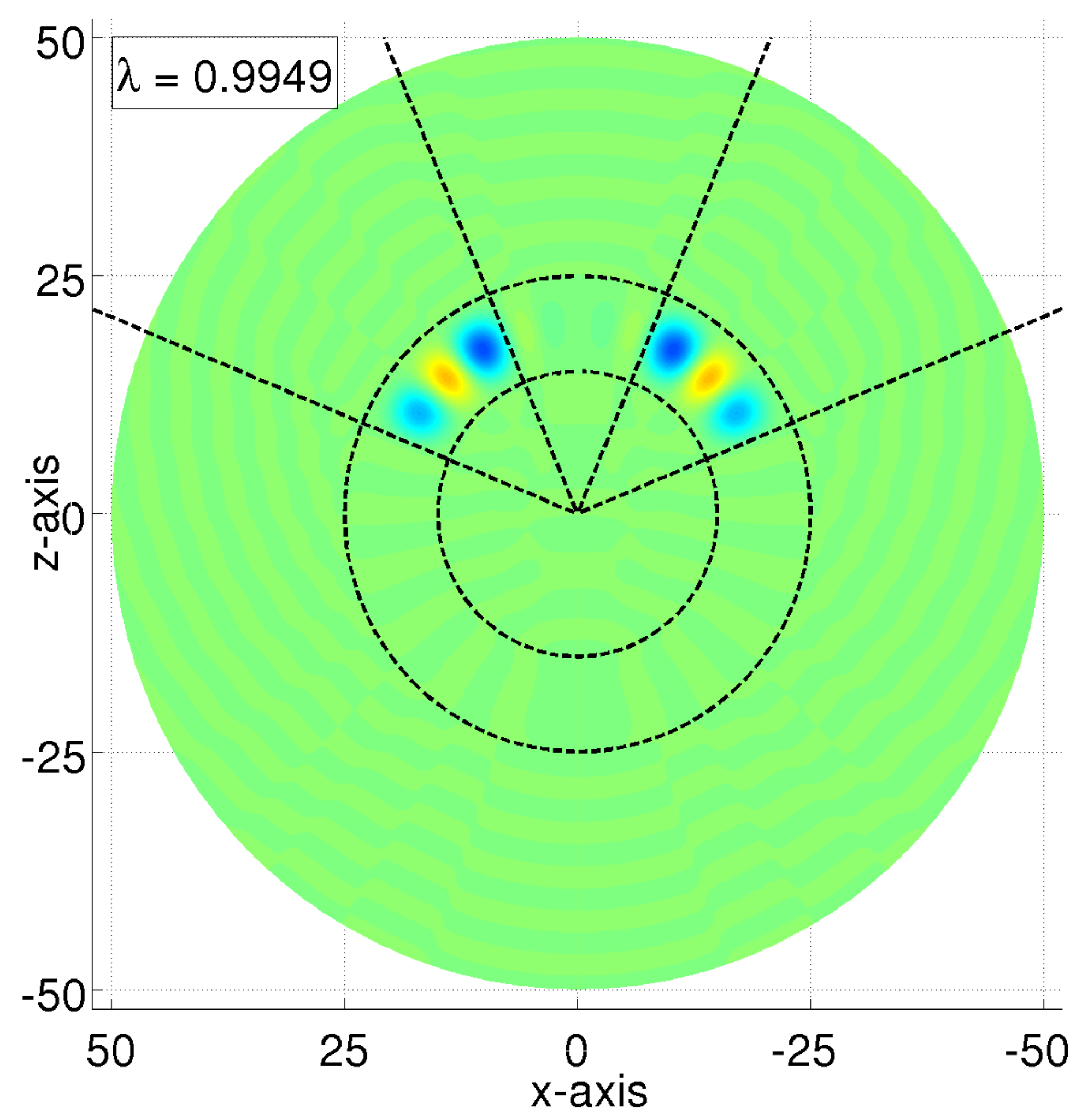}}
    \hspace{3mm}
    \subfloat{
        \includegraphics[scale=0.180]{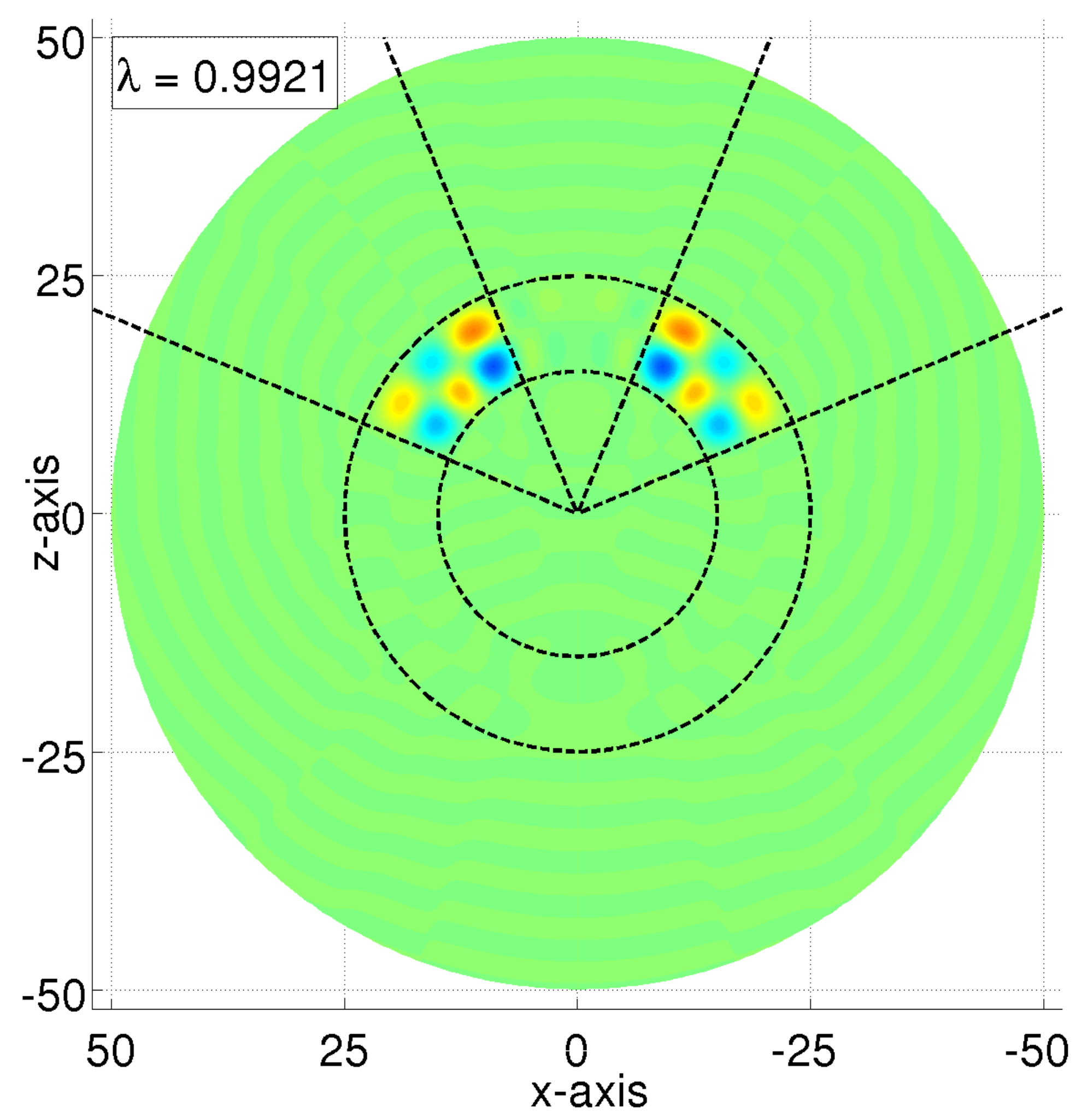}}

    \vspace{-3mm}
    \hspace{-14mm}
    \subfloat{
        \includegraphics[scale=0.180]{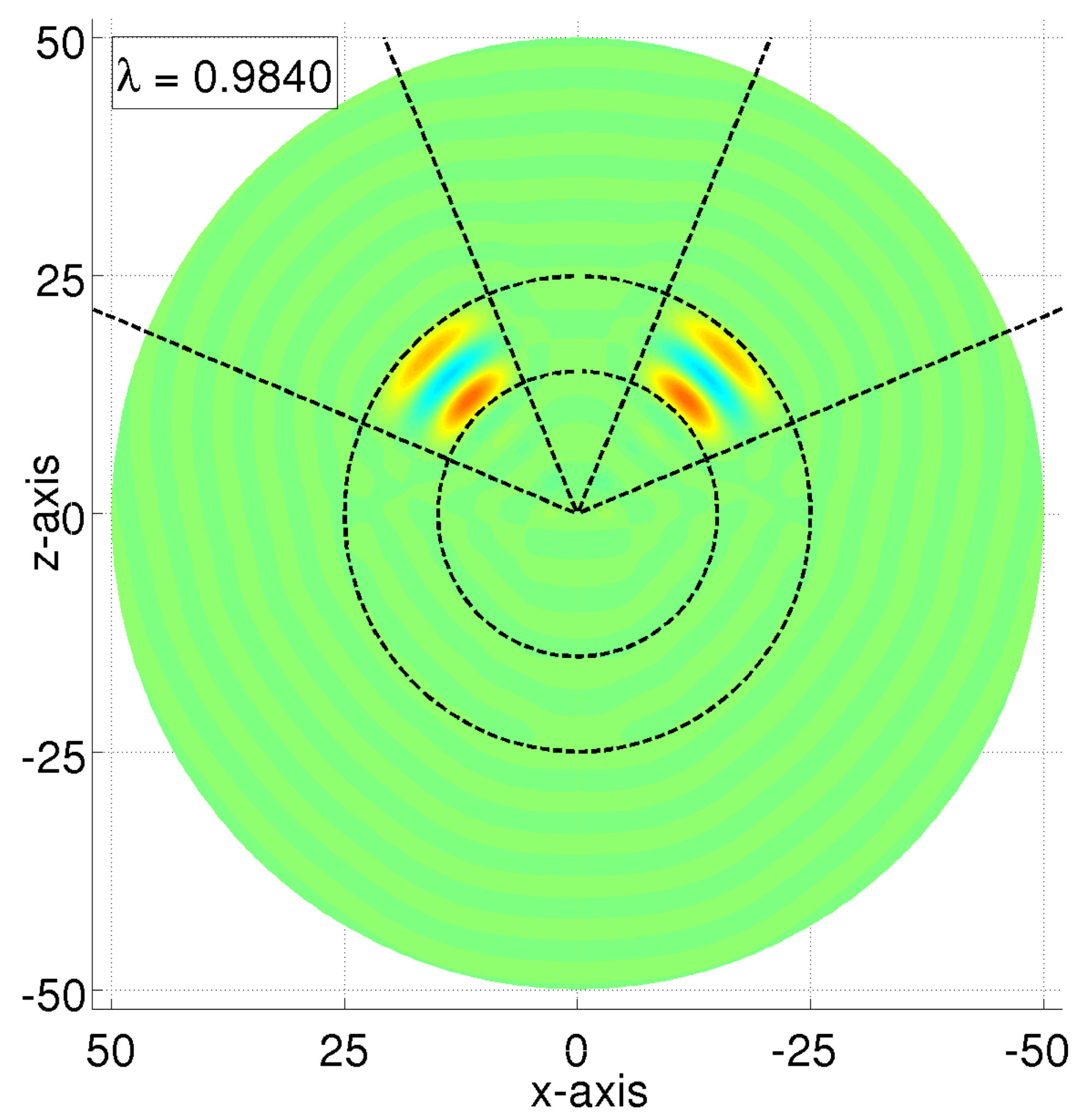}}
    \hspace{3mm}
    \subfloat{
        \includegraphics[scale=0.180]{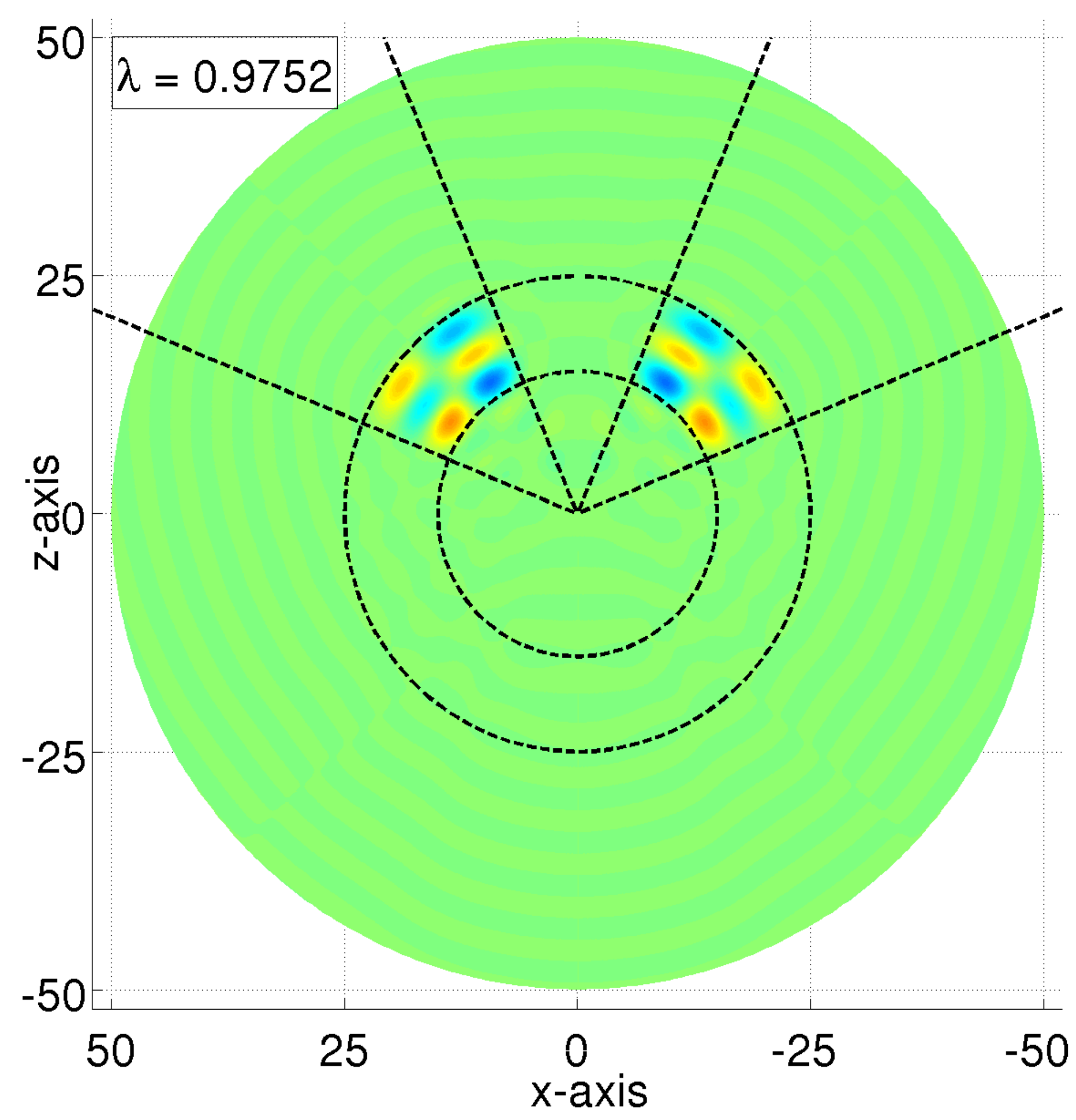}}
    \hspace{3mm}
    \subfloat{
        \includegraphics[scale=0.180]{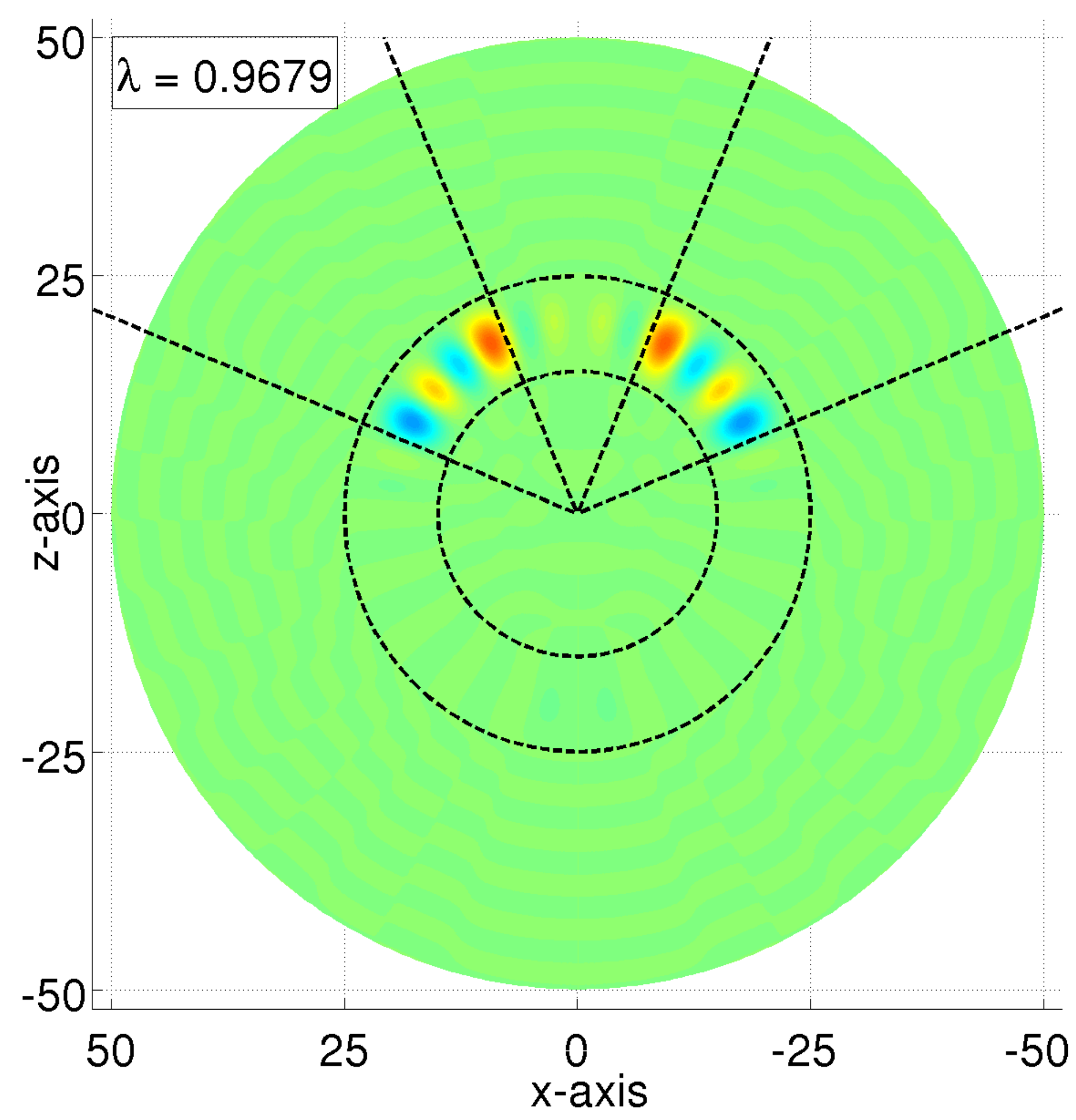}}

    \vspace{-3mm}
    \hspace{-14mm}
    \subfloat{
        \includegraphics[scale=0.180]{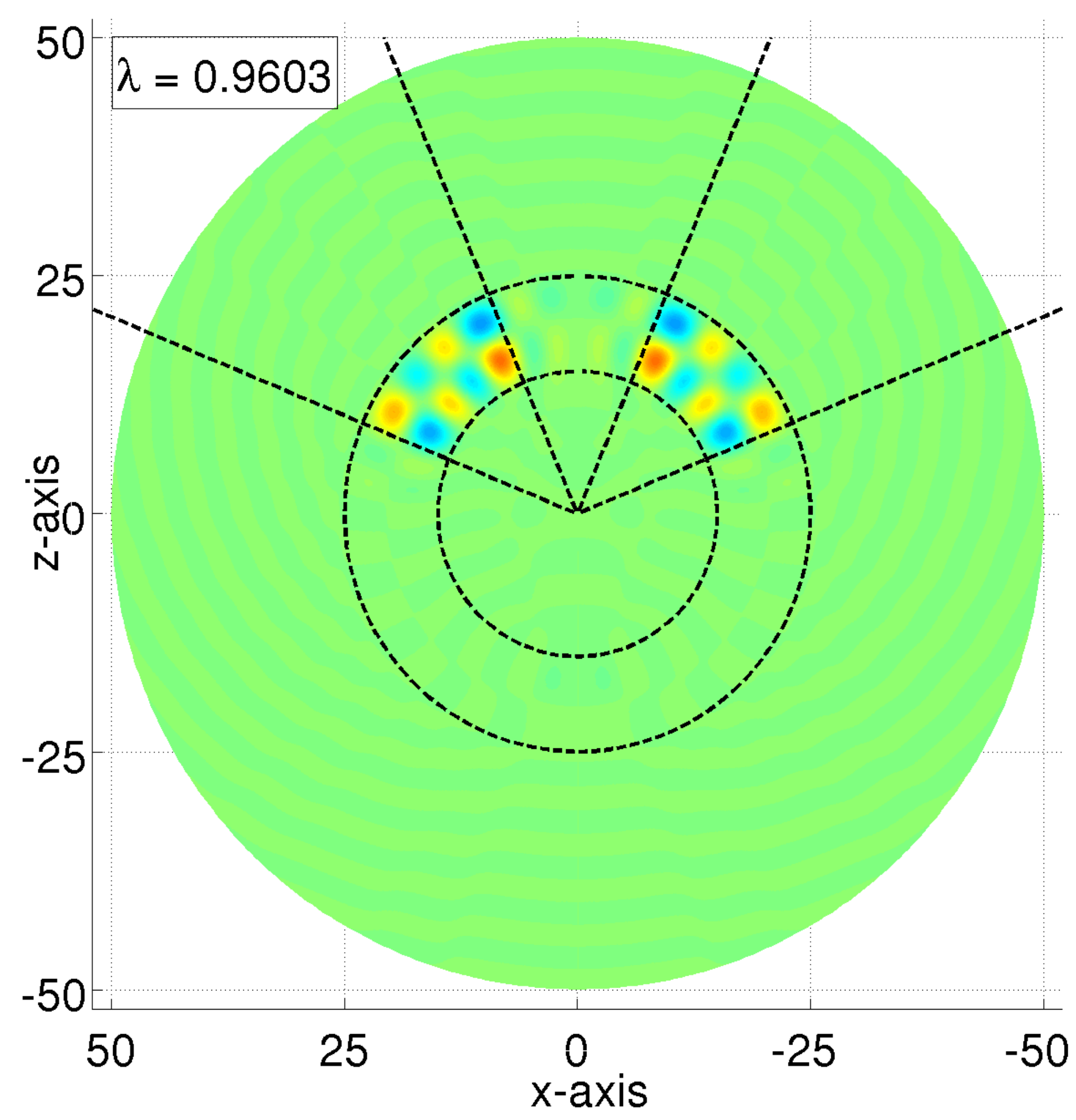}}
    \hspace{3mm}
    \subfloat{
        \includegraphics[scale=0.180]{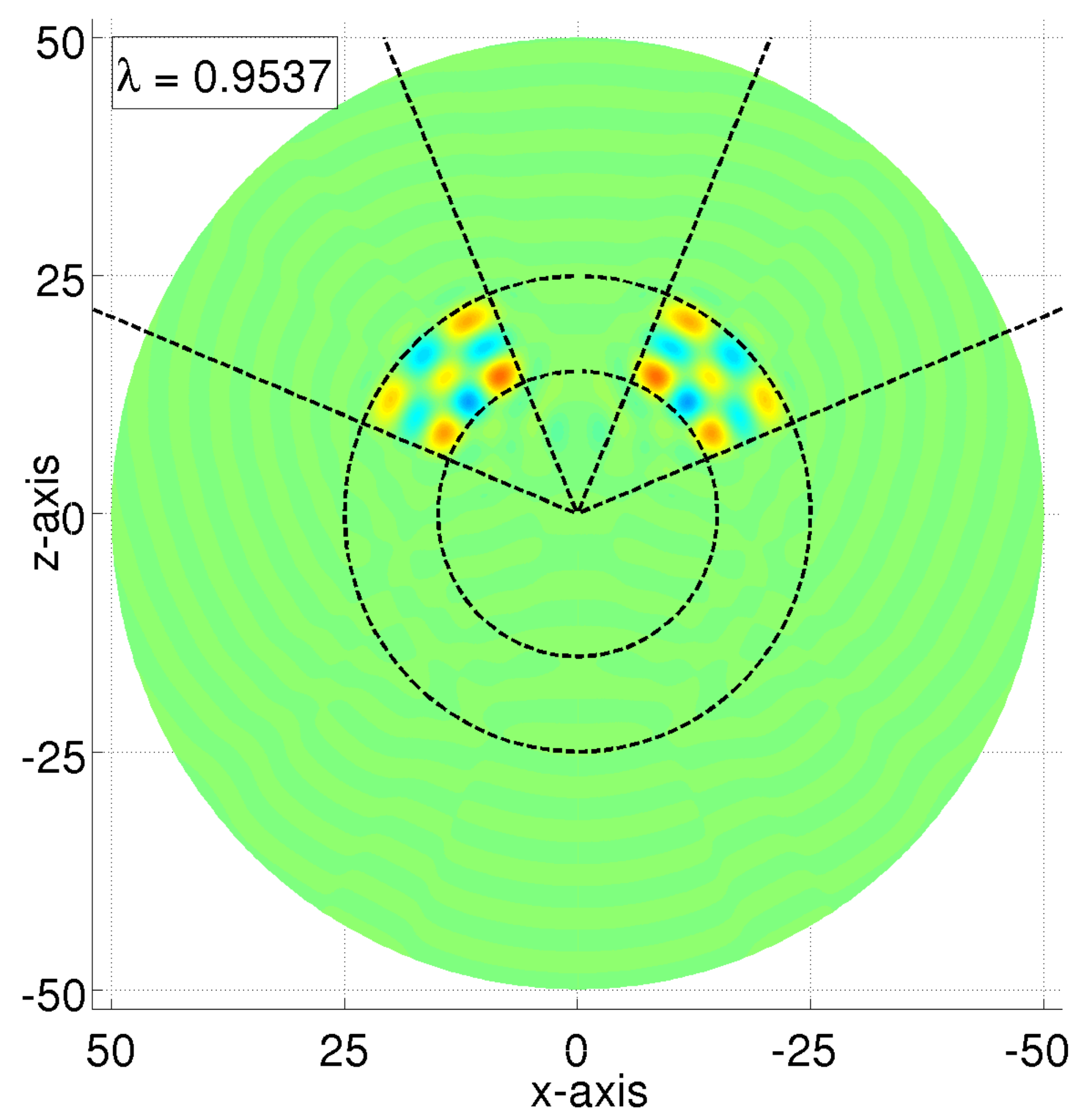}}
    \hspace{3mm}
    \subfloat{
        \includegraphics[scale=0.180]{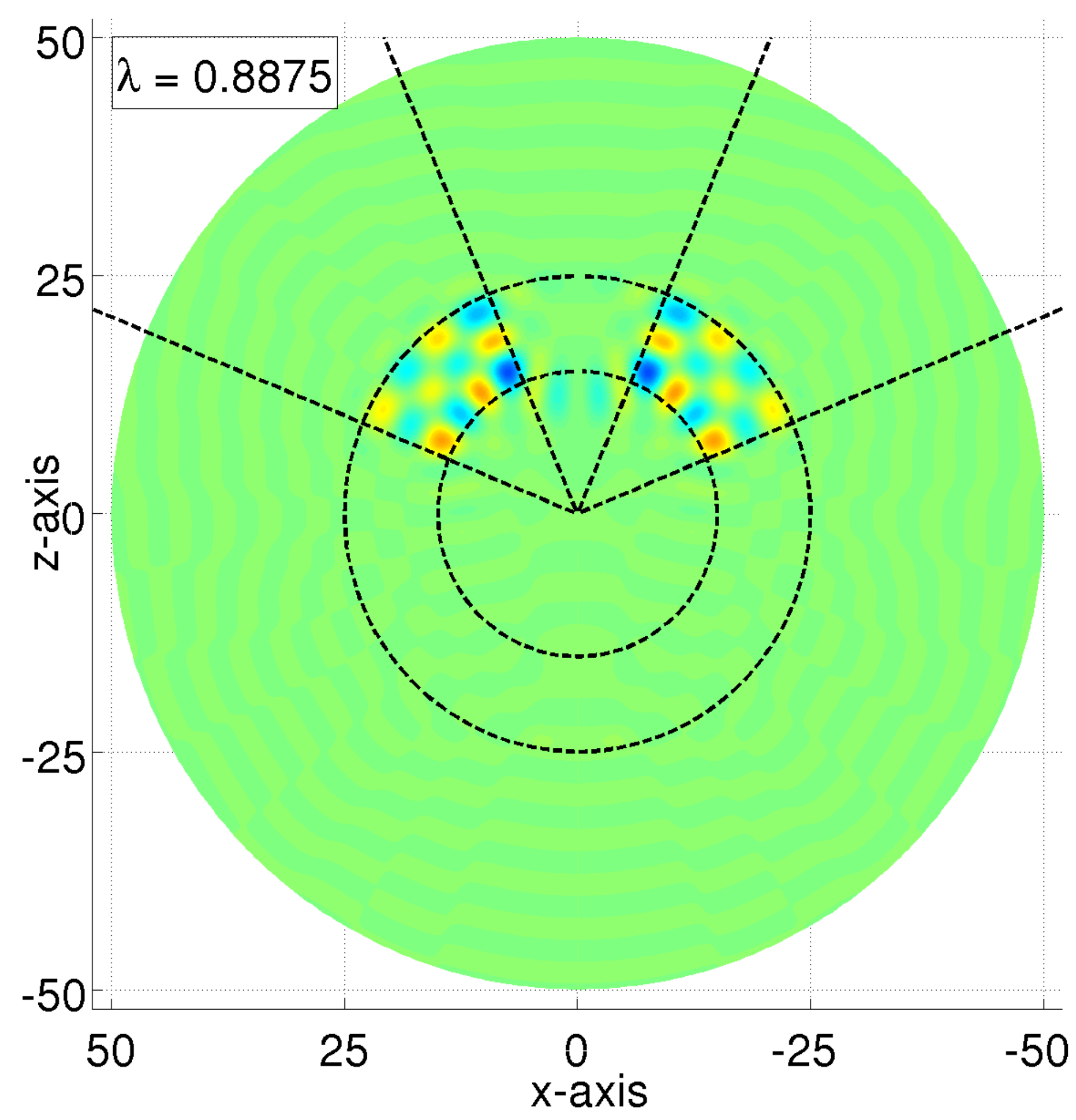}}


    \vspace{-3mm}

    \subfloat{
        \includegraphics[scale=0.25]{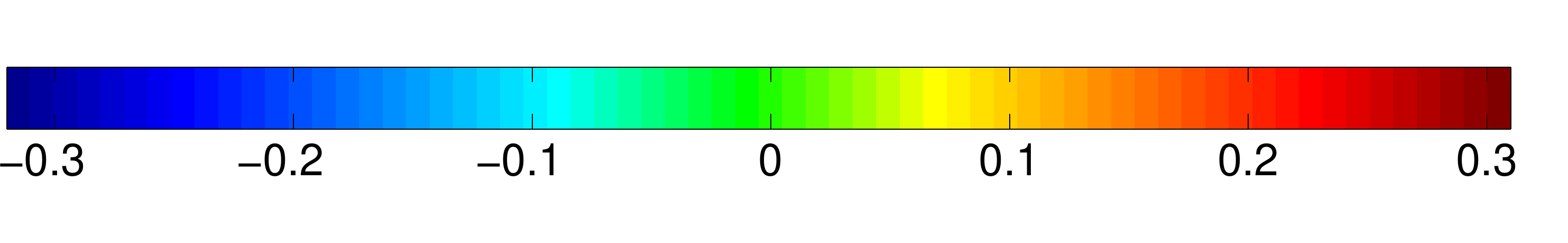}}
    \vspace{-3mm}

    \caption{\textbf{Fourier-Bessel band-limited spatially
        concentrated eigenfunctions} $f^m(r,\theta)$ given in
      \eqref{Eq:spatial_eigen_subproblem1}, obtained as solutions of
      the fixed-order eigenvalue problem in
      \eqref{Eq:symmetric1_sub_problem1} for $m=2$.  Each
      eigenfunction $f^m(r,\theta)$ is independent of $\phi$ and is
      plotted for $r\leq50$.  The spatial region of concentration is
      $R=\{ 15 \leq r \leq 25,\, \pi/8 \leq \theta \leq 3\pi/8,$ \mbox{$0
      \leq \phi <2 \pi \}$} and is azimuthally symmetric and radially
      independent.  The dependence of the eigenfunction in $\phi$ is
      $e^{im\phi}$, as given in \eqref{Eq:rtheta_to_ball}. The
      eigenfunctions are band-limited in the Fourier-Bessel domain
      within the spectral region $\tilde A_{(1.4,20)}$. The eigenvalue
      $\lambda$ associated with each eigenfunction is a measure of
      spatial concentration within the spatial region $R$.} \label{fig:eigen_bessel_spatial}
\end{figure*}

\begin{figure*}[!htb]
    \vspace{-5mm}
    \centering
    \hspace{-15mm}
    \subfloat{
        \includegraphics[scale=0.22]{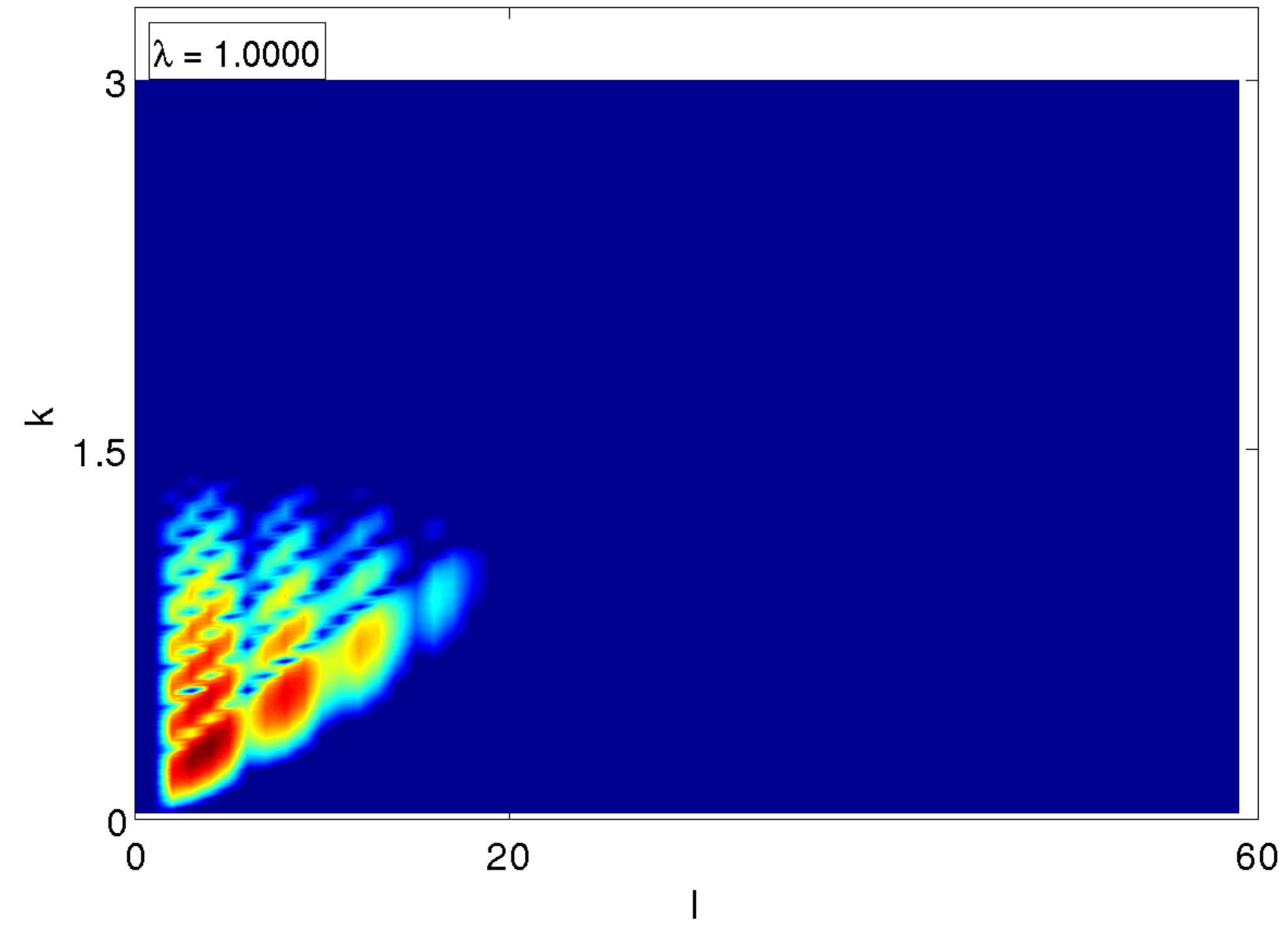}}
    \hspace{-2mm}
    \subfloat{
        \includegraphics[scale=0.22]{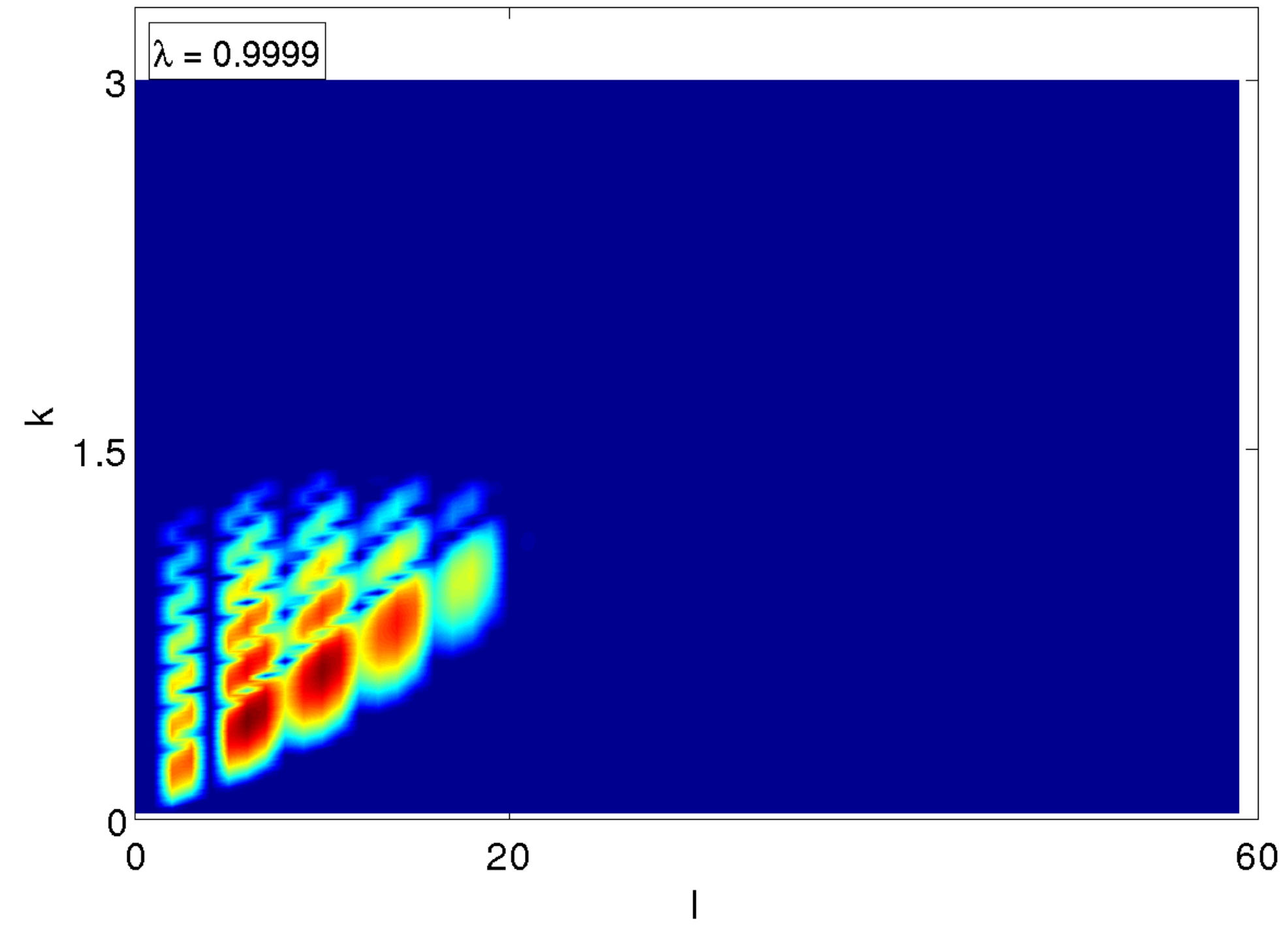}}
    \hspace{-2mm}
    \subfloat{
        \includegraphics[scale=0.22]{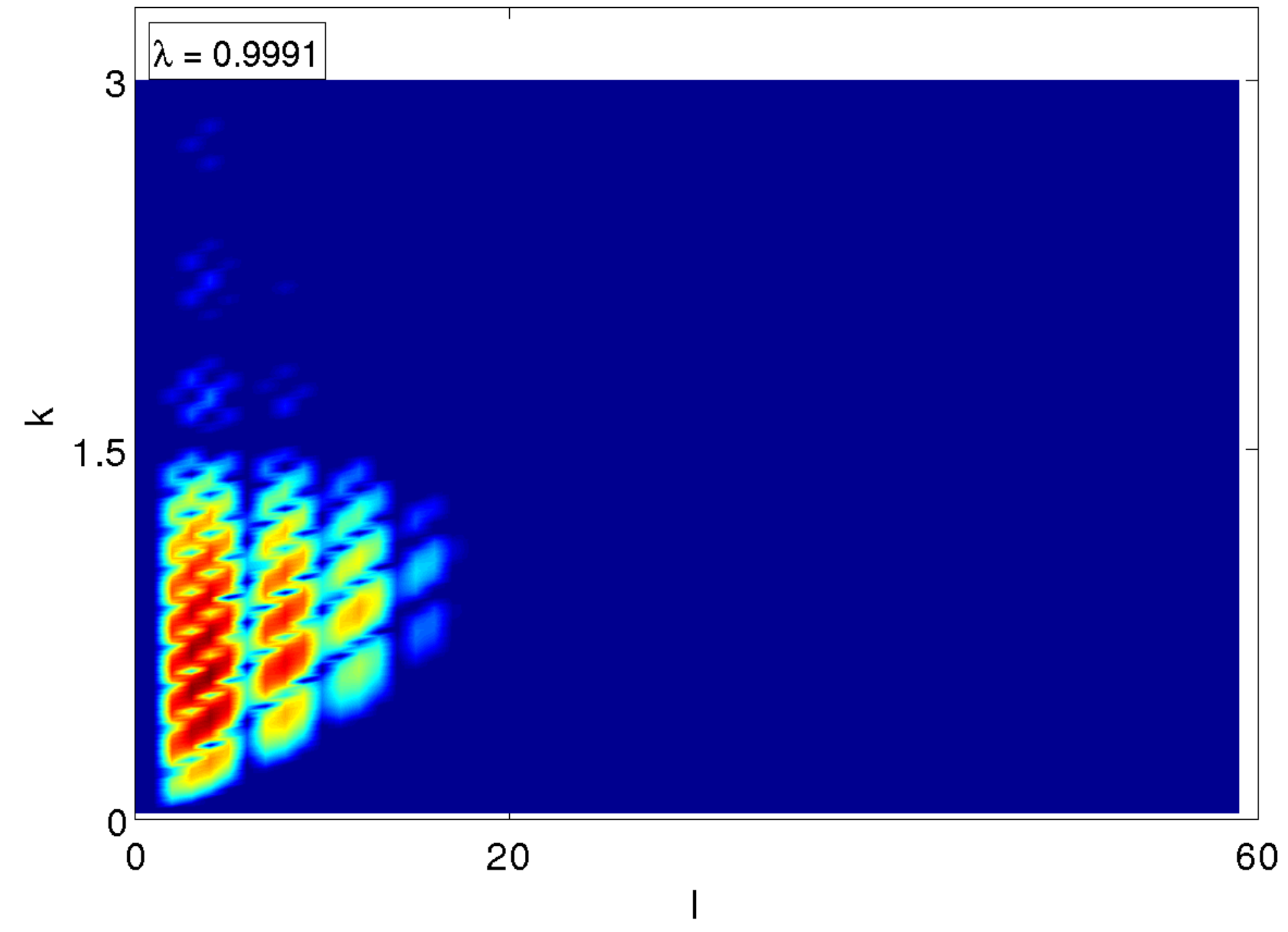}}

    \vspace{-3mm}
    \hspace{-15mm}
    \subfloat{
        \includegraphics[scale=0.22]{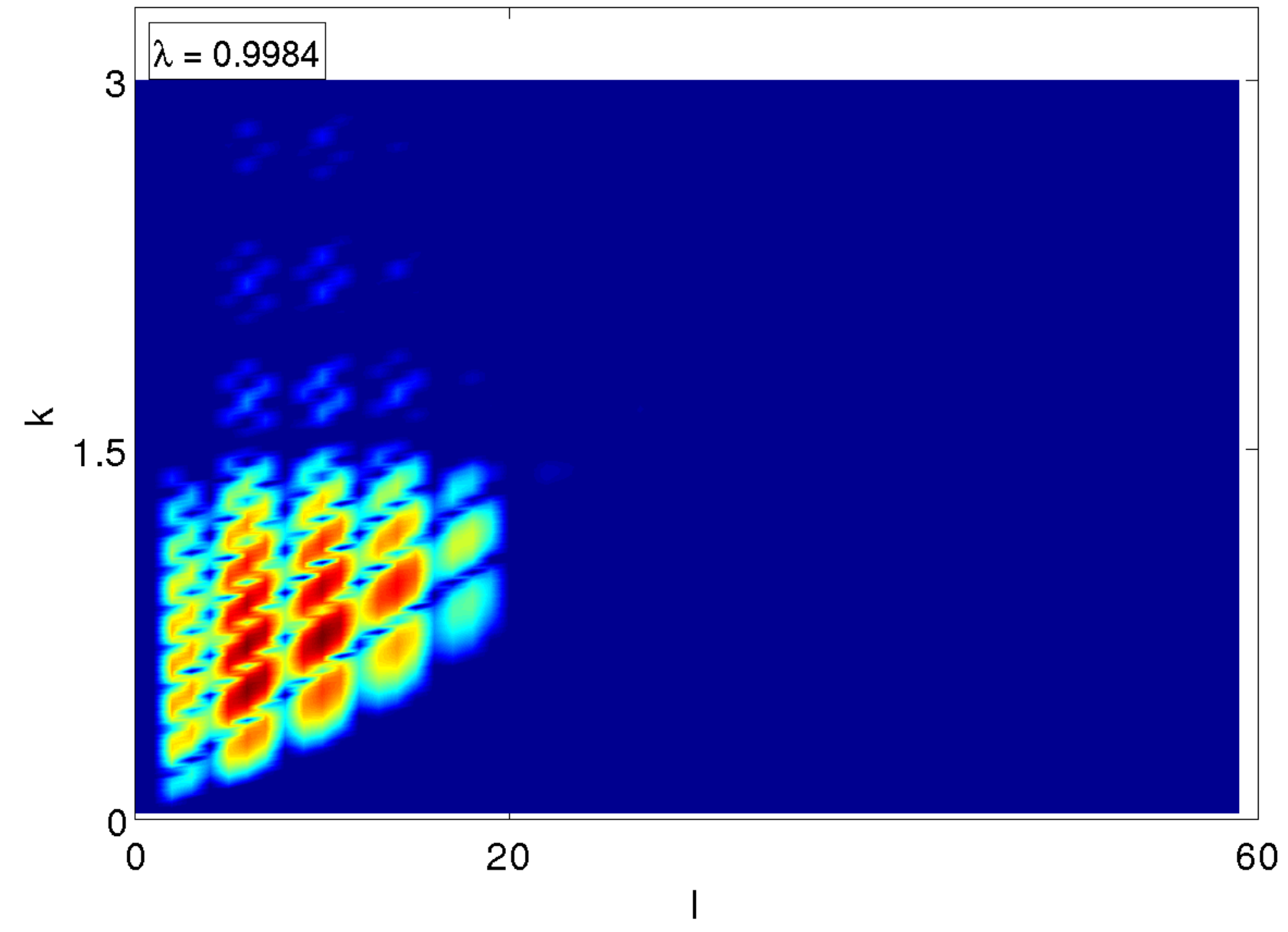}}
    \hspace{-2mm}
    \subfloat{
        \includegraphics[scale=0.22]{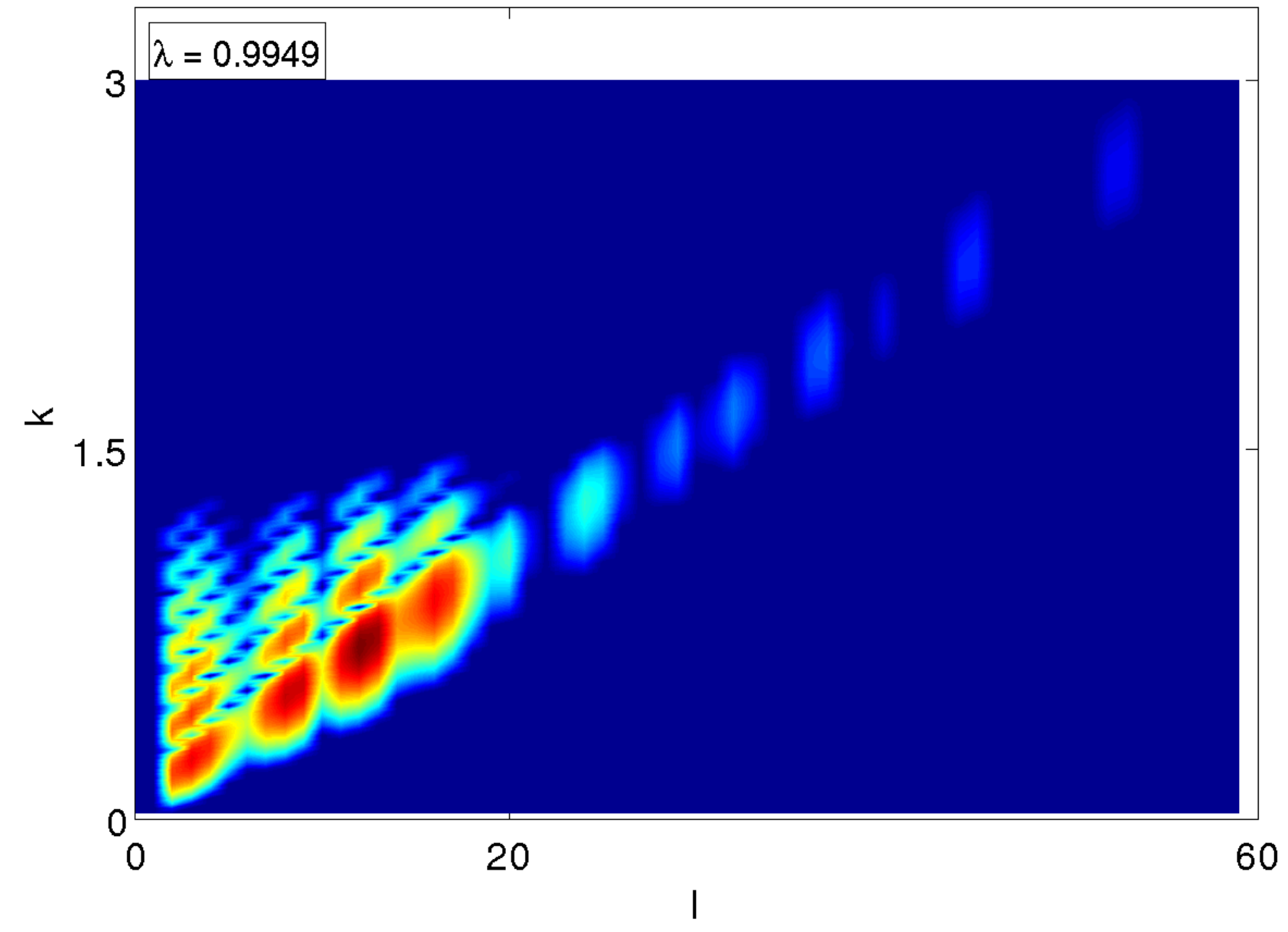}}
    \hspace{-2mm}
    \subfloat{
        \includegraphics[scale=0.22]{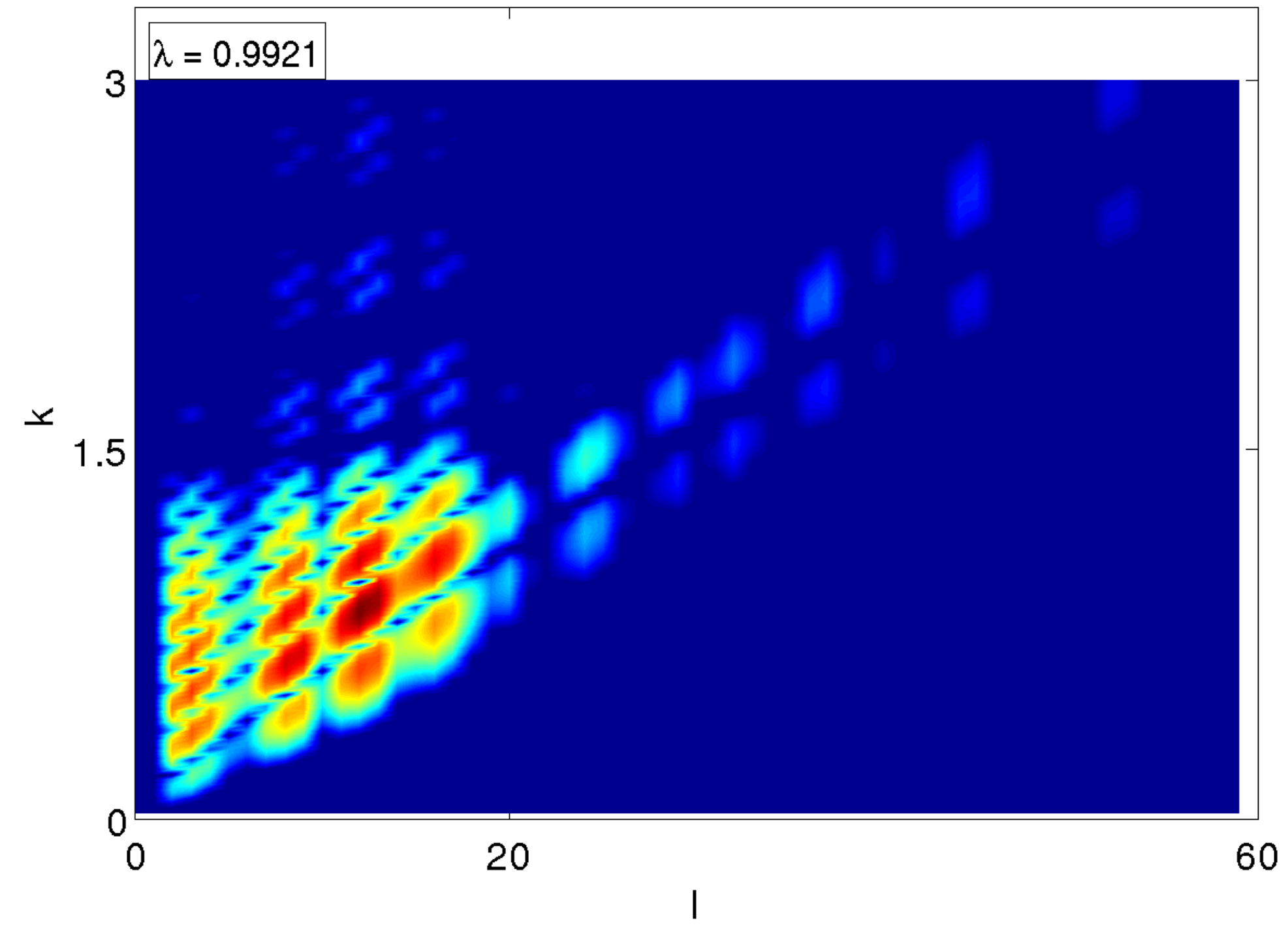}}

    \vspace{-3mm}
    \hspace{-15mm}
    \subfloat{
        \includegraphics[scale=0.22]{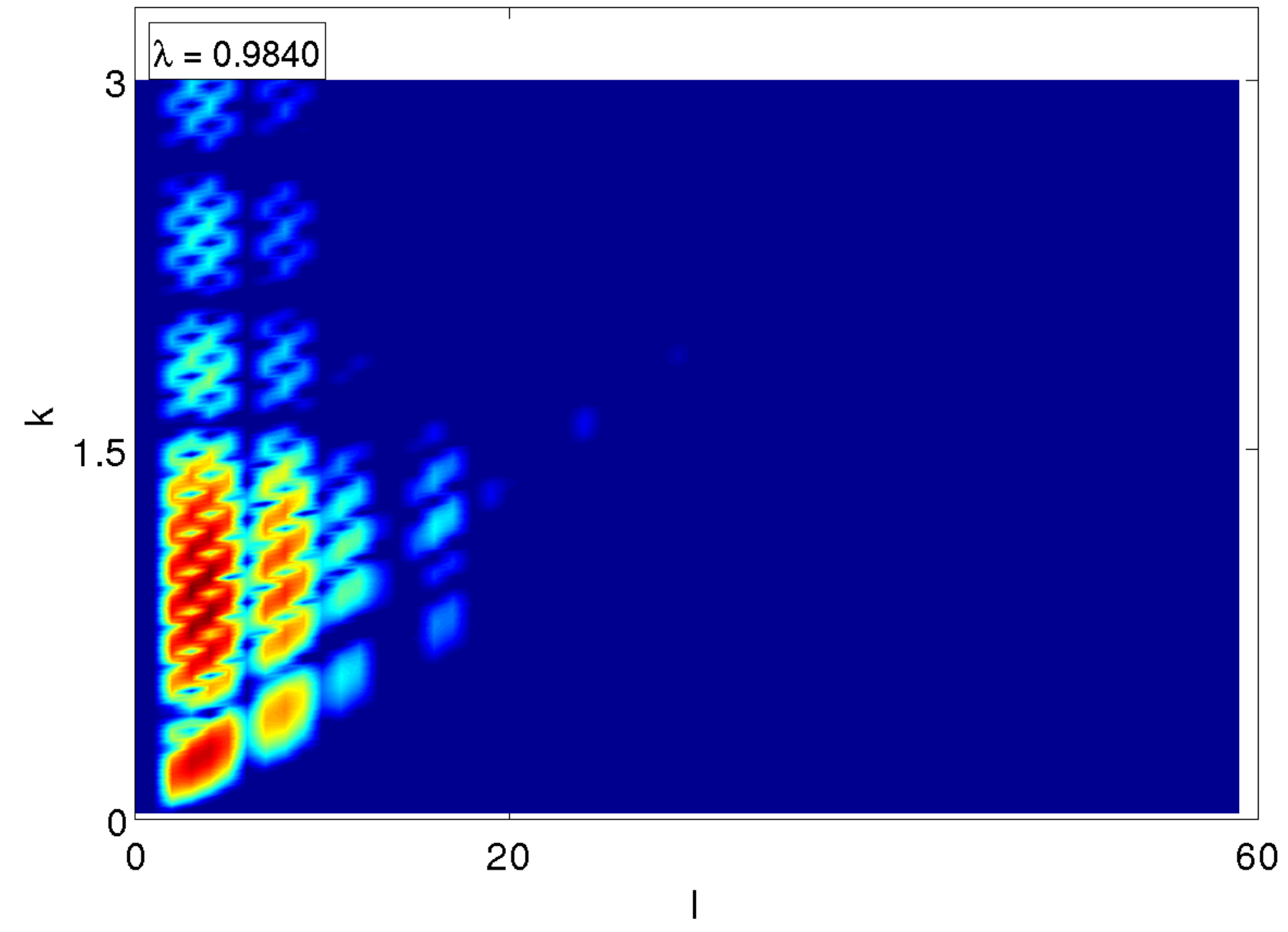}}
    \hspace{-2mm}
    \subfloat{
        \includegraphics[scale=0.22]{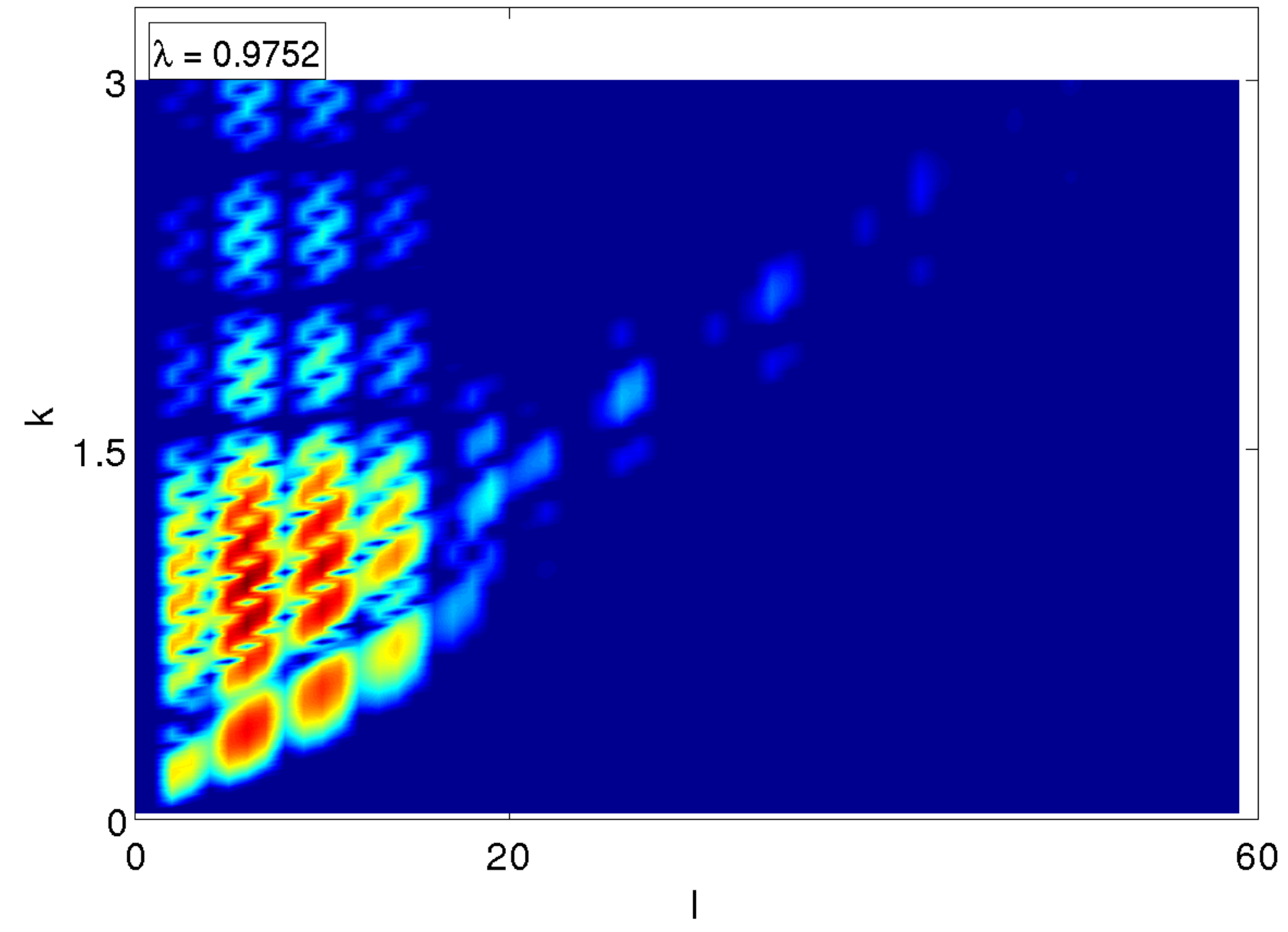}}
    \hspace{-2mm}
    \subfloat{
        \includegraphics[scale=0.22]{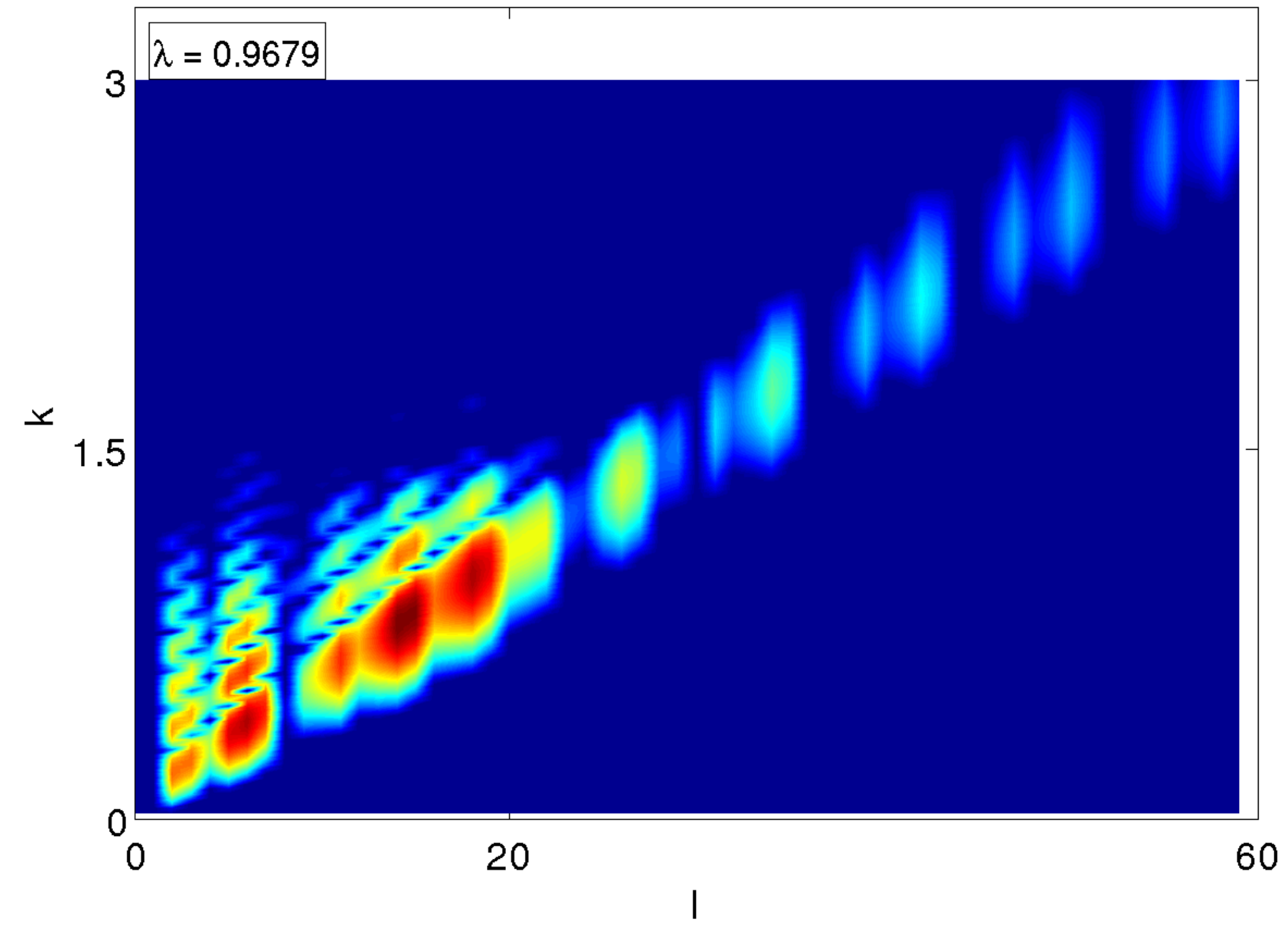}}

    \vspace{-3mm}
    \hspace{-15mm}
    \subfloat{
        \includegraphics[scale=0.22]{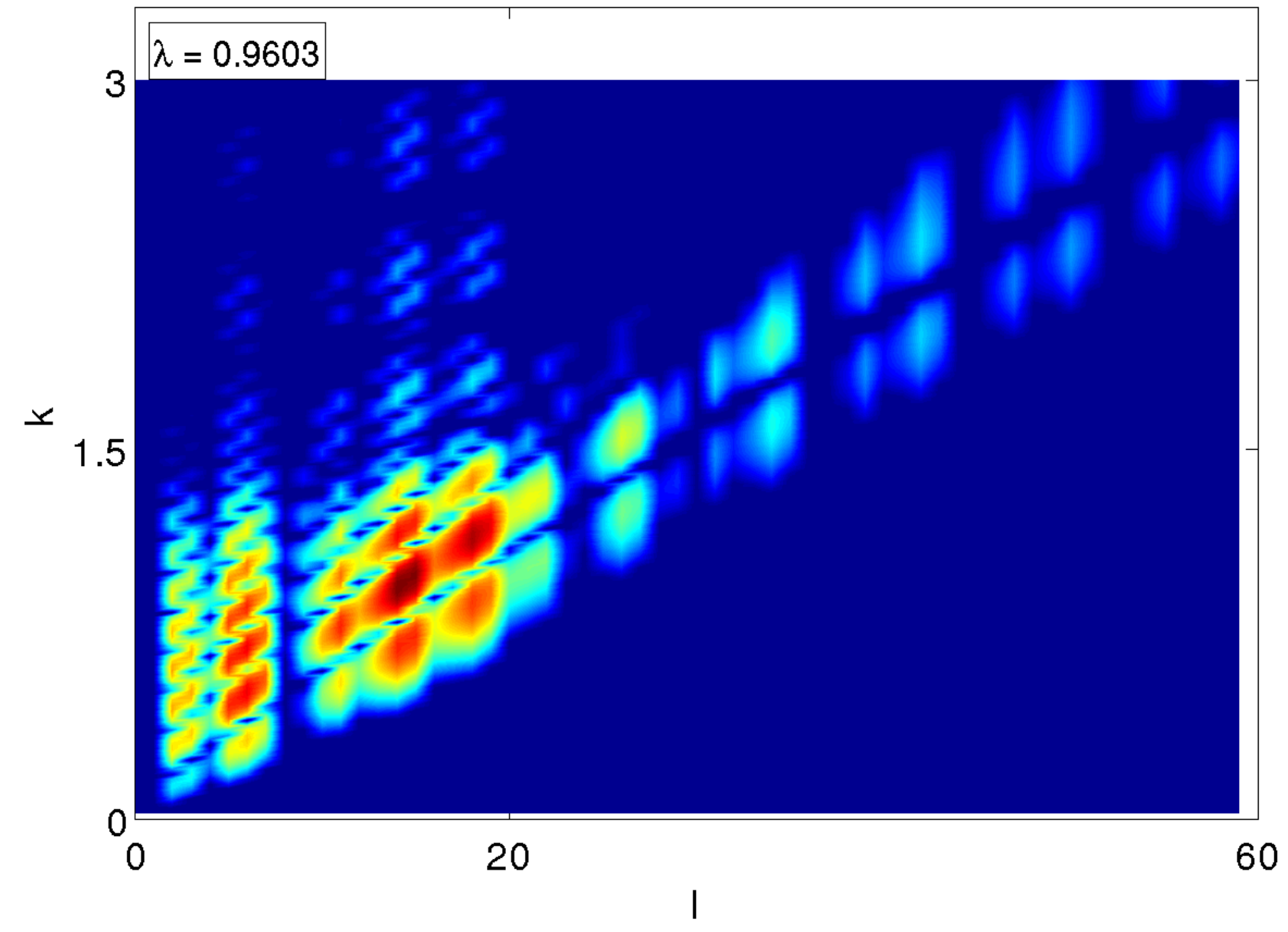}}
    \hspace{-2mm}
    \subfloat{
        \includegraphics[scale=0.22]{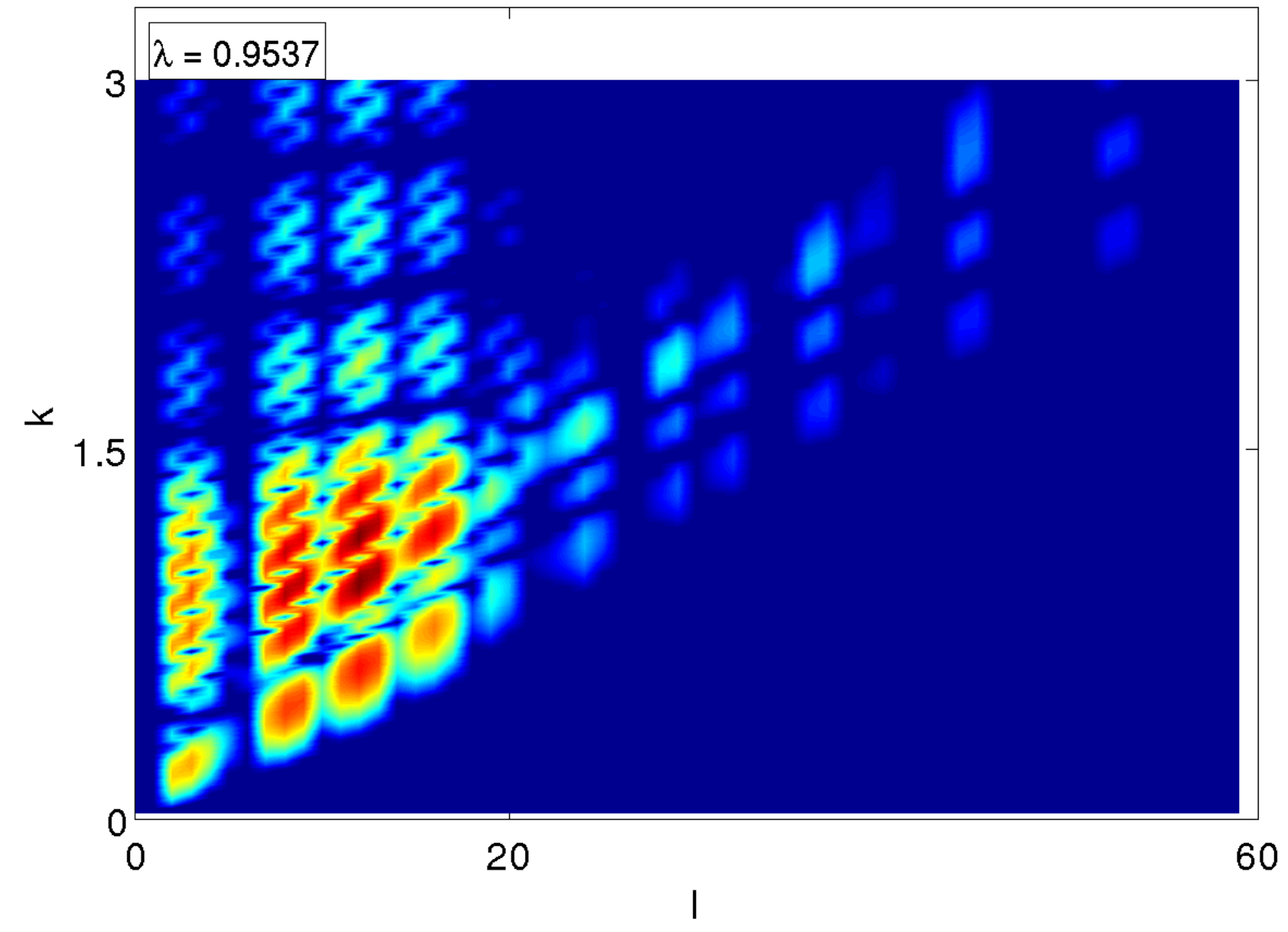}}
    \hspace{-2mm}
    \subfloat{
        \includegraphics[scale=0.22]{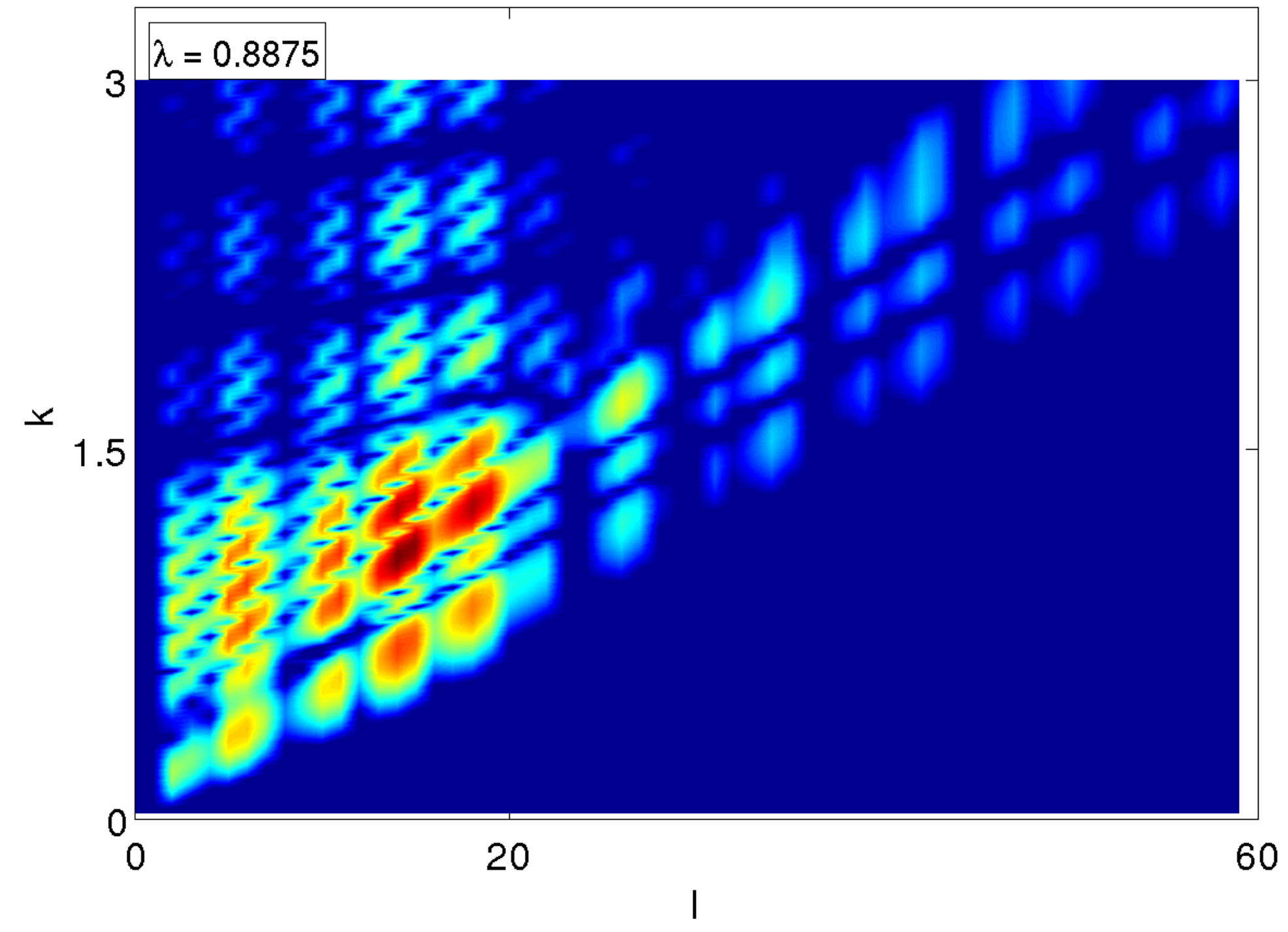}}


    \vspace{-3mm}

    \subfloat{
        \includegraphics[scale=0.28]{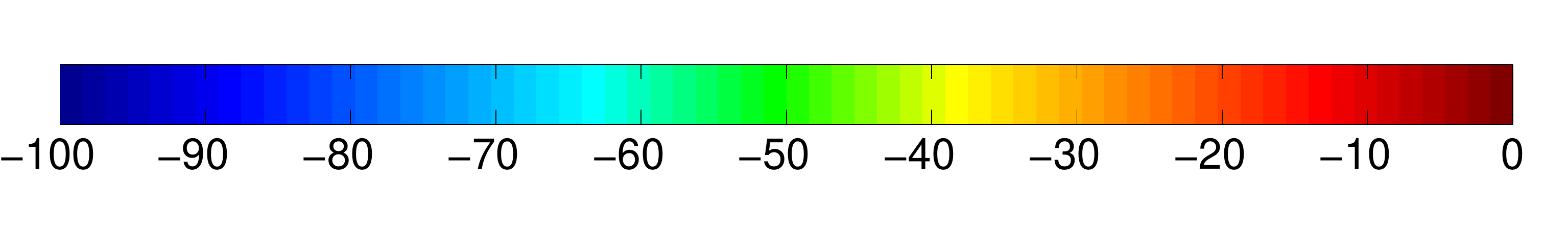}}

    \vspace{-3mm}

    \caption{\textbf{Fourier-Bessel spectral domain response of
        space-limited spectrally concentrated eigenfunctions}, given
      by their magnitude square $|g_{\ell m}(k)|^2$.  Each
      eigenfunction is limited in the azimuthally symmetric and
      radially independent spatial region $R=\{ 15\leq r \leq 25,\,
      \pi/8 \leq \theta \leq 3\pi/8, \, 0 \leq \phi <2\pi \}$ and
      spectrally concentrated in the Fourier-Bessel domain within the
      spectral region $\tilde A_{(1.4,20)}$. The eigenvalue $\lambda$
      associated with each eigenfunction is a measure of the spectral
      concentration within the Fourier-Bessel spectral domain
      $\tilde A_{(1.4,20)}$. The values of $|g_{\ell m}(k)|^2$ are
      plotted in decibels as $20\log|g_{\ell m}(k)|$, normalized to
      zero at the individual maxima of each
      eigenfunction. } \label{fig:eigen_bessel_spect}
\end{figure*}
%
%

\begin{figure*}[!htb]
    \centering
    \vspace{-10mm}
    \hspace{-14mm}
    \subfloat{
        \includegraphics[scale=0.181]{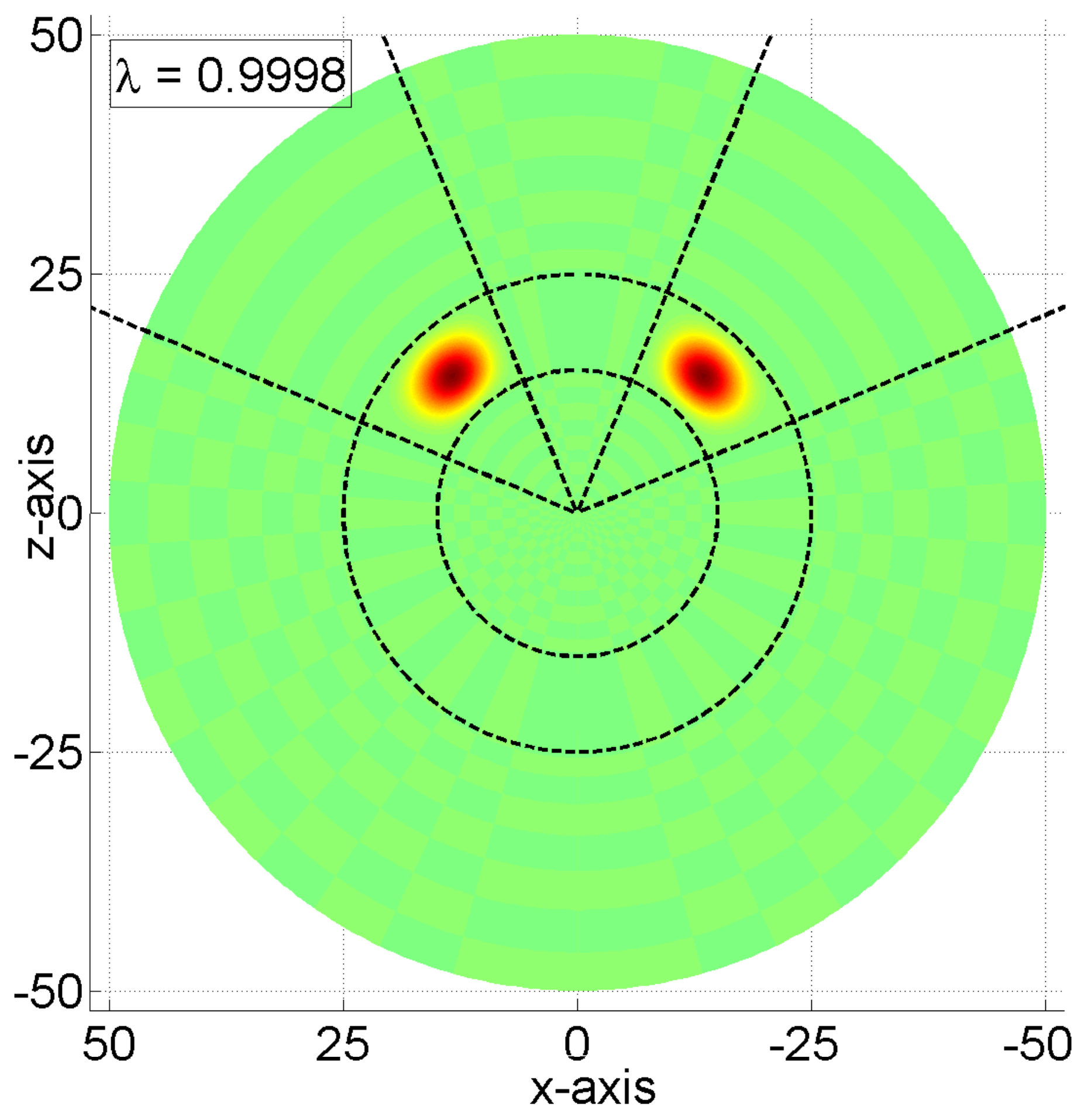}}
    \hspace{3mm}
    \subfloat{
        \includegraphics[scale=0.181]{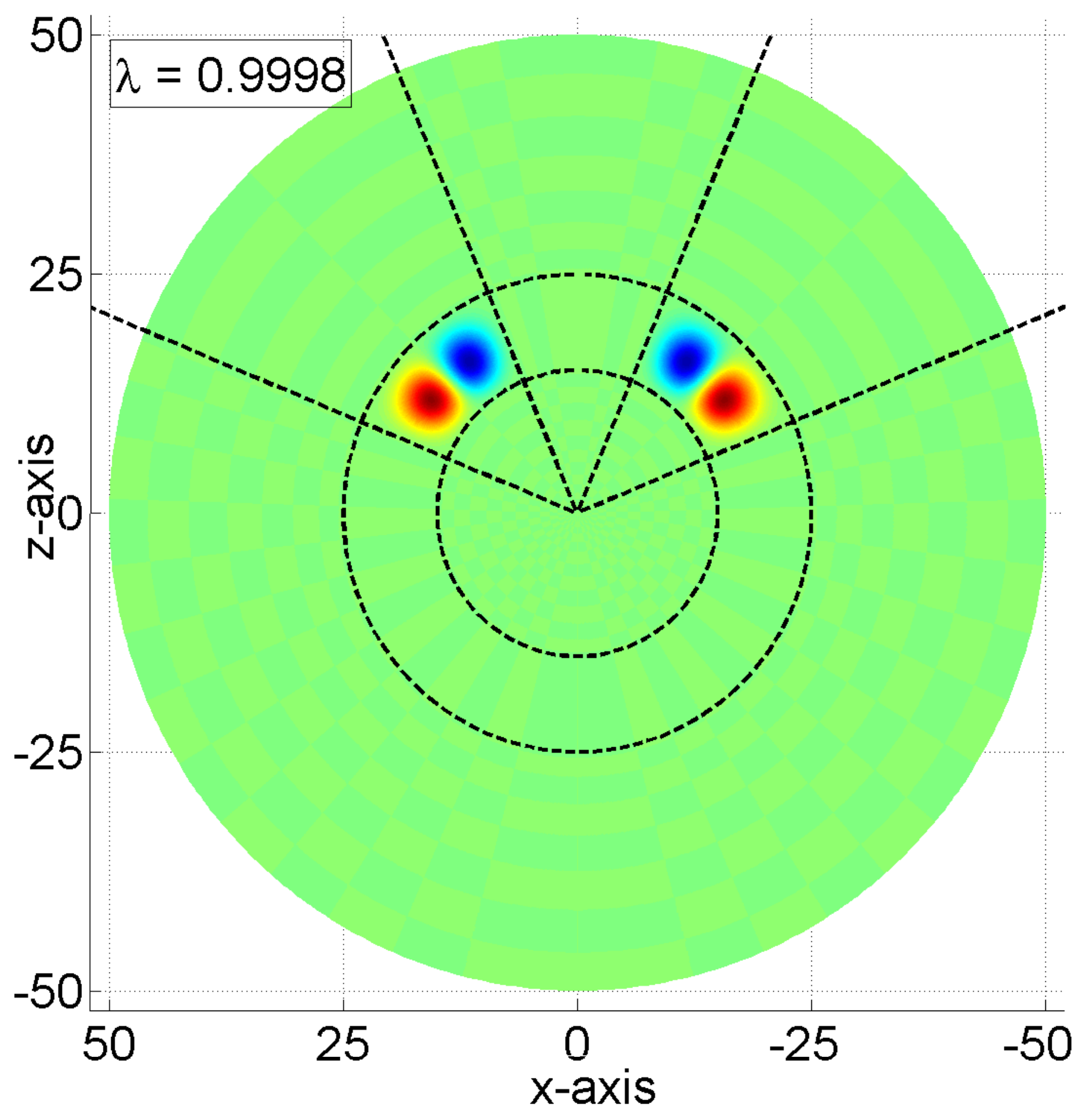}}
    \hspace{3mm}
    \subfloat{
        \includegraphics[scale=0.181]{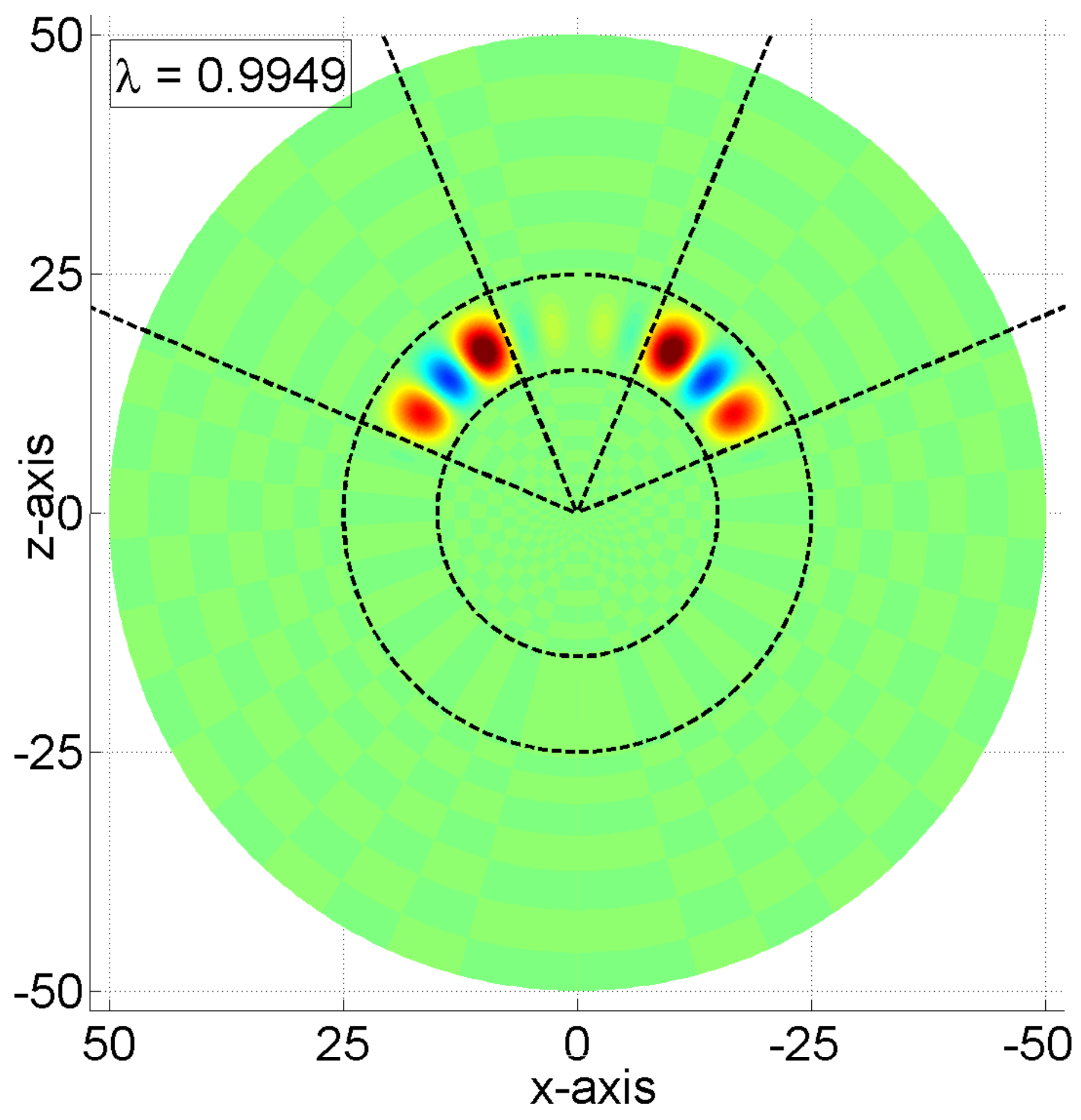}}

    \vspace{-3mm}
    \hspace{-14mm}
    \subfloat{
        \includegraphics[scale=0.181]{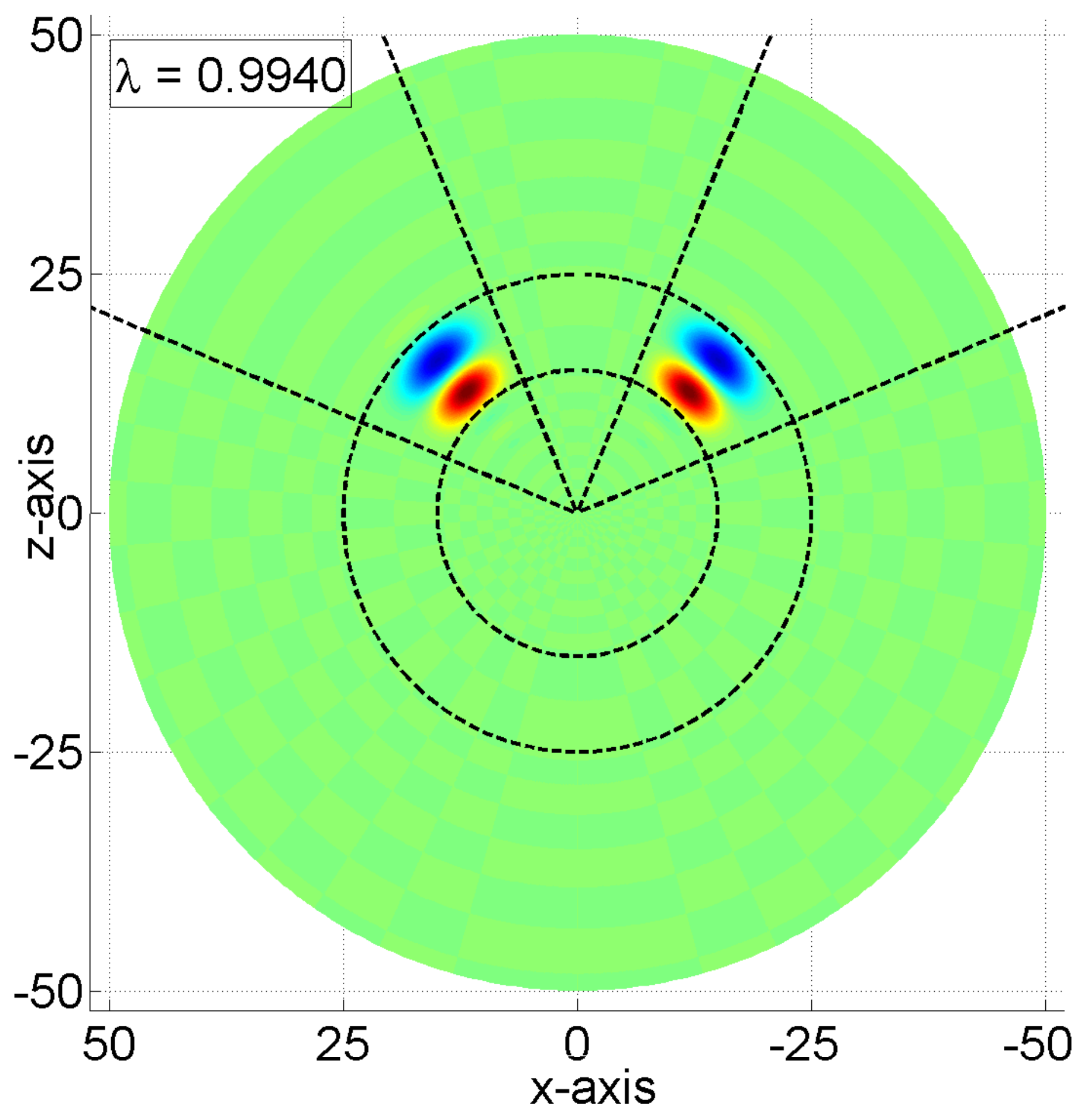}}
    \hspace{3mm}
    \subfloat{
        \includegraphics[scale=0.181]{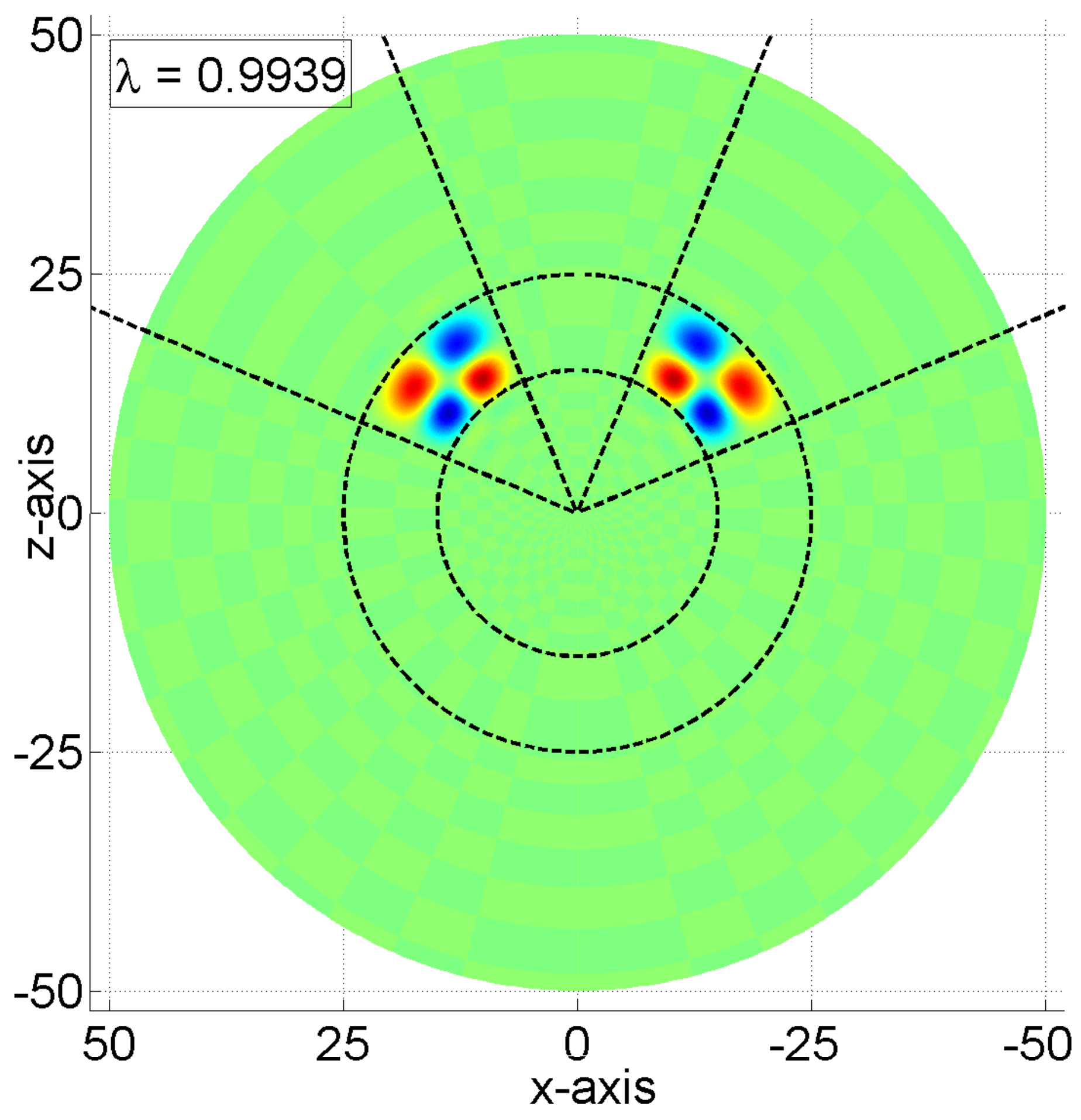}}
    \hspace{3mm}
    \subfloat{
        \includegraphics[scale=0.181]{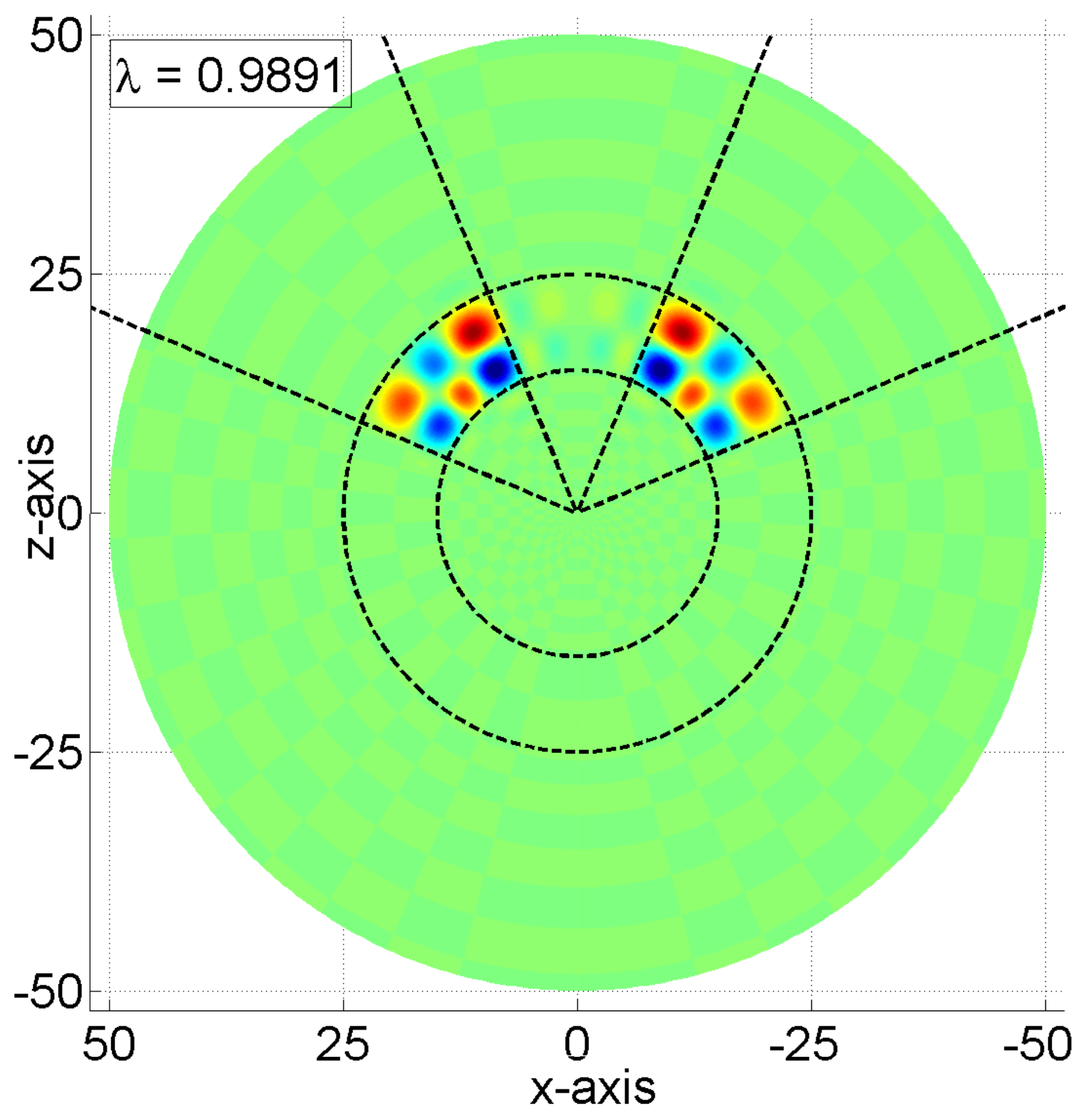}}

    \vspace{-3mm}
    \hspace{-14mm}
    \subfloat{
        \includegraphics[scale=0.181]{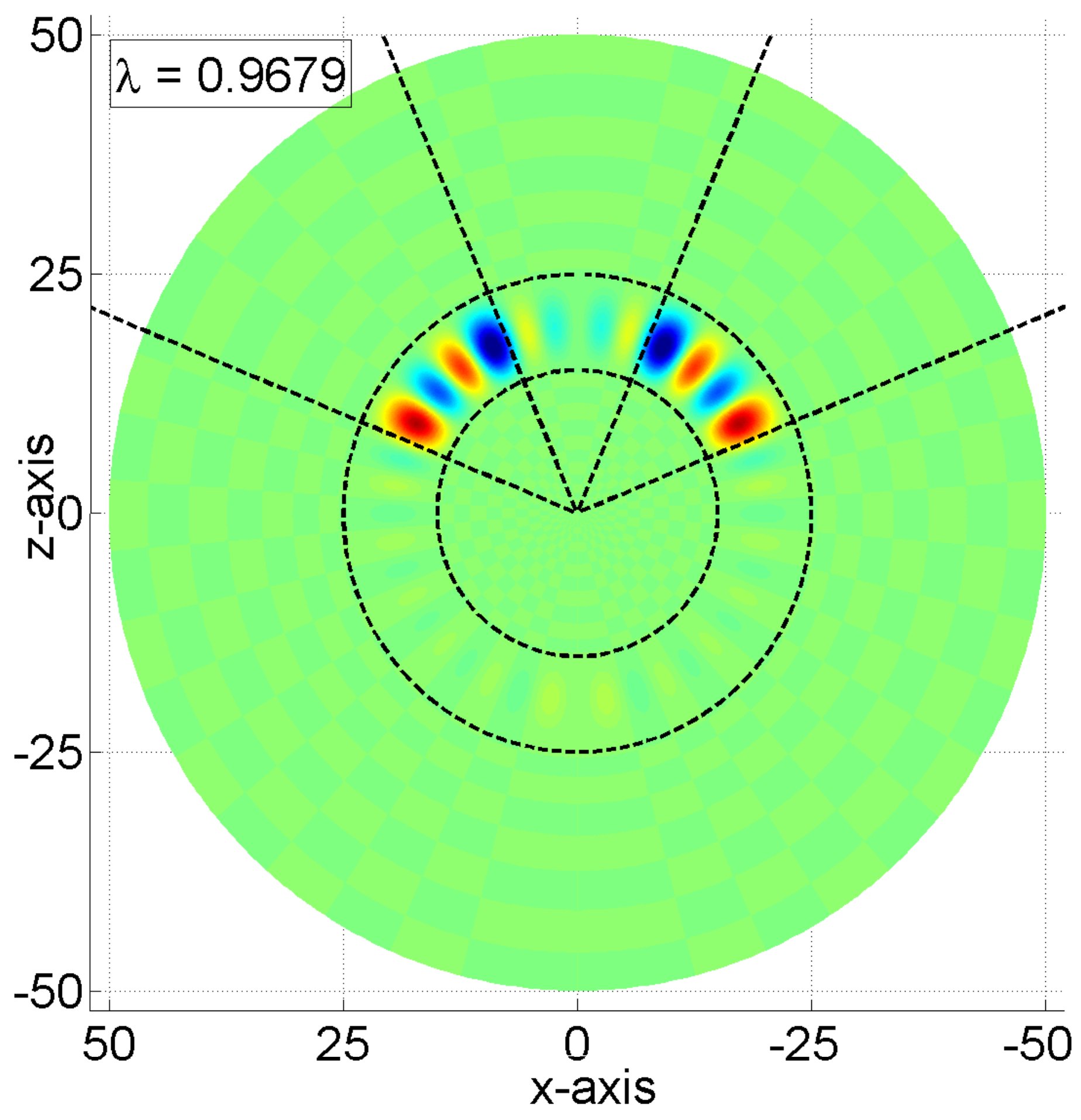}}
    \hspace{3mm}
    \subfloat{
        \includegraphics[scale=0.181]{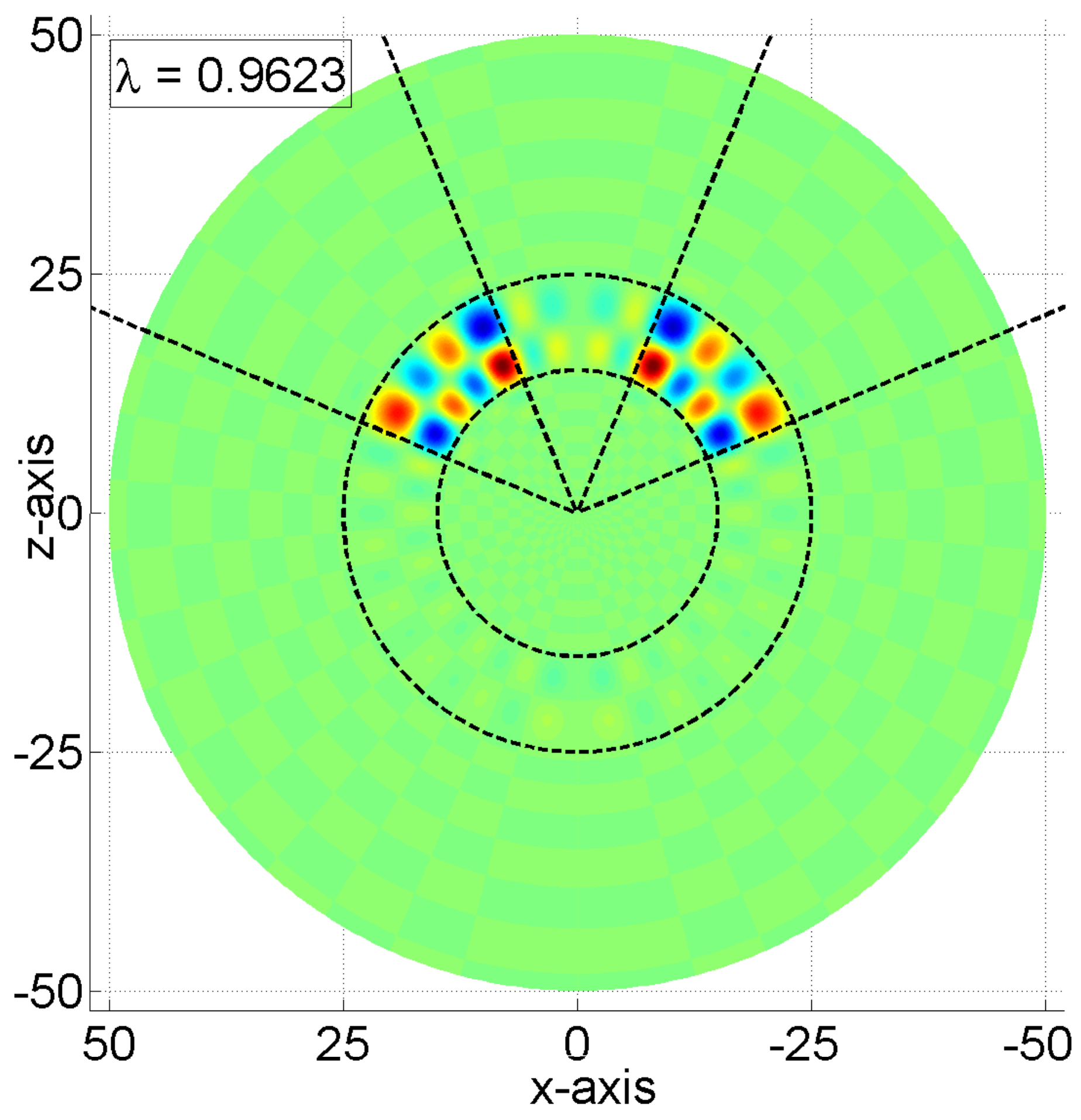}}
    \hspace{3mm}
    \subfloat{
        \includegraphics[scale=0.181]{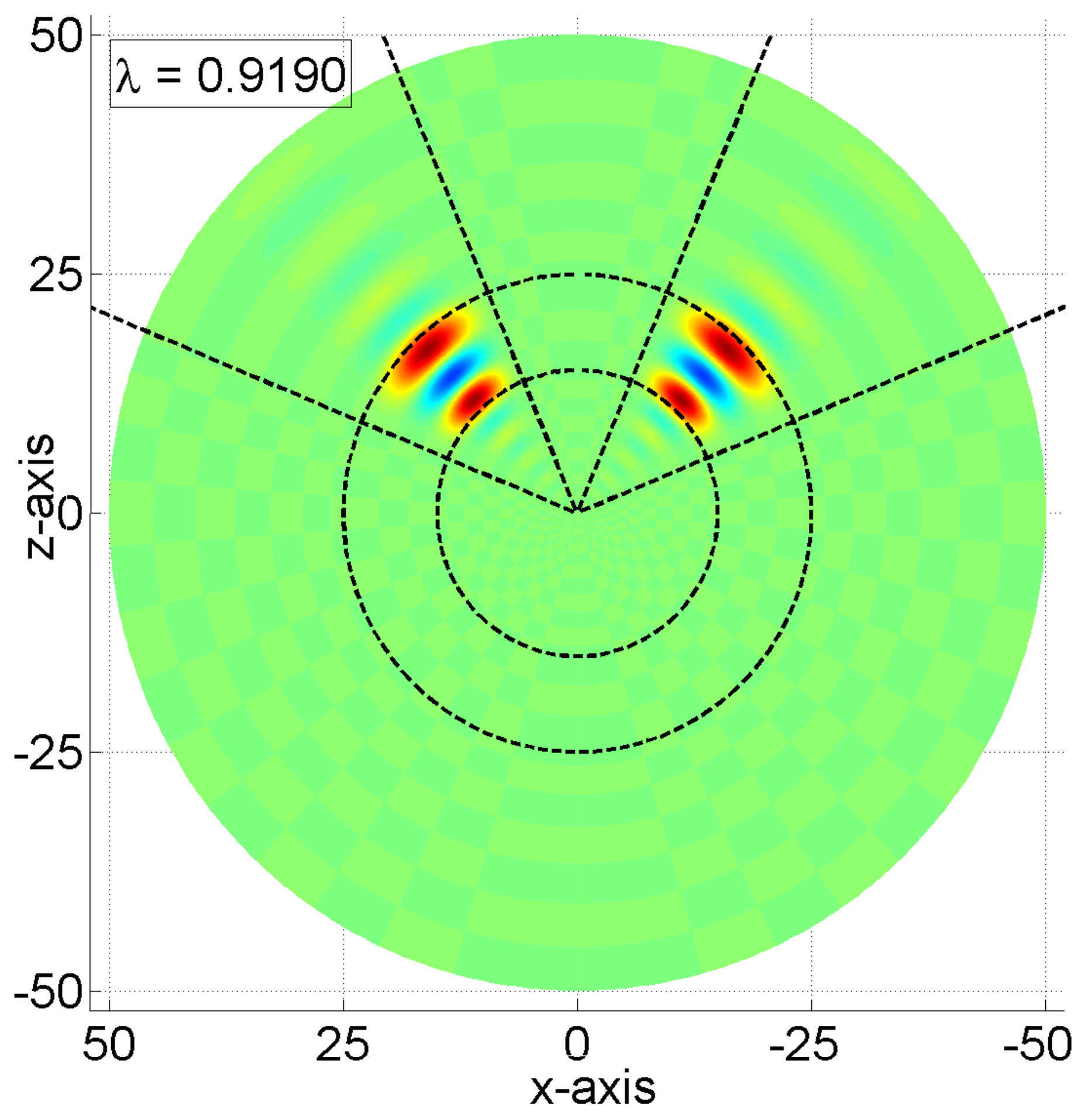}}

    \vspace{-3mm}
    \hspace{-14mm}
    \subfloat{
        \includegraphics[scale=0.181]{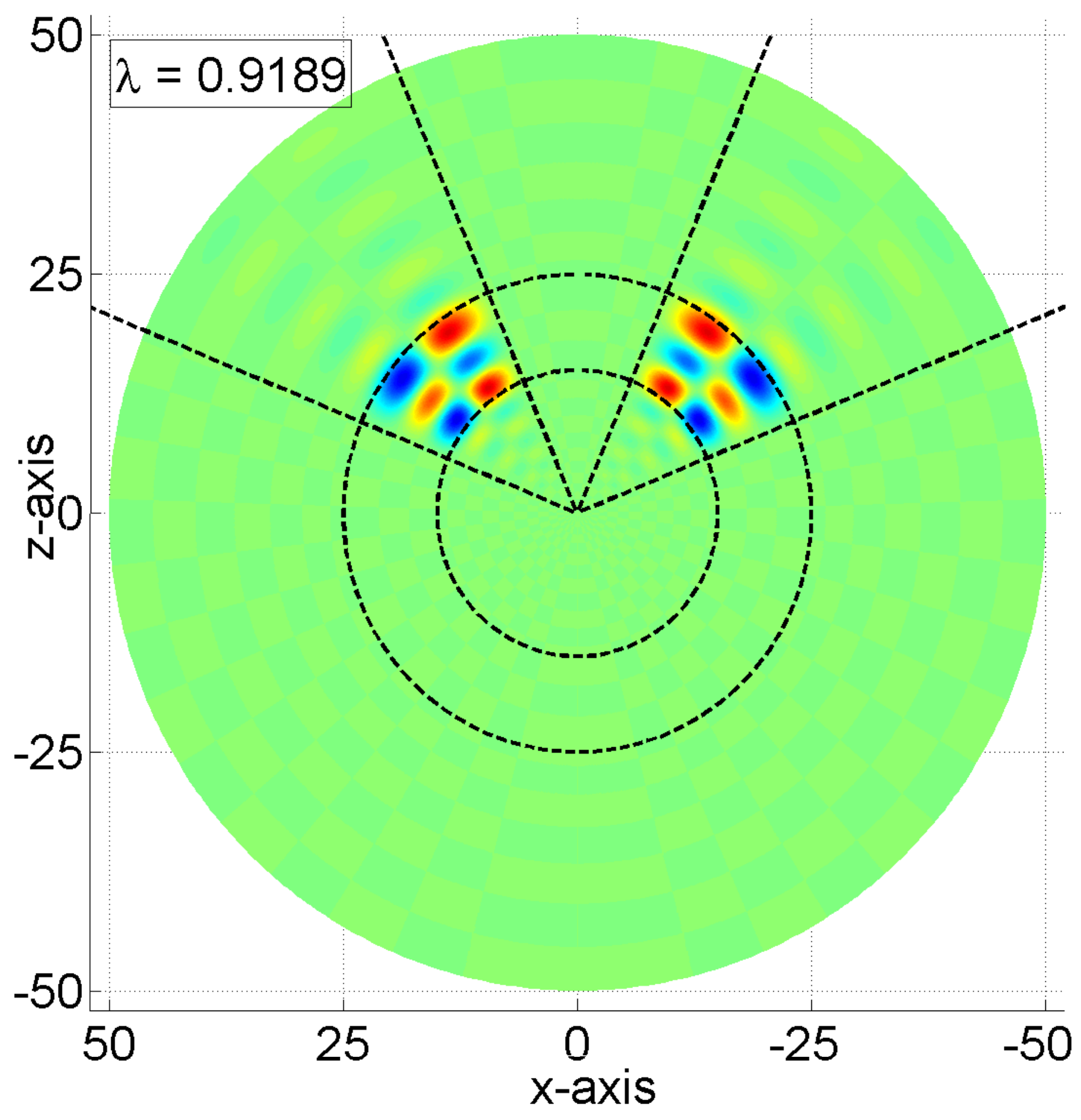}}
    \hspace{3mm}
    \subfloat{
        \includegraphics[scale=0.181]{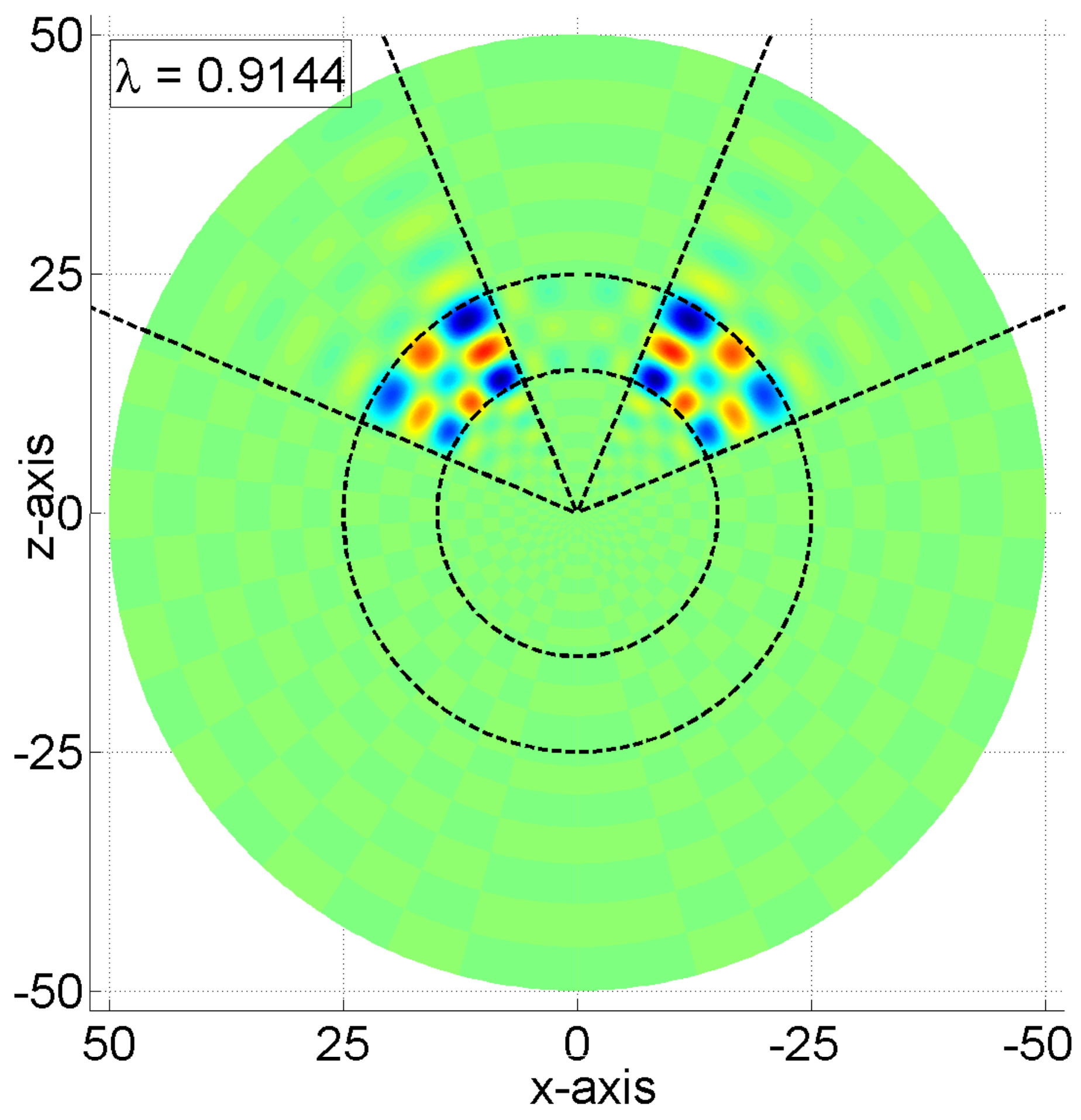}}
    \hspace{3mm}
    \subfloat{
        \includegraphics[scale=0.181]{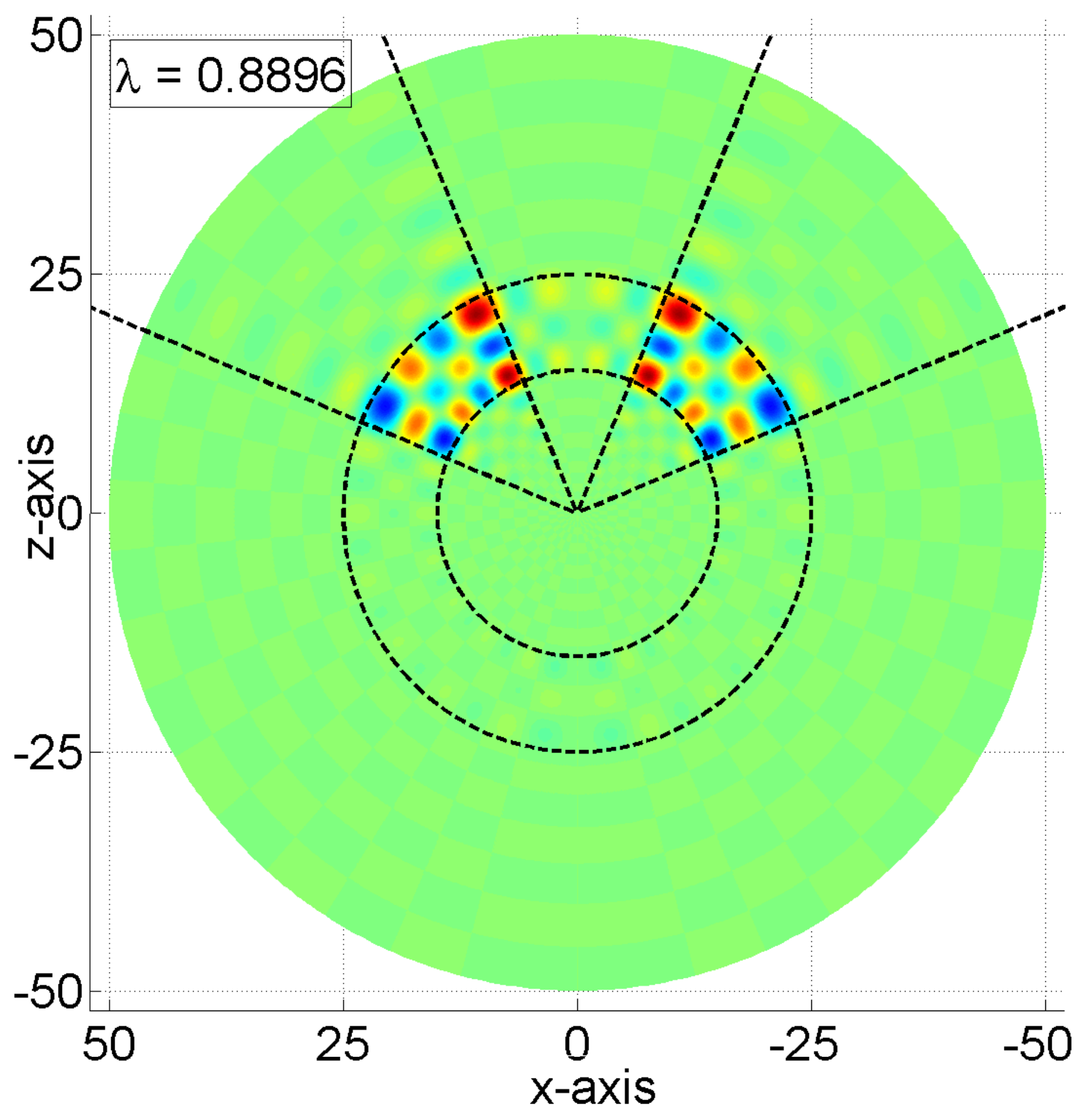}}


    \vspace{-3mm}

    \subfloat{
        \includegraphics[scale=0.25]{colorbar_bessel_spatial}}
    \vspace{-3mm}

    \caption{\textbf{Fourier-Laguerre band-limited spatially
        concentrated eigenfunctions} $f^m(r,\theta)$ given in
      \eqref{Eq:spatial_eigen_subproblem2}, obtained as solutions of
      the fixed-order eigenvalue problem in
      \eqref{Eq:symmetric1_sub_problem2} for $m=2$.  Each
      eigenfunction $f^m(r,\theta)$ is independent of $\phi$ and is
      plotted for $r\leq50$.  The spatial region of concentration is
      $R=\{ 15 \leq r \leq 25,\, \pi/8 \leq \theta \leq 3\pi/8,$
      \mbox{$0 \leq \phi <2 \pi \}$} and is azimuthally symmetric and
      radially independent.  The dependence of the eigenfunction in
      $\phi$ is $e^{im\phi}$, as given in
      \eqref{Eq:rtheta_to_ball}. The eigenfunctions are band-limited
      in the Fourier-Laguerre domain within the spectral region
      $A_{(30,20)}$. The eigenvalue $\lambda$ associated with each
      eigenfunction is a measure of spatial concentration within the
      spatial region $R$. } \label{fig:eigen_laguerre_spatial}
\end{figure*}
\begin{figure*}[!htb]
    \vspace{-12mm}
    \centering
    \hspace{-15mm}
    \subfloat{
        \includegraphics[scale=0.215]{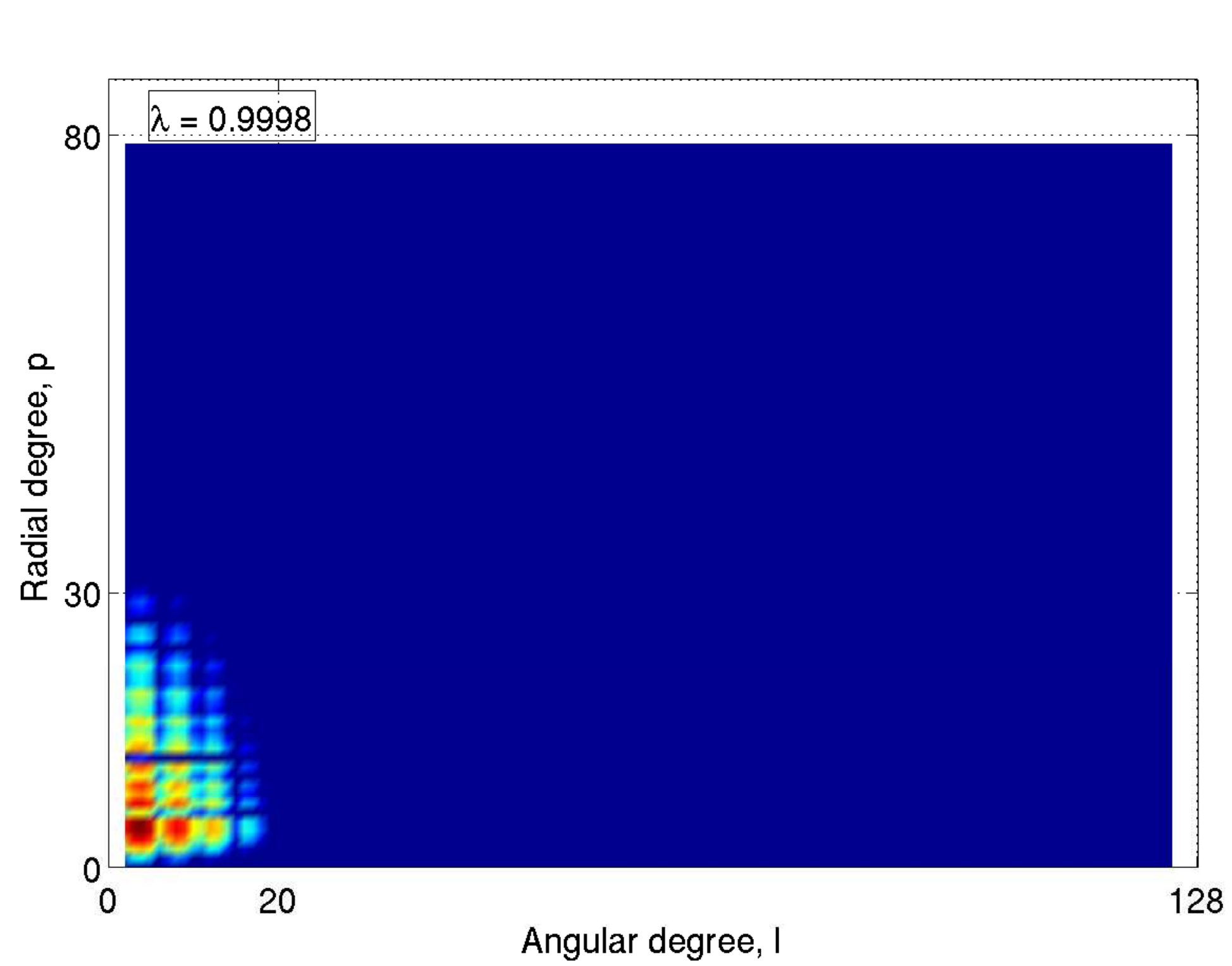}}
    \hspace{-2mm}
    \subfloat{
        \includegraphics[scale=0.215]{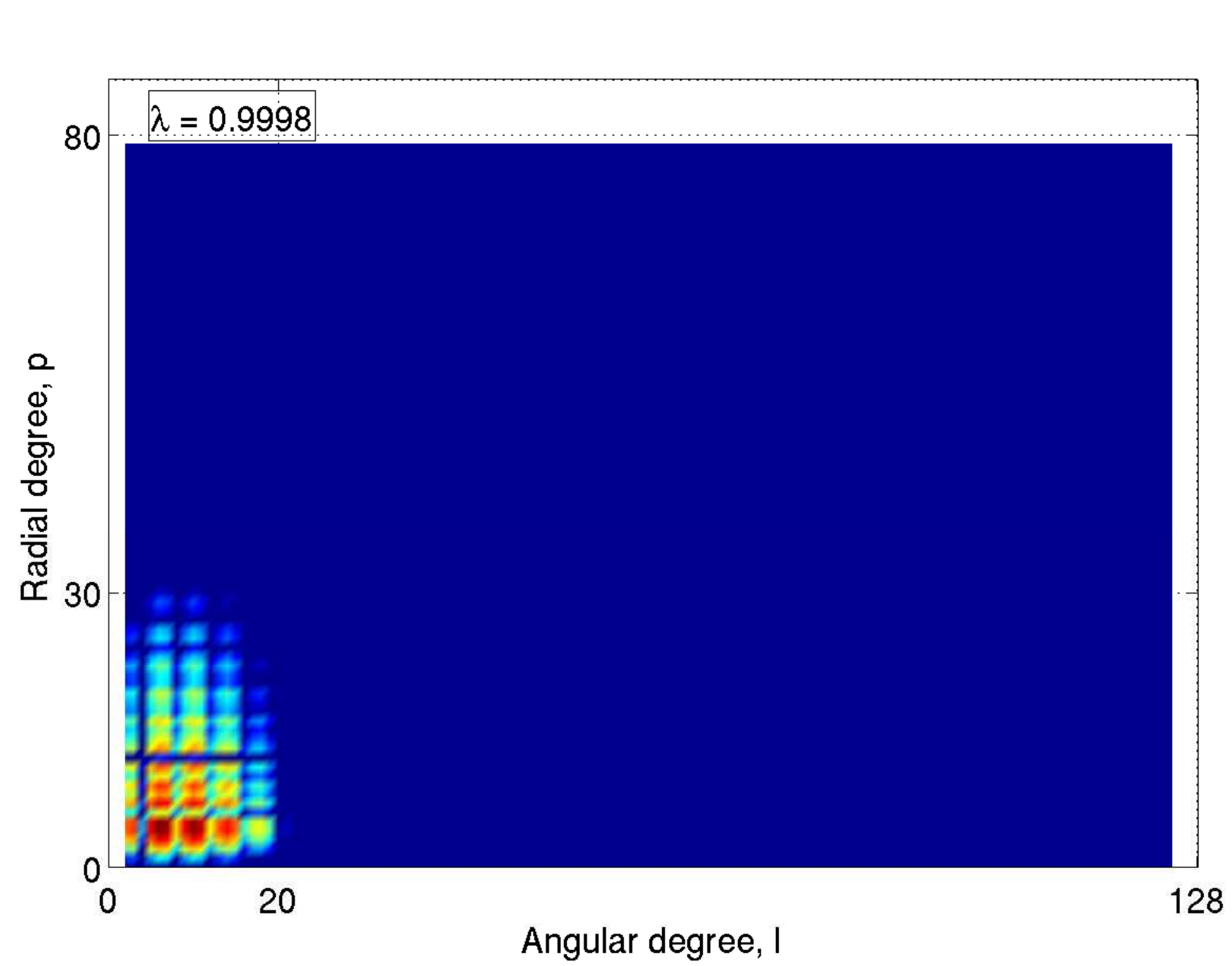}}
    \hspace{-2mm}
    \subfloat{
        \includegraphics[scale=0.215]{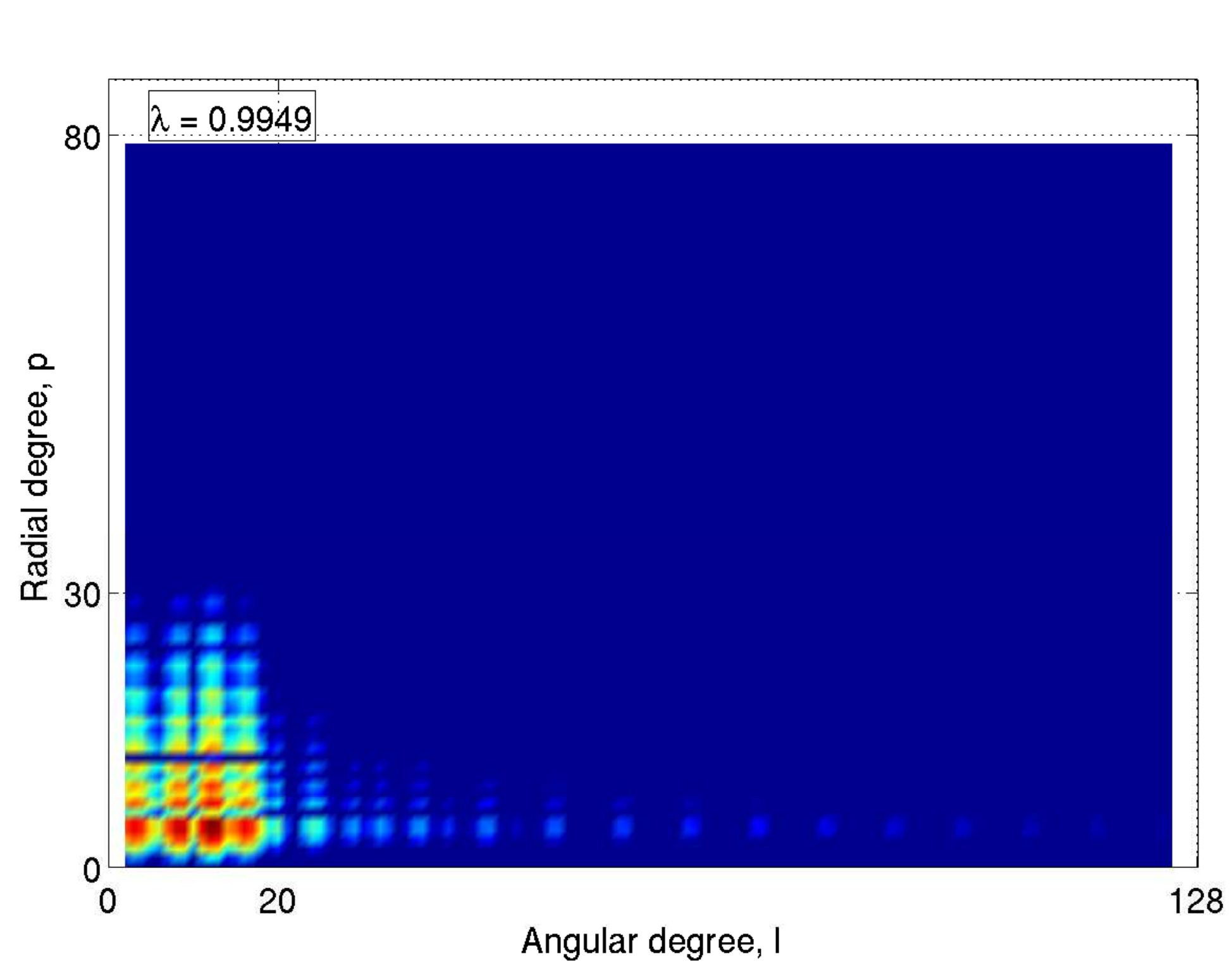}}

    \vspace{-3mm}
    \hspace{-15mm}
    \subfloat{
        \includegraphics[scale=0.215]{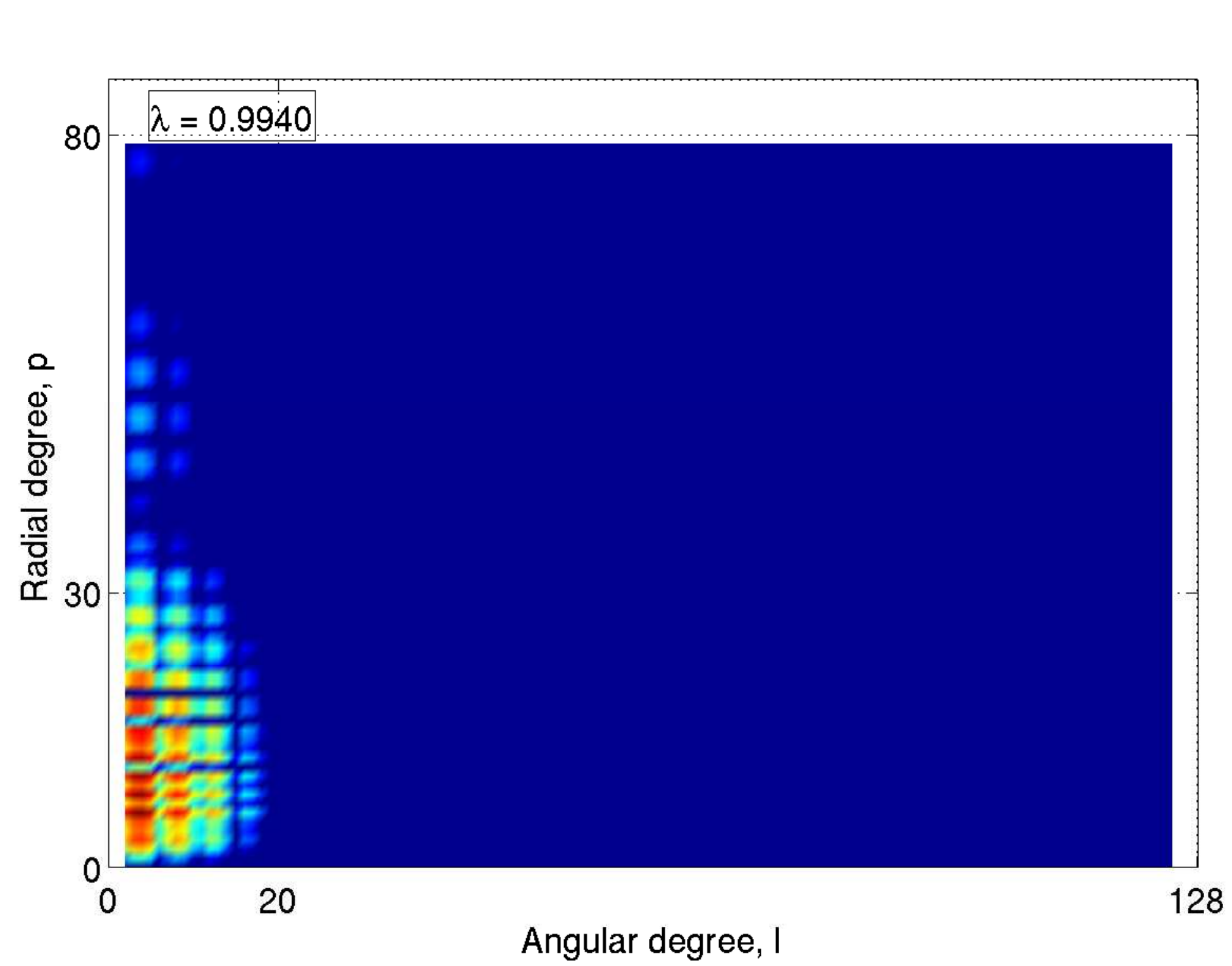}}
    \hspace{-2mm}
    \subfloat{
        \includegraphics[scale=0.215]{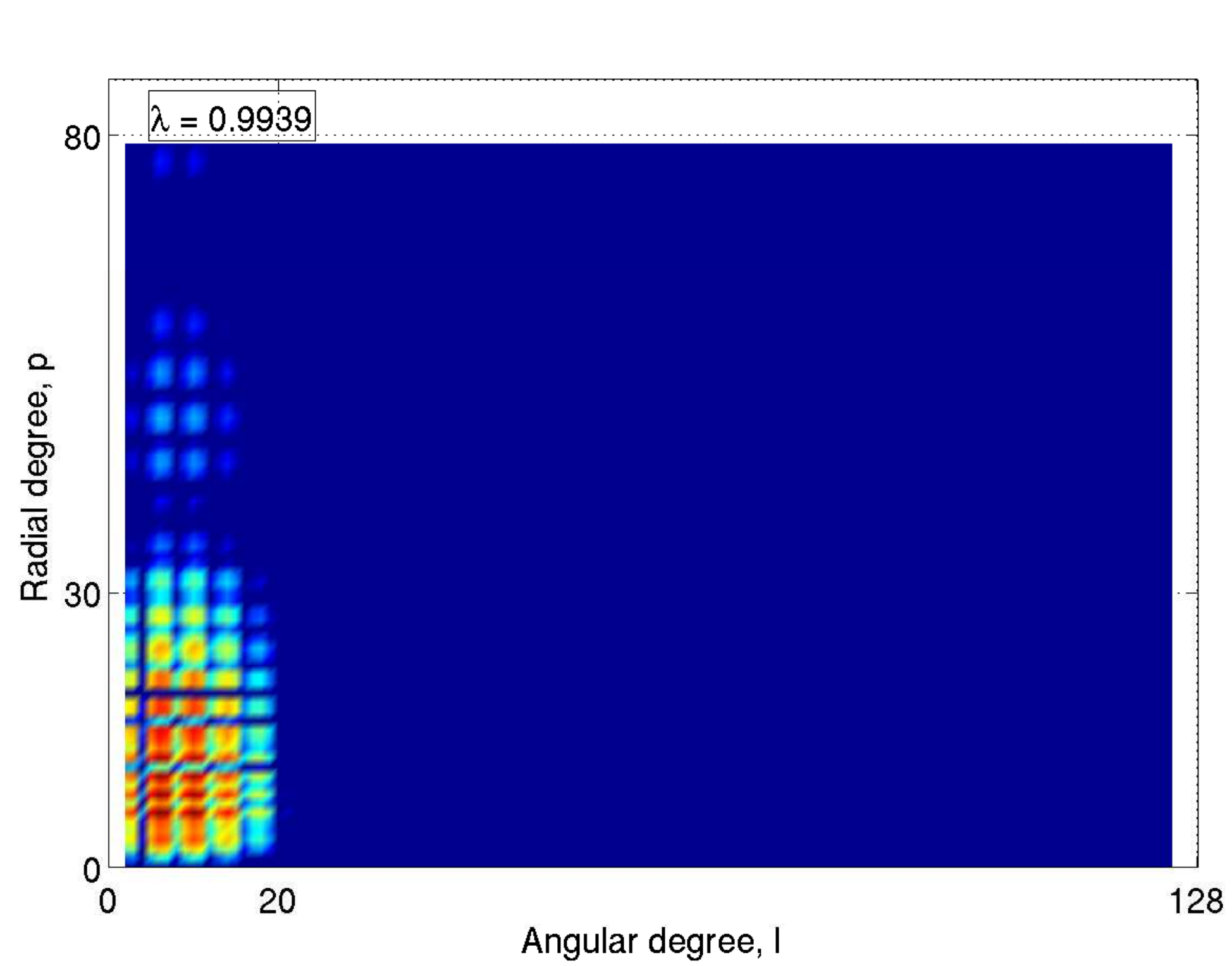}}
    \hspace{-2mm}
    \subfloat{
        \includegraphics[scale=0.215]{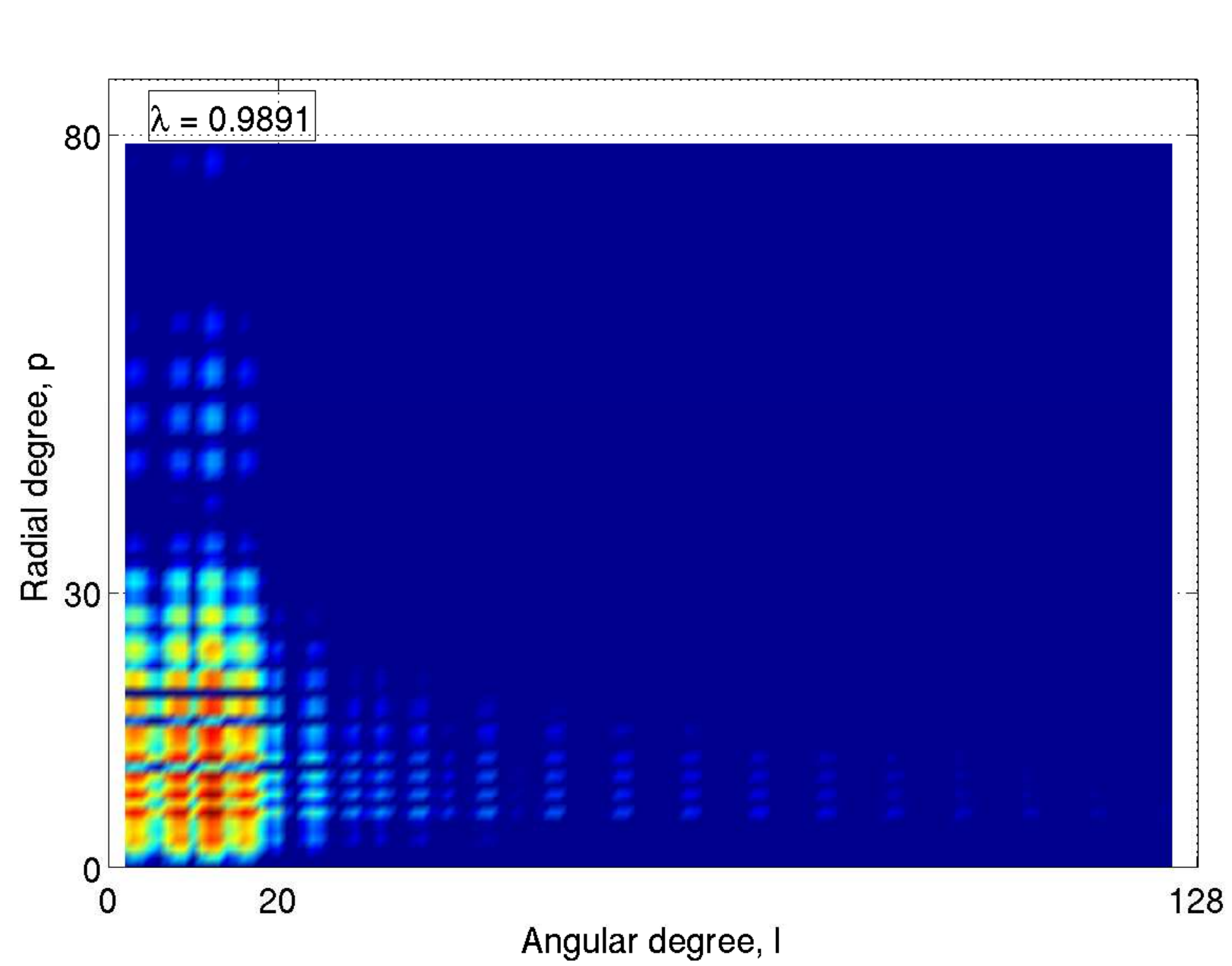}}

    \vspace{-3mm}
    \hspace{-15mm}
    \subfloat{
        \includegraphics[scale=0.215]{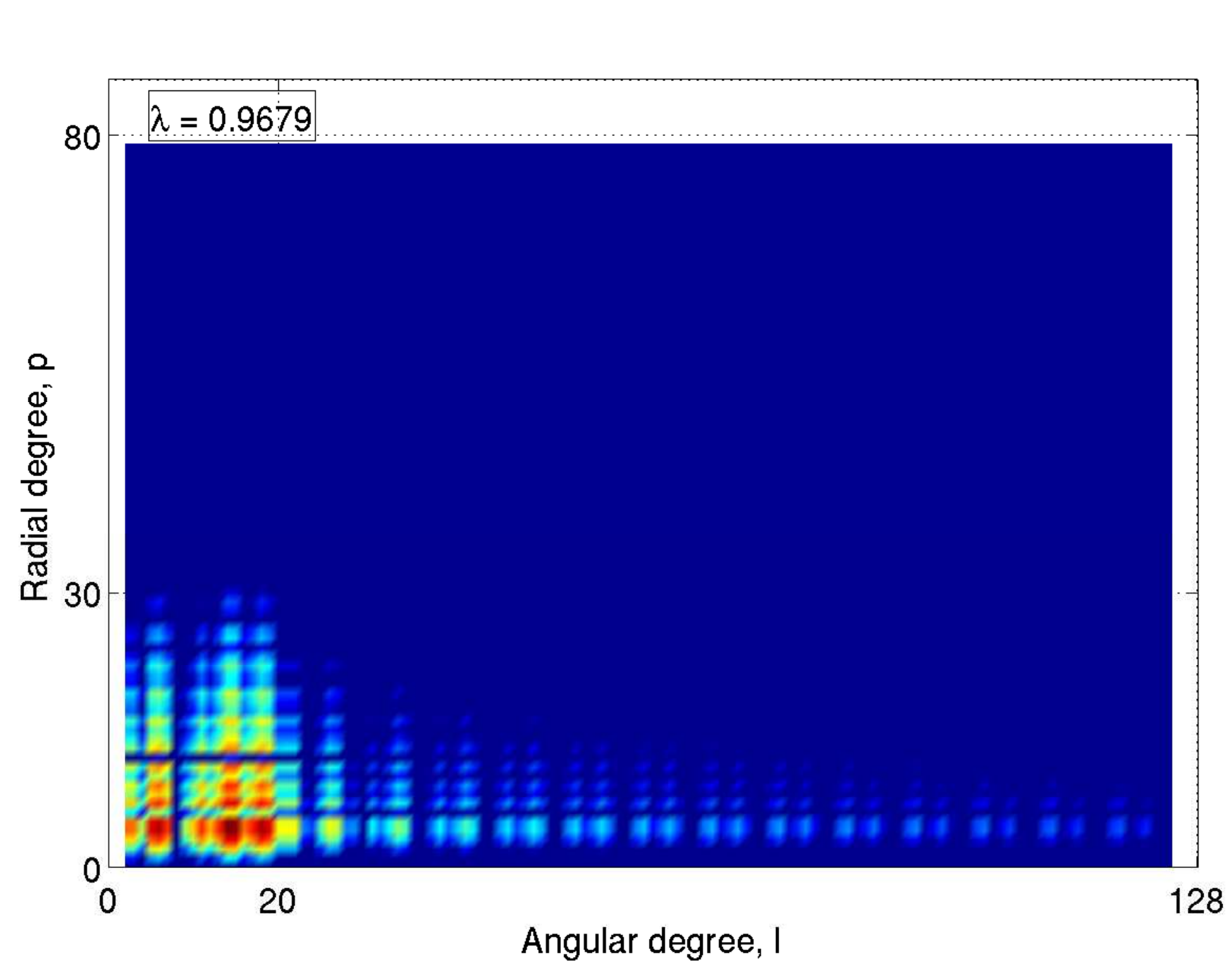}}
    \hspace{-2mm}
    \subfloat{
        \includegraphics[scale=0.215]{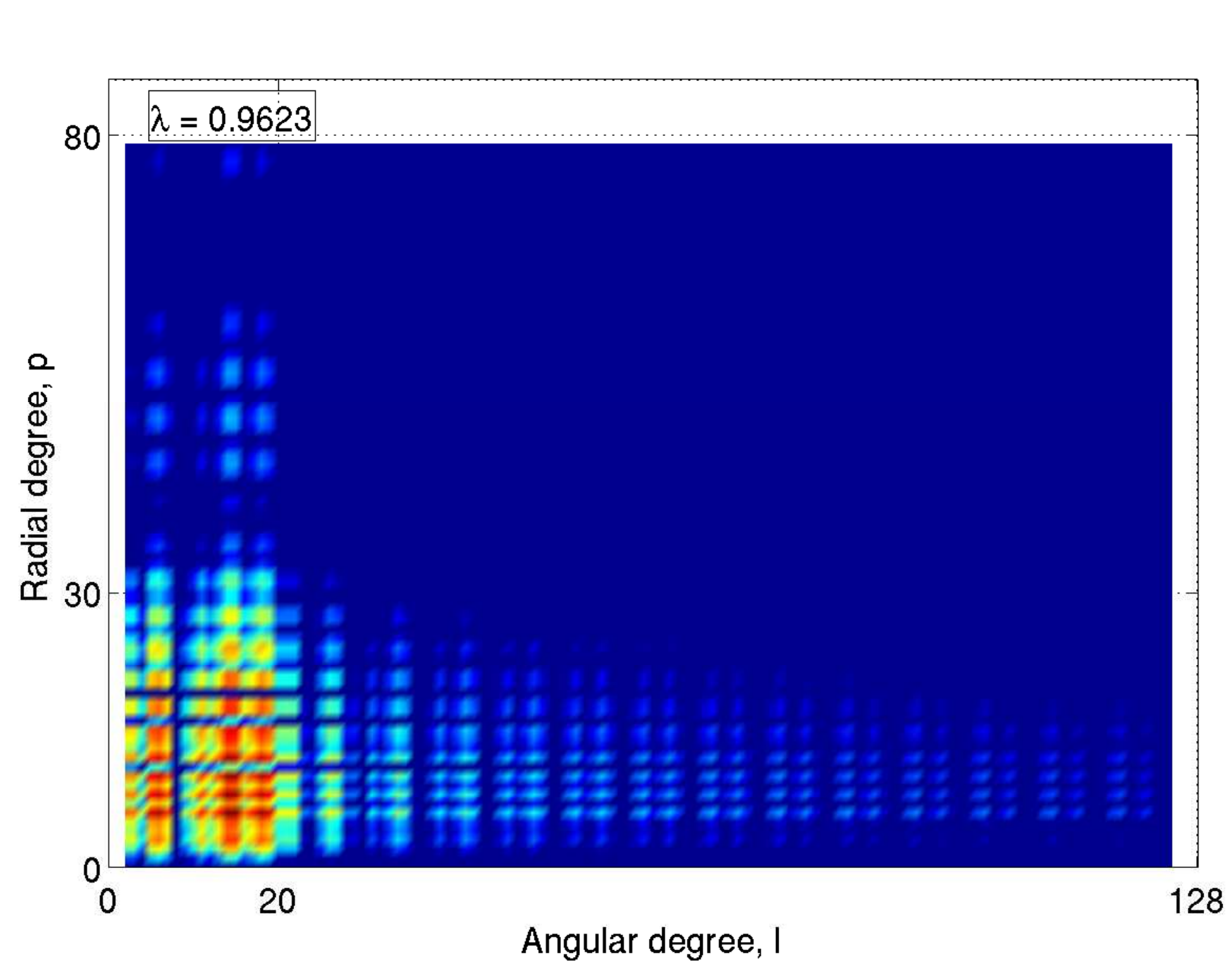}}
    \hspace{-2mm}
    \subfloat{
        \includegraphics[scale=0.215]{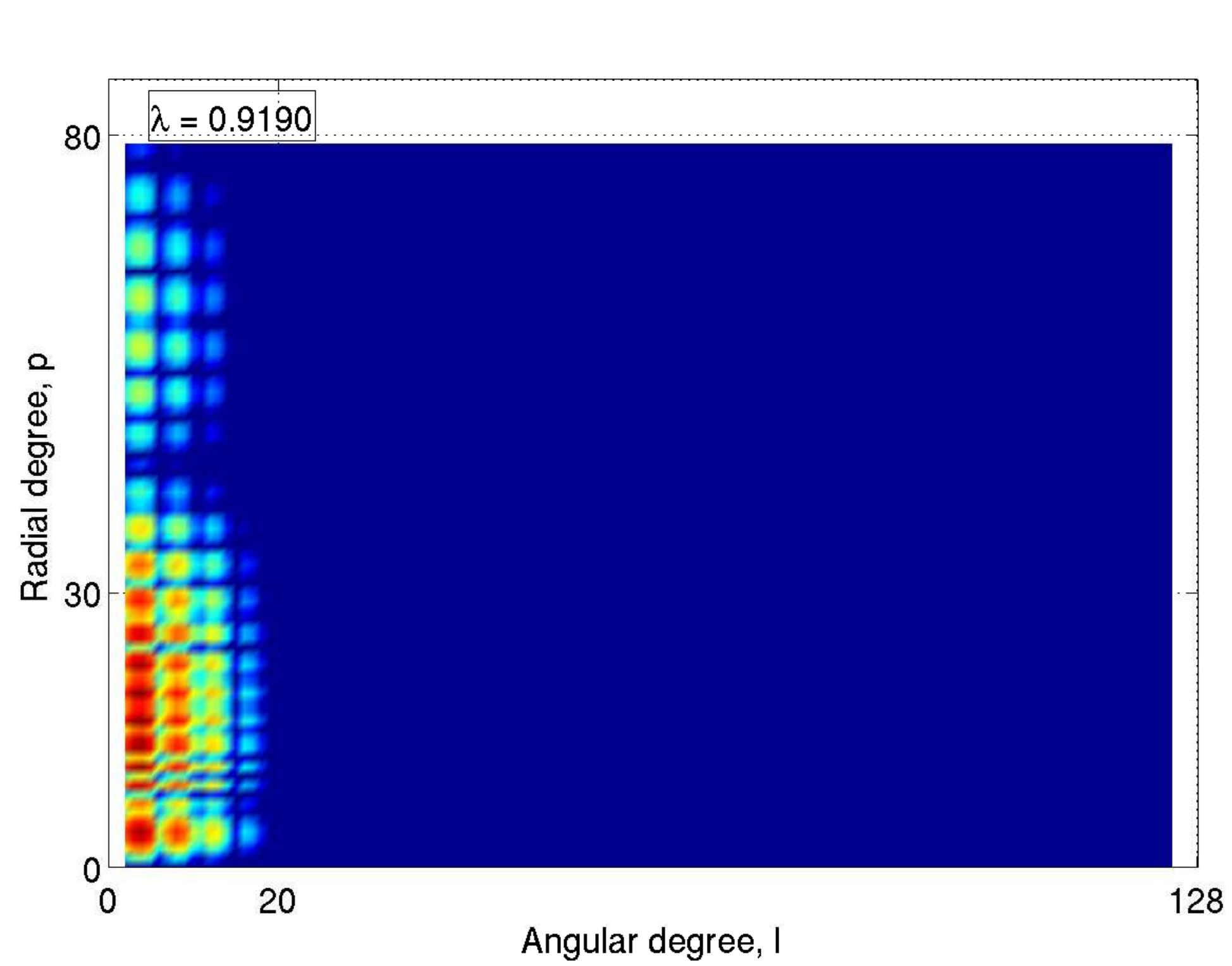}}

    \vspace{-3mm}
    \hspace{-15mm}
    \subfloat{
        \includegraphics[scale=0.215]{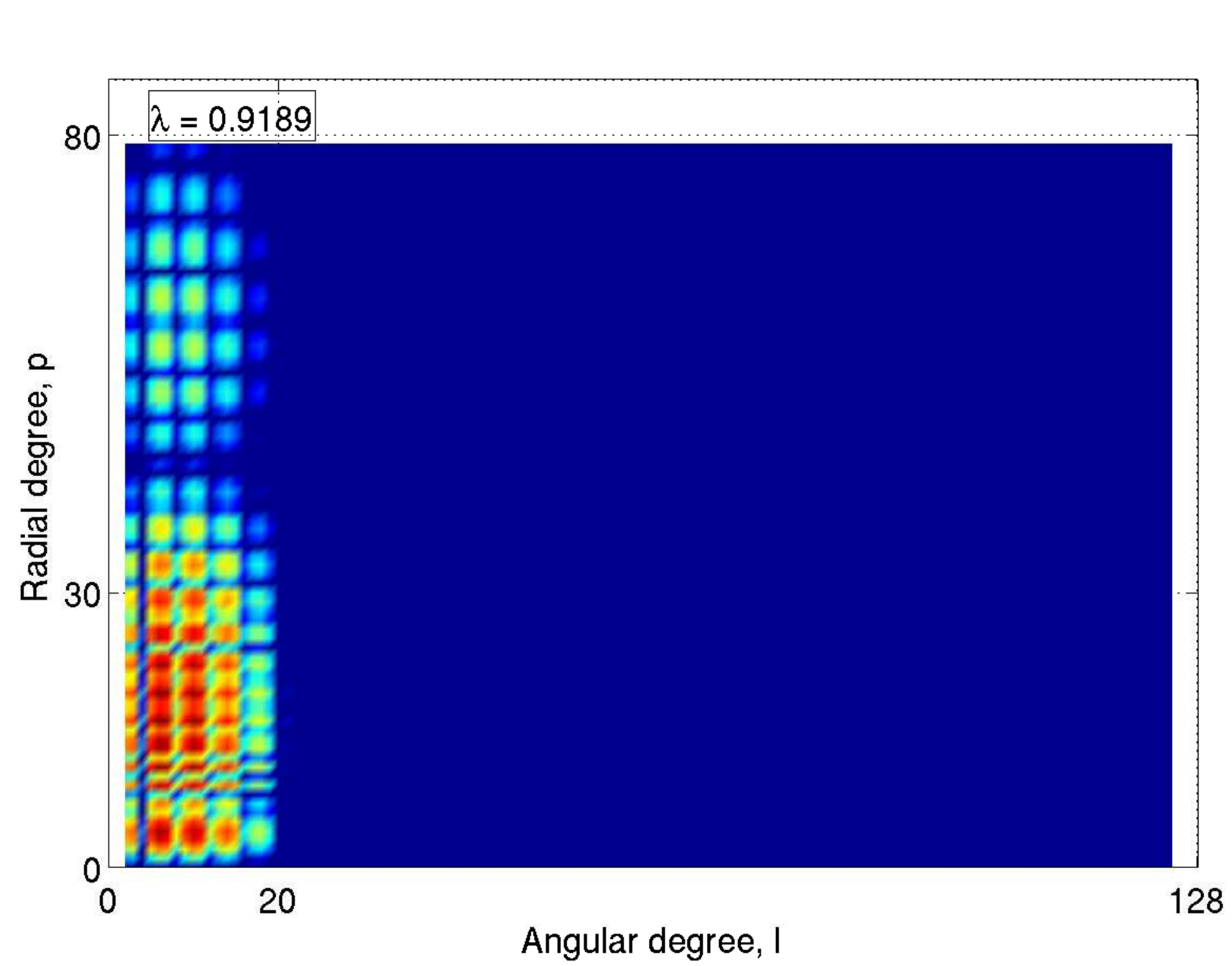}}
    \hspace{-2mm}
    \subfloat{
        \includegraphics[scale=0.215]{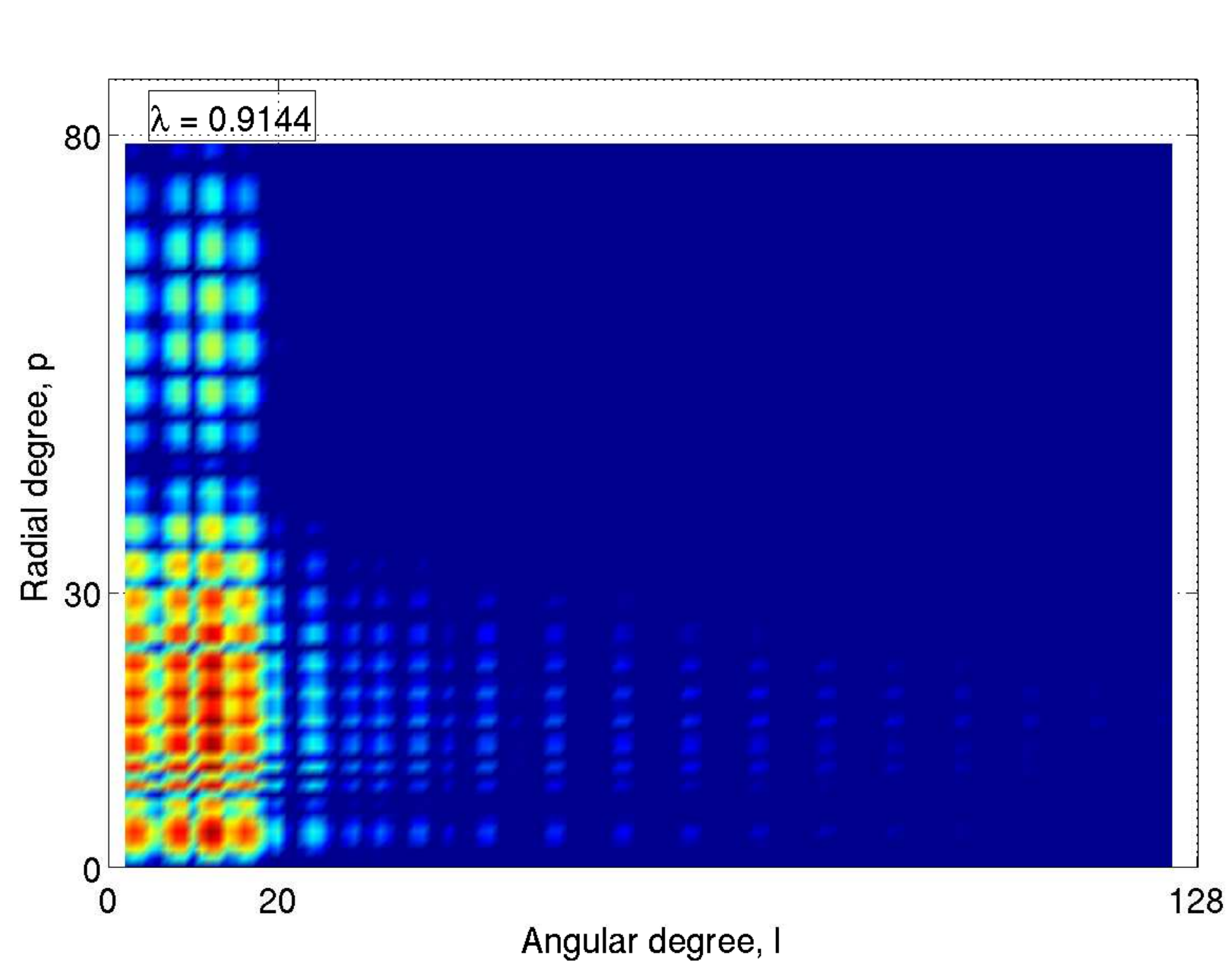}}
    \hspace{-2mm}
    \subfloat{
        \includegraphics[scale=0.215]{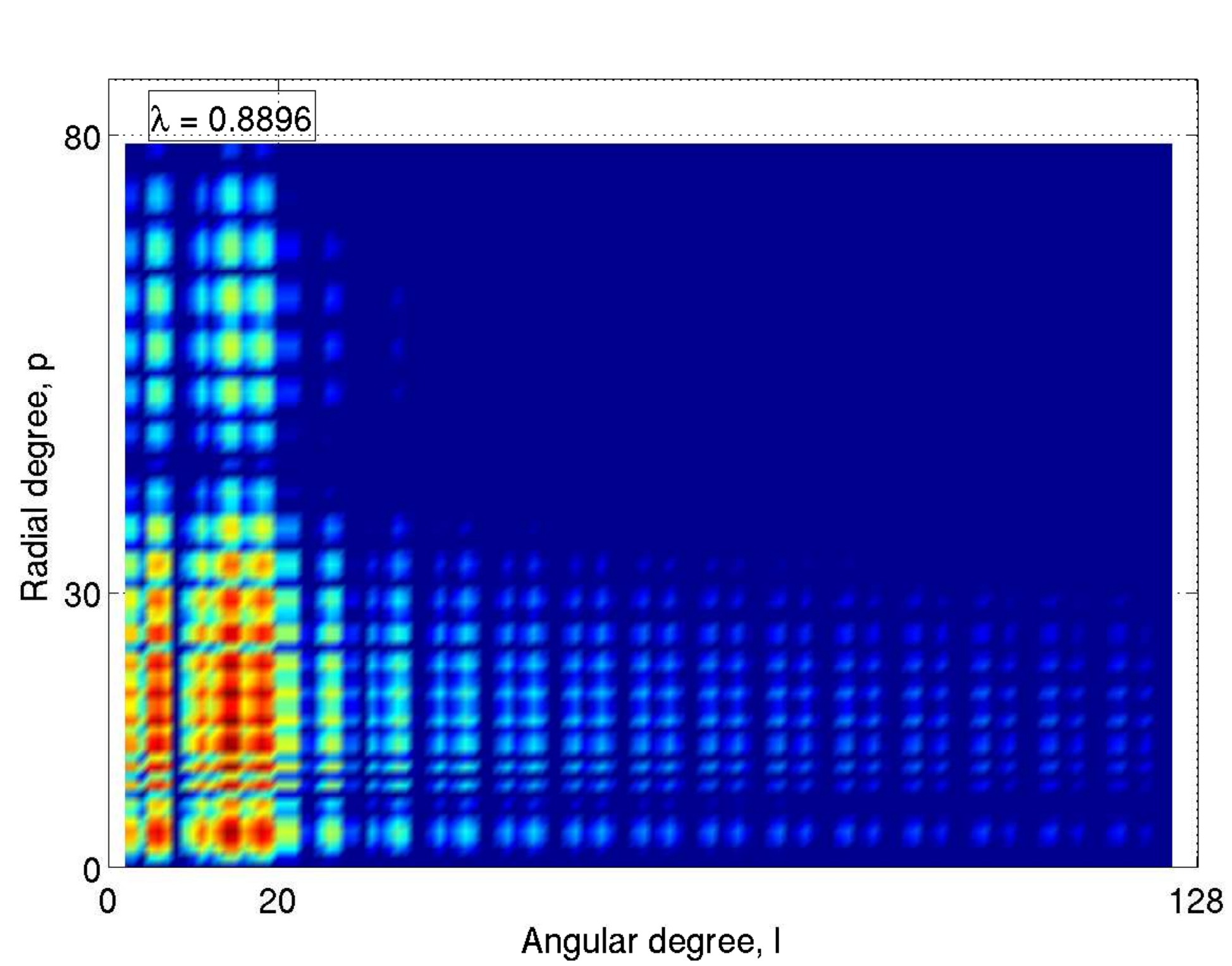}}


    \vspace{-3mm}

    \subfloat{
        \includegraphics[scale=0.28]{colorbar_spect}}

    \vspace{-3mm}

    \caption{\textbf{Fourier-Laguerre spectral domain response of
        space-limited spectrally concentrated eigenfunctions}, given
      by their magnitude square $|g_{\ell m p}|^2$.  Each
      eigenfunction is limited in the azimuthally symmetric and
      radially independent spatial region $R=\{ 15\leq r \leq 25,\,
      \pi/8 \leq \theta \leq 3\pi/8, \, 0 \leq \phi <2\pi \}$ and
      spectrally concentrated in the Fourier-Laguerre domain within the
      spectral region $A_{(30,20)}$. The eigenvalue $\lambda$
      associated with each eigenfunction is a measure of the spectral
      concentration within the Fourier-Laguerre spectral domain
      $A_{(30,20)}$. The values of $|g_{\ell m p}|^2$ are
      plotted in decibels as $20\log|g_{\ell m}(k)|$, normalized to
      zero at the individual maxima of each
      eigenfunction. } \label{fig:eigen_laguerre_spect}
\end{figure*}

\subsection{Illustration}

As an example, we compute band-limited and space-limited eigenfunctions
for the azimuthally and radially symmetric spatial region $R=\{ 15\leq r
\leq 25,\, \pi/8 \leq \theta \leq 3\pi/8, \, 0 \leq \phi <2\pi\}$
and spectral regions $\tilde A_{KL} = \tilde A_{(1.4,20)}$ in the Fourier-Bessel spectral
domain and $A_{PL} = A_{(30,20)}$ in the Fourier-Laguerre spectral
domain.
The band-limit $K=1.4$ for the Fourier-Bessel domain is chosen such that, for the spatial region under consideration,
the spherical Shannon numbers $\tilde N_{KL}$ given in \eqref{Eq:sum_eigen_values_FB} and $N_{PL}$ given in \eqref{Eq:sum_eigen_values_FL} are the same, i.e.~$\tilde N_{KL}\approx N_{PL}$.  The spherical Shannon numbers for this region are plotted in \figref{fig:shannon_comparison}a and b respectively for different values of band-limits $K$ and $P$.
Since the region is azimuthally symmetric, the
eigenfunctions band-limited in Fourier-Bessel and Fourier-Laguerre
spectral domains can be found by solving the fixed-order eigenvalue
problems formulated in \eqref{Eq:symmetric1_sub_problem1} and
\eqref{Eq:symmetric1_sub_problem2} respectively. When computing
band-limited functions in the Fourier-Bessel domain, we discretize
the spectrum along $k$ with a step size of $\Delta k = 0.02$. With such
a discretization, the problem in \eqref{Eq:symmetric1_sub_problem1}
reduces to a finite dimensional algebraic eigenvalue problem. In
order to find band-limited eigenfunctions in the Fourier-Laguerre domain
as a solution of \eqref{Eq:symmetric1_sub_problem2}, we use the
decomposition into the two subproblems formulated in
\eqref{Eq:symmetric2_FL_subproblem1} and
\eqref{Eq:symmetric2_FL_subproblem2} to separately maximize the
radial and angular concentration respectively. From the
band-limited eigenfunction $f$, the space-limited eigenfunction $g$ can be obtained by limiting the band-limited eigenfunction $f$ in the spatial domain using \eqref{Eq:ftog_in_spatial}. Alternatively, the space-limited eigenfunction $g$ can be obtained in the spectral domain from the
band-limited eigenfunction $f$ using \eqref{Eq:relation_g_f_spectral1}.

The spectrum of eigenvalues is plotted in
\figref{fig:eigen_spectrum_lag_sphere}a and
\figref{fig:eigen_spectrum_lag_sphere}b for the
Fourier-Bessel and Fourier-Laguerre problems respectively, where the
spherical Shannon number $\tilde N_{KL}$ or $N_{PL}$ is also
indicated, which confirms that the spherical Shannon number is indeed a
good approximation of the number of significant eigenvalues. It is worth noting that $\tilde N_{KL}$ is computed using the analytic expressions \eqref{Eq:sum_eigen_values_FB}, which illustrates that the number of significant eigenfunctions is independent of the choice of discretization along $k$ and thus the infinite dimensional subspace $\lballss{K}{L}$ is, in fact, spanned by the finite $\tilde N_{KL}$ eigenfunctions. The
spectrum of eigenvalues $\lambda^1$ and $\lambda^2$ is plotted in
\figref{fig:eigen_spectrum_lag_sphere}c and
\figref{fig:eigen_spectrum_lag_sphere}d respectively, where we have
also indicated $N^P$ and $N_L$. For the concentration problem in the Fourier-Laguerre domain, the measure of radial concentration~($\lambda^1$) is independently controlled by the radial band-limit $P$ and angular concentration ($\lambda^2$) is controlled by the angular band-limit $L$. Consequently, concentration in the spatial region $R\subset\bll$ is given by $\lambda = \lambda^1\lambda^2$, due to which, the transition of the spectrum of eigenvalues $\lambda$ from unity to zero is not smooth. However, the spectrum of each $\lambda^1$ and $\lambda^2$ does have a smooth transition. For the problem in the Fourier-Bessel domain, the transition is smooth due to the spectral coupling between radial and angular components, as discussed earlier.

For the Fourier-Bessel spectral domain,
the first $12$ most concentrated band-limited eigenfunctions $f^{(m)}(r,\theta)$ given in \eqref{Eq:spatial_eigen_subproblem1} for $m=2$
are shown in the spatial domain in \figref{fig:eigen_bessel_spatial}.  The associated space-limited functions concentrated in the spectral domain are shown in \figref{fig:eigen_bessel_spect}, computed using \eqref{Eq:relation_g_f_spectral1} (note that these functions are not band-limited and so their spectral representation is non-zero outside of the plotted domain).
For the Fourier-Laguerre domain, the
first $12$ most concentrated band-limited eigenfunctions $f^{(m)}(r,\theta)$ given in \eqref{Eq:spatial_eigen_subproblem2} for $m=2$
are shown in the spatial domain in \figref{fig:eigen_laguerre_spatial}.   The associated space-limited functions concentrated in the spectral domain are shown in \figref{fig:eigen_laguerre_spect}, computed using \eqref{Eq:relation_g_f_spectral1} (again, note that these functions are not band-limited). For space-limited eigenfunctions, it can be noted that the leakage of concentration~(energy) outside the Fourier-Bessel spectral region of interest is diagonally spread due to the spectral coupling between the radial and angular components for the Fourier-Bessel functions. In contrast, the leakage outside the Fourier-Laguerre spectral region of interest has horizontally and/or vertically spread due to the separability of the radial and angular components of the Fourier-Laguerre functions.

\section{Representation of Signal in Slepian Basis}\label{sec:Applications}

The spatially localised band-limited Fourier-Bessel or
Fourier-Laguerre eigenfunctions form a complete orthonormal basis for
signals band-limited in the respective domain.  We call the
band-limited eigenfunctions that arise as solutions of the eigenvalue
problems \eqref{Eq:eig_harmonic_spectral1} and
\eqref{Eq:eig_harmonic_spectral2_matrix_formulation} a \emph{Slepian
basis}. By completeness, any band-limited signal can be represented in
the Slepian basis constructed with the same band-limit.  A
Fourier-Bessel band-limited signal $h \in
\mathcal{\tilde H}_{KL}\subset \lball$ can thus be represented by
\begin{align}\label{Eq:decomposition_FB}
h(\bv{r}) &= \int_{\rplus} \ddr{k} \sum_{\ell=0}^{L-1}
\sum_{m=-\ell}^\ell {h}_{\ell m}(k) X_{\ell m}(k,\bv{r}) \nonumber \\ &= \sum_{\alpha=1}^\infty \shcc{h}{\alpha} f^\alpha(\bv{r}),
\quad f^\alpha,\, h  \in \mathcal{\tilde H}_{KL}
\end{align}
and a Fourier-Laguerre band-limited signal $h  \in \mathcal{H}_{PL}
\subset \lball$ can be represented by
\begin{align}\label{Eq:decomposition_FL}
h(\bv{r}) &= \sum_{p=0}^{P-1} \sum_{{\ell}=0}^{L-1}
        \sum_{m=-{\ell}}^{\ell} \shc{h}{\ell m p}  Z_{\ell m p}(\bv{r}) \nonumber \\ &= \sum_{\alpha=1}^{PL^2} \shcc{h}{\alpha} f^\alpha(\bv{r}),
\quad f^\alpha,\, h  \in \mathcal{H}_{PL},
\end{align}
where $\shcc{h}{\alpha}$ denotes the Slepian coefficient given by
\begin{align}\label{Eq:decomposition_coefficient}
\shcc{h}{\alpha} = \int_{\bll}\ddv{r} h(\bv{r}) \left(f^\alpha(\bv{r})\right)^\ast
\end{align}
(we use a single subscript to denote Slepian coefficients, similar to
the notation used to denote spherical Laguerre coefficients; however,
the case intended should be obvious from the context). Note that the signal $h$ and the eigenfunctions $f^\alpha$ belong to the same subspace.
We note that the Slepian coefficient in
\eqref{Eq:decomposition_coefficient} can be also be obtained through the spectral domain representation of both signal
and Slepian basis function, that is, by
\begin{align}\label{Eq:eigen_coefficient_from_spectral}
\shcc{h}{\alpha} &= \sum_{\ell=0}^{L-1} \sum_{m=-\ell}^\ell \int_{k=0}^K \ddr{k} \, h_{\ell m}(k)  \left(\fbc{f}{\ell m}{\alpha}(k)\right)^\ast,\quad\, f^\alpha,\, h  \in \mathcal{\tilde H}_{KL}
\nonumber
\\
\shcc{h}{\alpha} &= \sum_{p=0}^{P-1} \sum_{\ell=0}^{L-1} \sum_{m=-\ell}^\ell  h_{\ell mp}  \left(f_{\ell m p}^\alpha\right)^\ast,\quad\, f^\alpha,\, h  \in \mathcal{H}_{PL}.
\end{align}

When the signal is spatially localized and concentrated in some
spatial region $R$, the summation over the Slepian basis functions in
\eqref{Eq:decomposition_FB} and \eqref{Eq:decomposition_FL} can be
truncated at the spherical Shannon numbers $\tilde N_{KL}$ and
$N_{PL}$, respectively, since the Shannon numbers give the
approximate number of eigenfunctions concentrated in the region $R$.
Consequently, a signal concentrated in $R$ can be approximated
accurately by
\begin{align}\label{Eq:decomposition_approximation}
h(\bv{r}) &\approx \sum_{\alpha=1}^{\tilde N_{KL}} \shcc{h}{\alpha} f^\alpha(\bv{r}),
\quad  h  \in \mathcal{\tilde H}_{KL},
\\
h(\bv{r}) &\approx \sum_{\alpha=1}^{N_{PL}} \shcc{h}{\alpha} f^\alpha(\bv{r}),
\quad  h  \in \mathcal{H}_{PL},
\end{align}
where the same region $R$ is used to construct the eigenfunctions.  It
also follows that the Slepian coefficients associated with the
remaining eigenfunctions (which are concentrated in the region $\bll
\backslash R$) are close to zero.  Consequently, a concentrated signal
is \emph{sparse}\footnote{Formally the signal is \emph{compressible} since
  the Slepian coefficients associated with the remaining
  eigenfunctions are small and not identically zero.} when represented
in the Slepian basis, where the degree of sparsity is quantified by
the spherical Shannon number.  The representation of a concentrated
signal in the Slepian basis, truncated at the spherical Shannon
number, can thus give considerable savings in terms of the number of
coefficients required to represent the signal accurately.  The
accuracy of the approximate representation of the signal within the
spatial region $R$ can be quantified by the ratio of the energy of the
approximate representation to the energy of the exact representation
within the region $R$. We define such an accuracy measure by
\begin{align}\label{Eq:Quality_approximation}
Q(J) &= \frac{ \mathlarger{\int}_R \ddv{\bv{r}} \bigg|\sum\limits_{\alpha=1}^{J} \shcc{h}{\alpha} f^\alpha(\bv{r})\bigg|^2 }{\mathlarger{\int}_R \ddv{\bv{r}} \bigg| \sum\limits_{\alpha=1}^{\infty} \shcc{h}{\alpha} f^\alpha(\bv{r})\bigg|^2}
= \frac{\sum\limits_{\alpha=1}^{J} \lambda_\alpha \big|\shcc{h}{\alpha}\big|^2}{\sum\limits_{\alpha=1}^{\infty} \lambda_\alpha \big|\shcc{h}{\alpha}\big|^2},
\end{align}
where $J = \tilde N_{KL}$ for $f^\alpha,\, h \in \mathcal{\tilde
  H}_{KL}$ and $J = N_{PL}$ for $f^\alpha,\, h \in \mathcal{H}_{PL}$.
In obtaining the second equality we have used the orthogonality
relationships given in \eqref{Eq:properties_eig_FB_spatial} and
\eqref{Eq:properties_eig_FL_spatial}.  The closer the value of $Q(J)$
to unity, the greater the accuracy of the approximate representation of
signal within the spatial region $R$.  Since $\lambda_\alpha \approx
1$ for $\alpha\leq J$ and $\lambda_\alpha \approx 0$ for $\alpha>J$,
the accuracy of the approximation within the spatial region $R$ is
high.

\begin{figure*}[!t]
    \centering
    \subfloat[Angular region]{
        \includegraphics[scale=0.35]{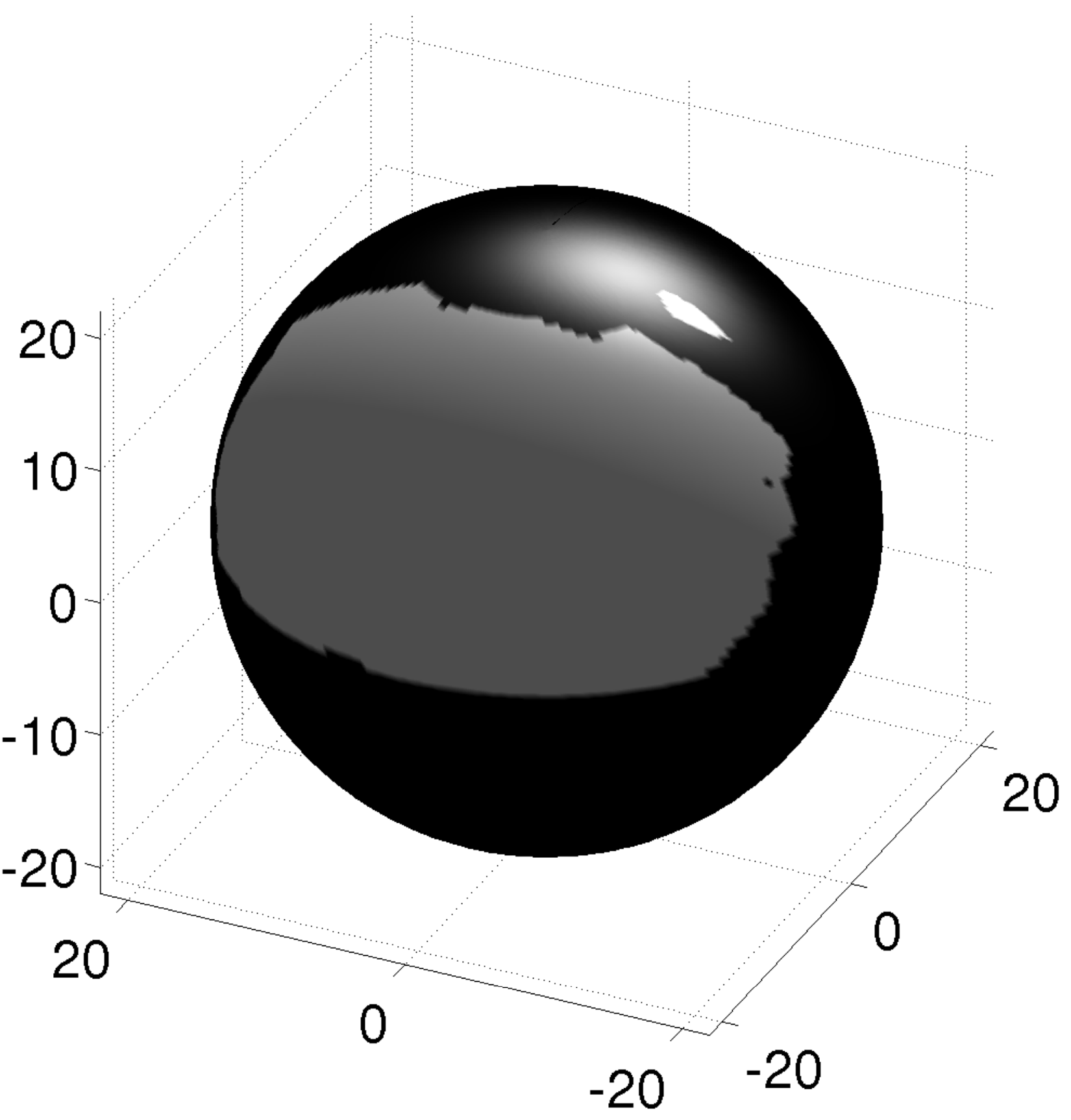}}
    \hspace{8mm}
    \subfloat[Radial region]{
        \includegraphics[scale=0.28]{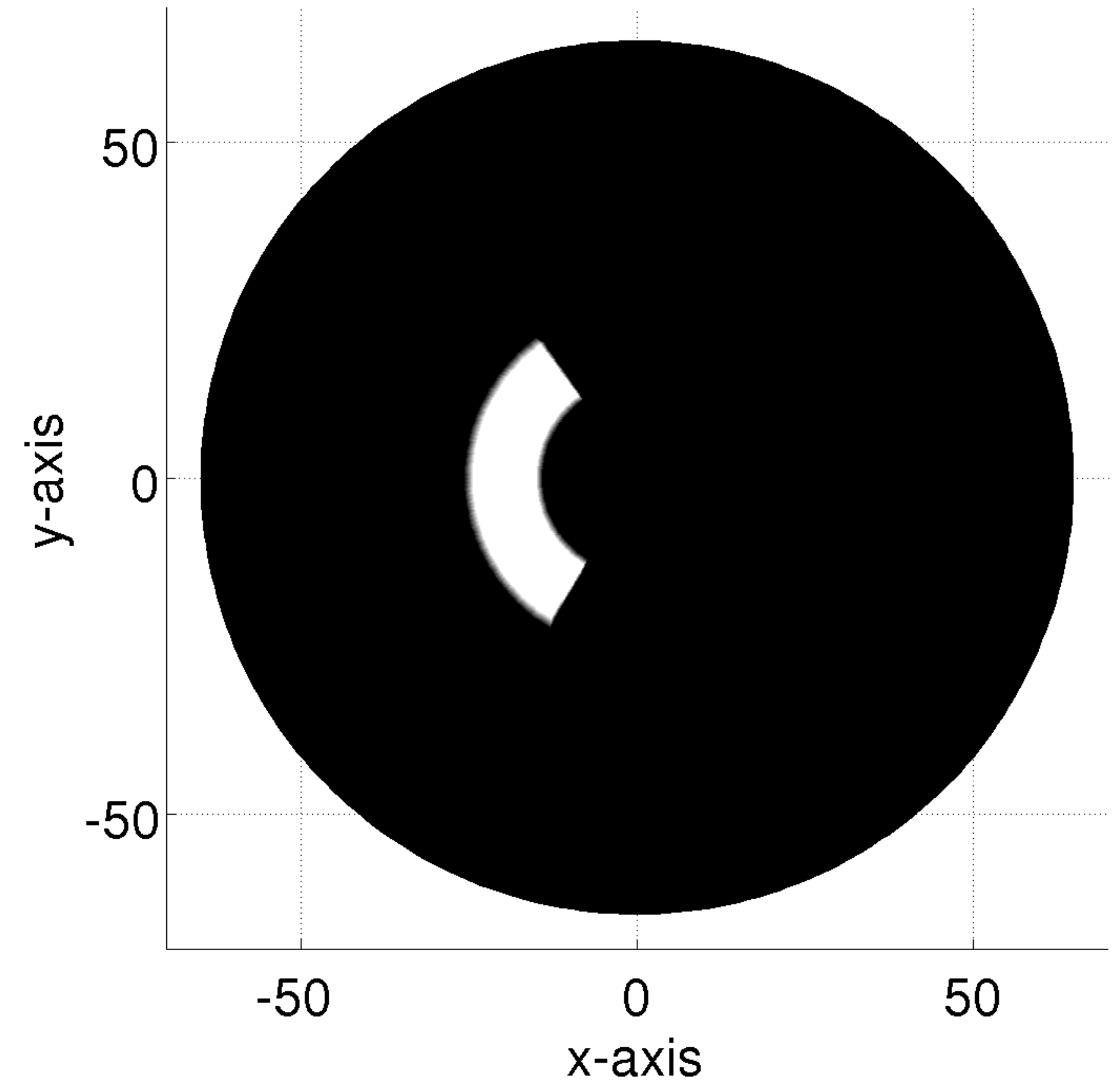}}
 \vspace{5mm}

    \hspace{-15mm}
    \caption{\textbf{Non-trivial spatial region} constructed from (a) the
      SDSS DR7 quasar mask on the sphere (shown on the sphere of
      radius $r=20$) and (b) a radial profile limited to the
      interval $15 \leq r \leq 25$ (shown in the $x$-$y$ plane).}
    \label{fig:mask_spatial}
\end{figure*}

\begin{figure*}[!t]
    \centering
    \subfloat[Slices through the ball]{
        \includegraphics[scale=0.37]{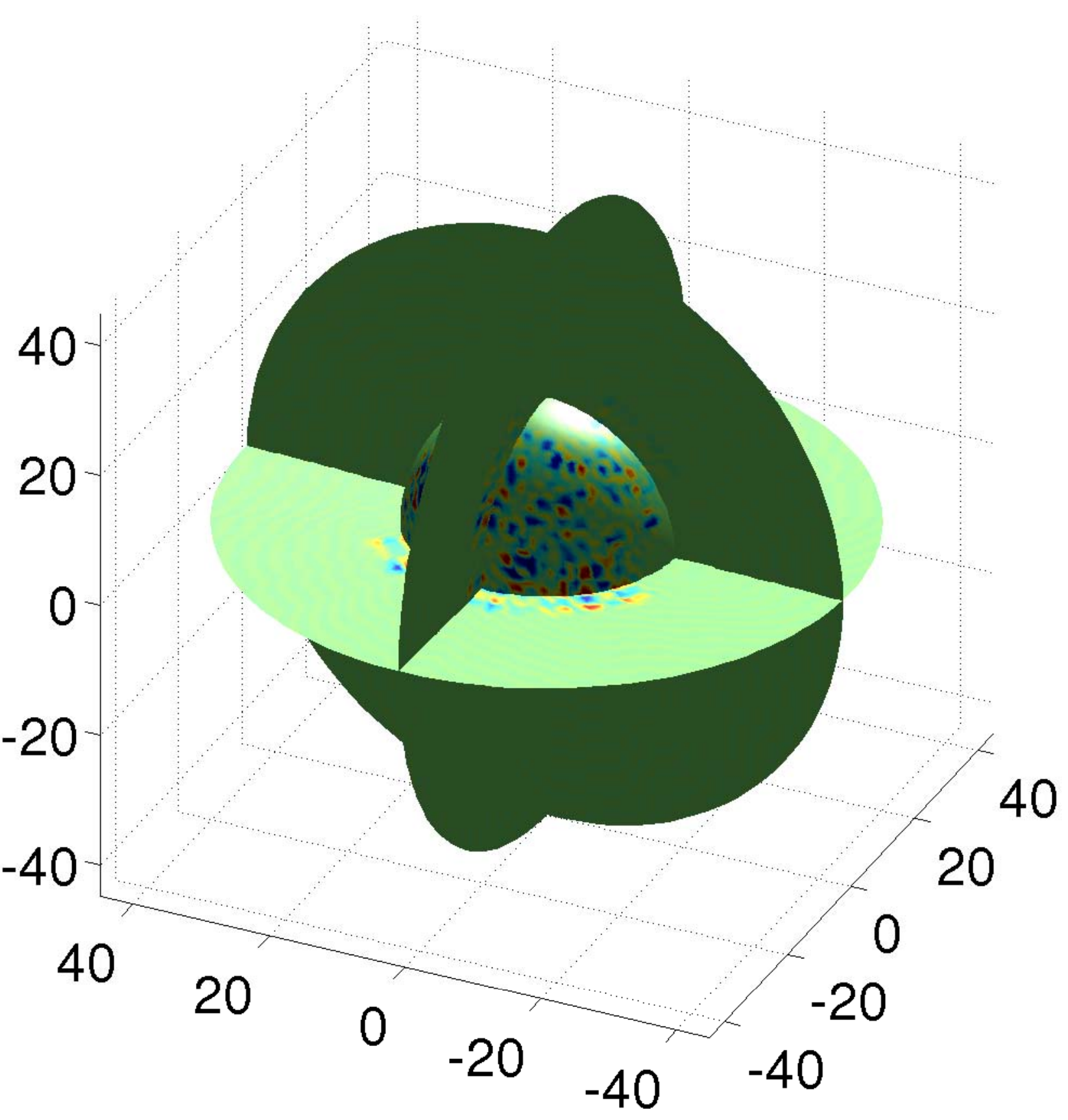}}
    \hspace{3mm}
    \subfloat[Spherical shell]{
        \includegraphics[scale=0.37]{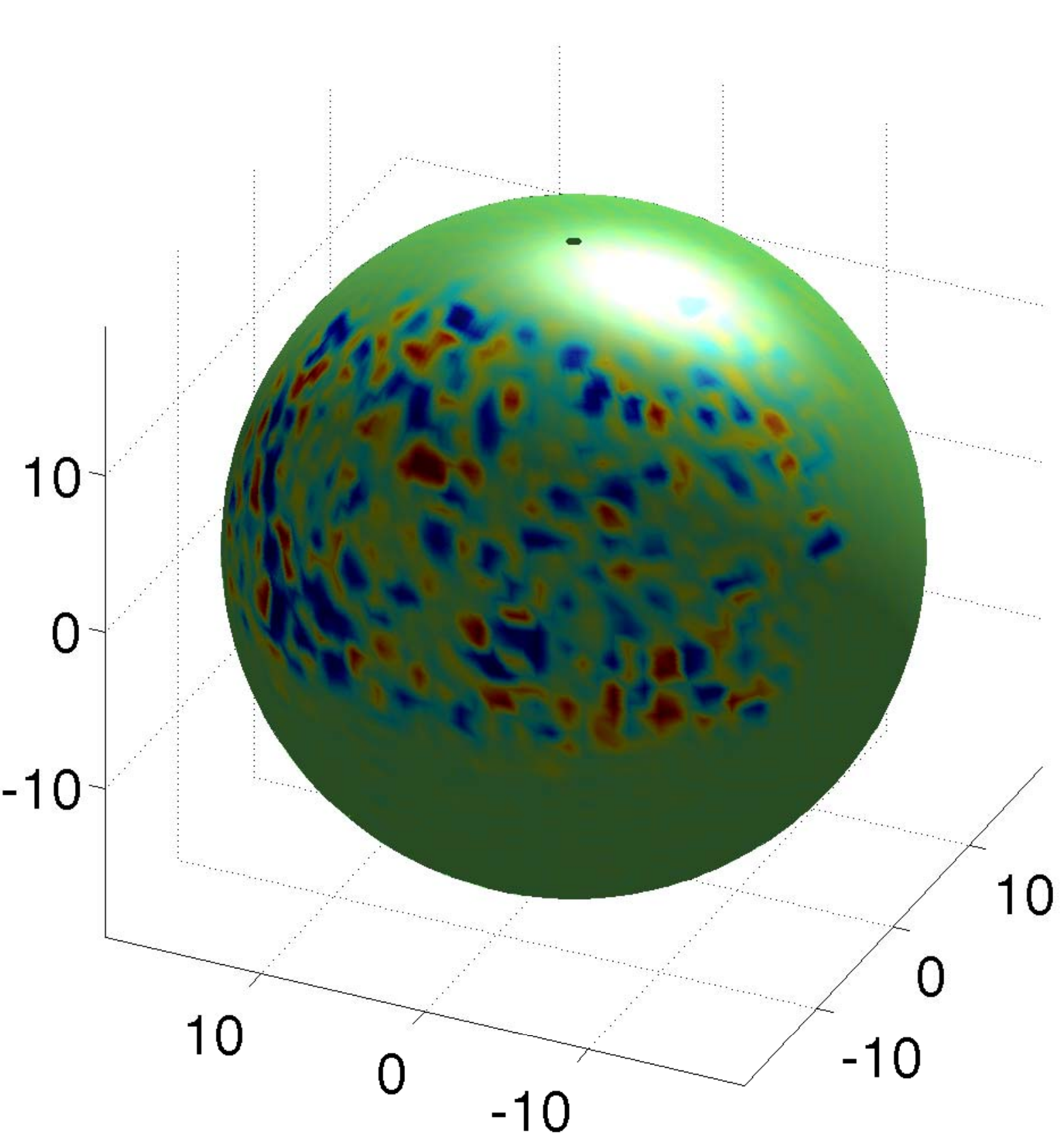}}
 \vspace{5mm}

    \subfloat{
        \includegraphics[scale=0.3]{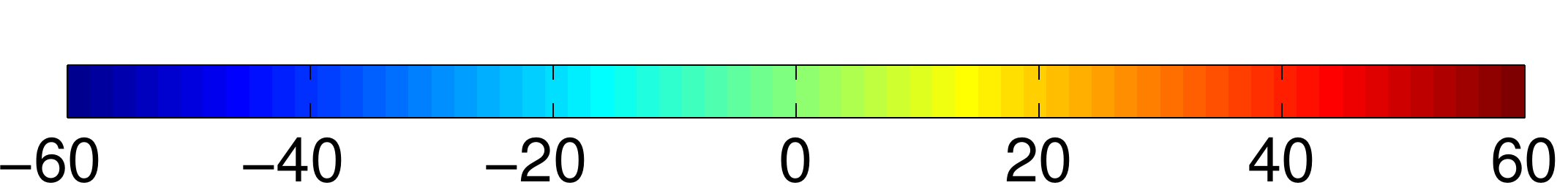}}

    \hspace{-15mm}
    \caption{\textbf{Simulated dark matter density test signal on the
        ball}, band-limited and spatially concentrated in the region
      shown in \figref{fig:mask_spatial}.  The test signal $h \in
      \mathcal{H}_{PL} \subset \lball$ is band-limited in the
      Fourier-Laguerre domain at angular degree $L=72$ and radial
      degree $P=144$, and spatially concentrated in the angular region
      defined by the SDSS DR7 quasar mask and radially in the interval
      $15 \leq r \leq 25$. The test signal is plotted on (a) slices
      through the ball and (b) on a spherical shell at $r=20$. }
    \label{fig:signal_spatial}
\end{figure*}

\subsection{Illustration}

We present an example to illustrate the fact that the representation
of a spatially concentrated band-limited signal in the Slepian basis
is sparse. We consider a test signal obtained from a simulation of the
dark matter distribution of the Universe, observed over a partial
field-of-view.  Specifically, the test signal is extracted from the
full-sky Horizon Simulation~\cite{Teyssier}, an $N$-body simulation
that covers a 1Gpc periodic box of 70 billion dark matter particles
generated from the cosmological concordance model derived from 3-year
Wilkinson Microwave Anisotropy Probe~(WMAP)
observations~\cite{Spergel:2007}.  We limit the data in the spatial
domain with an Sloan Digital Sky Survey (SDSS)
DR7\footnote{\url{http://www.sdss.org/dr7/}} quasar binary mask on the
sphere, denoted by $\rsphere \subset \untsph$, and in the interval
$15\leq r \leq 25$ along the radial line (note that the units of
radius for the data under consideration are Mpc). The resulting mask
is shown in \figref{fig:mask_spatial}. The masked signal is
band-limited at $L=72$ and $P=144$ to obtain the
spatially concentrated band-limited test signal $h \in
\mathcal{H}_{PL}$ shown in \figref{fig:signal_spatial}.

\begin{figure}[t]
    \centering
    \includegraphics[scale=0.67]{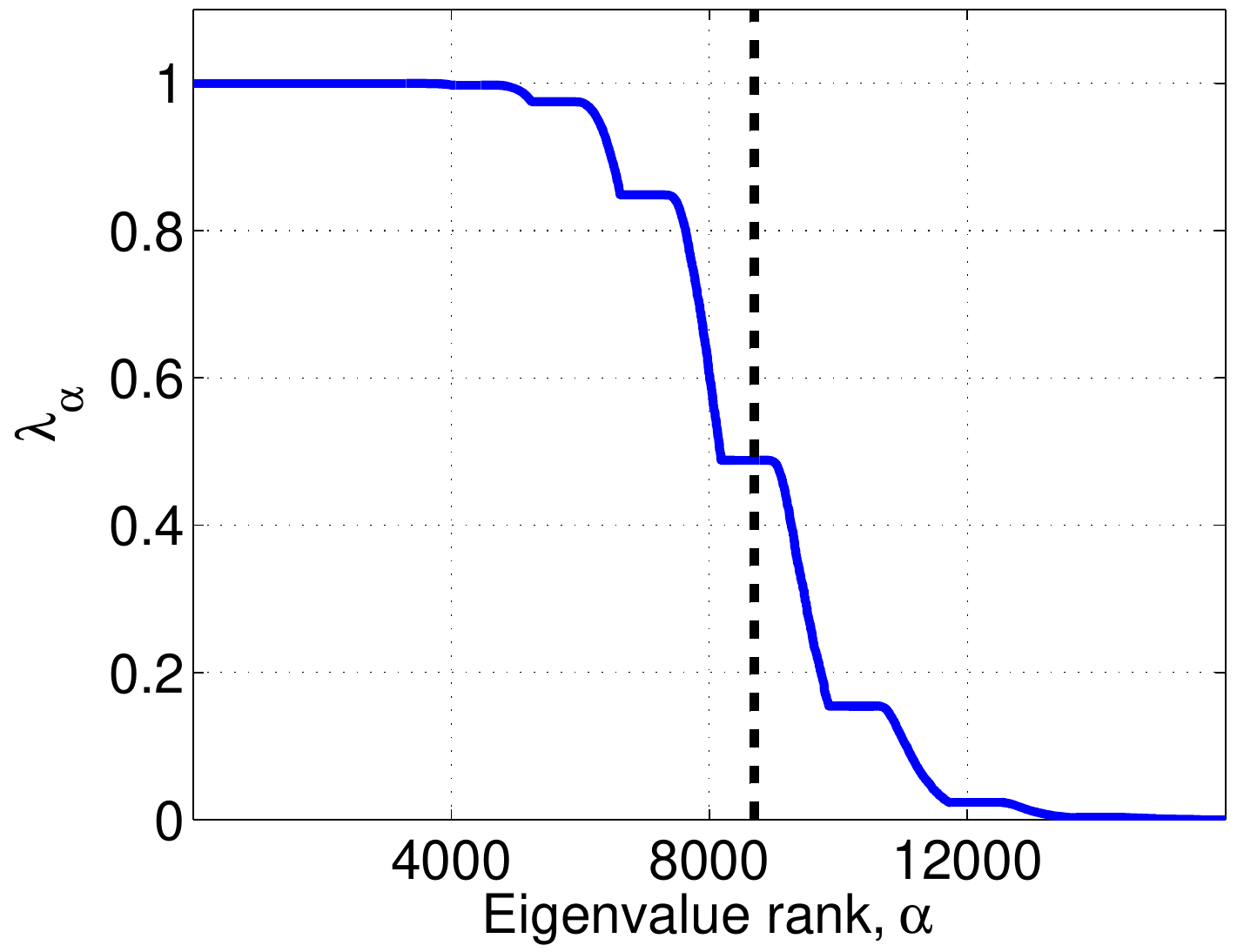}
    \caption{\textbf{Spectrum of eigenvalues $\lambda_\alpha$ for the
        non-trivial region} associated with the band-limited spatially
      concentrated eigenfunctions. The angular band-limit $L=72$ and
      radial band-limit $P=144$ are assumed, while the spatial region
      shown in \figref{fig:mask_spatial} is considered, i.e.\ defined
      by the SDSS DR7 quasar mask on the sphere and the radial
      interval $15 \leq r \leq 25$. The horizontal axis is truncated
      to show only the first 16,000 eigenvalues~(out of the total of
      746,496 eigenvalues). The spherical Shannon number $\lfloor
      N_{PL} \rfloor = 8,696$ is also indicated by the vertical dashed
      line. Note, again, that the transition of the eigenvalue
      spectrum from unity to zero is not smooth since the radial and
      angular components are interspersed.  The spherical Shannon
      number again estimates the location of the transition in the
      eigenvalue spectrum accurately.}
    \label{fig:illustration_eigen_values}
\end{figure}

The region considered is not symmetric but is radially independent,
hence we compute the band-limited eigenfunctions separately on the
radial line over the interval $15\leq r \leq 25$ with band-limit
$P=144$ and on the sphere over the region $\rsphere$ with band-limit
$L=72$. In total, $PL^2 = 746,496$ eigenfunctions $f^\alpha \in
\mathcal{H}_{PL}$ are computed. The spectrum of eigenvalues is plotted
in \figref{fig:illustration_eigen_values} for the first 16,000
eigenvalues only, where the spherical Shannon number $\lfloor N_{PL}
\rfloor = 8,696 $ is also indicated, which is computed using
\eqref{Eq:Shannon_symmetric_not_radial}.

\begin{figure}[t]
    \centering
    \includegraphics[scale=0.72]{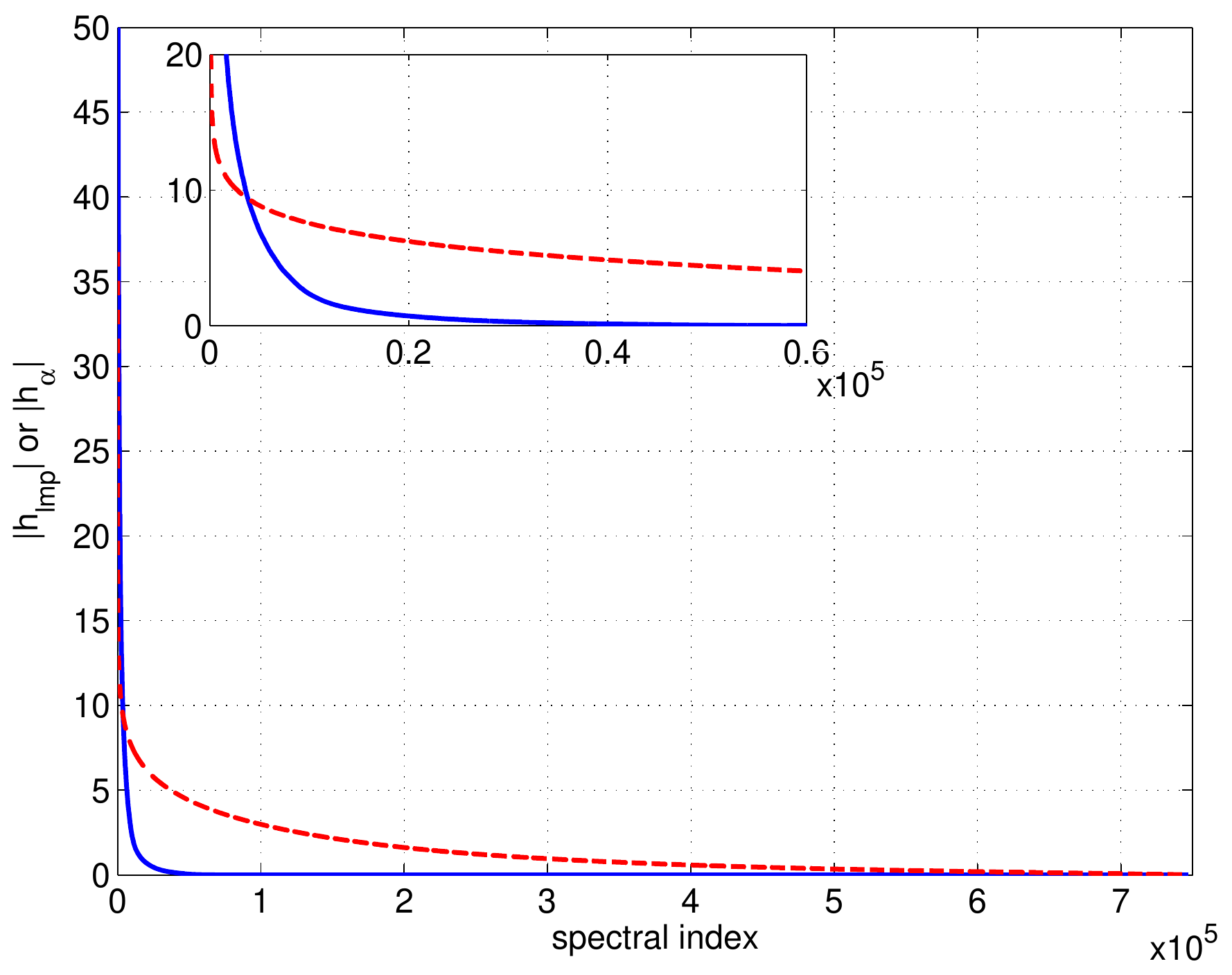}
    \caption{\textbf{Spectral decay of the Fourier-Laguerre and
        Slepian coefficients} of the band-limited spatially
      concentrated test signal shown in
      \figref{fig:signal_spatial}. The absolute values of both the
      Fourier-Laguerre $h_{\ell m p}$ and Slepian $\shcc{h}{\alpha}$
      coefficients are first sorted in descending order and then
      plotted. Since the band-limited signal $h$ is spatially
      concentrated, it has a sparse representation in the Slepian basis.
      }
    \label{fig:illustration_eigen_expansion}
\end{figure}

Since the eigenfunctions $f^\alpha$ serve as a complete basis for the
subspace $\mathcal{H}_{PL}$, the signal $h$ can be alternatively
represented in the Slepian basis as given in
\eqref{Eq:decomposition_FL}, where we compute the Slepian coefficients
$\shcc{h}{\alpha}$ using
\eqref{Eq:eigen_coefficient_from_spectral}. Recall that
$\mathcal{H}_{PL}$ is a $PL^2$ dimensional subspace, thus there are
$PL^2$ spectral coefficients $h_{\ell m p}$ or Slepian coefficients
$\shcc{h}{\alpha}$.  The decay of the Fourier-Laguerre coefficients
$h_{\ell m p}$ and the Slepian coefficients $\shcc{h}{\alpha}$ is
compared in \figref{fig:illustration_eigen_expansion}, where both are
first sorted and then plotted.  As expected, the Slepian coefficients
decay much more rapidly than the Fourier-Laguerre coefficients.  The
spatially concentrated signal indeed has a very sparse representation
in the Slepian basis; it can be represented accurately within the
spatial region of interest with only $\lfloor N_{PL} \rfloor = 8,696$ Slepian
coefficients, compared to the full dimensionality of the subspace of
$PL^2 = 746,496$.  The energy ratio captured by the approximate
representation, defined by the accuracy measure of
\eqref{Eq:Quality_approximation}, is $Q(N_{PL}) =99.73 \%$, again
indicating that the approximation is very accurate within the spatial
region of interest.

Our code to solve the eigenvalue problems that result from the Slepian
spatial-spectral concentration problem on the ball, and used to
performed the illustrations presented here, will be made public
following the review of this manuscript.  In the interim, the code
can be obtained from the authors.

\section{Conclusions}\label{sec:conclusions}

We have formulated and solved the Slepian spatial-spectral
concentration problem on the ball. We consider two domains for the
spectral characterization of the signal defined on the ball, namely
the Fourier-Bessel domain and the Fourier-Laguerre domain.  The
Fourier-Laguerre domain is considered in addition to the standard
Fourier-Bessel domain since the former has a number of practical
advantages. The orthogonal families of band-limited spatially
concentrated functions and of space-limited spectrally concentrated
functions can be computed as solutions of eigenvalue problems. The
spatially and spectrally concentrated eigenfunctions that arise as
solutions of these eigenvalue problems coincide with each other inside
both the spatial and spectral regions of interest. The eigenvalue
associated with each eigenfunction is a measure of both the spatial
concentration of the band-limited function and of the spectral
concentration of the space-limited function. The number of
well-concentrated~(significant) eigenfunctions depends on the
spherical Shannon number, which also serves as the dimension of the
space of functions that can be concentrated in both the spatial and
spectral domains at the same time. When the spatial region of interest
is rotationally symmetric and/or radially independent, the eigenvalue
problem decomposes into subproblems, which reduces the computational
burden significantly.  The family of concentrated eigenfunctions can
be used to form an orthonormal basis, or Slepian basis, which provides
a sparse representation of concentrated functions.

Just as the Slepian basis in the one-dimensional Euclidean domain and
other geometries have proven to be extremely valuable, we hope the
Slepian eigenfunctions on the ball developed in this work will prove
useful in a variety of applications in various fields of science and
engineering~(e.g., geophysics, cosmology and planetary science), where
signals are inherently defined on the ball. Some applications where
the proposed orthogonal family of eigenfunctions~(Slepian basis) are
likely to be of use are the estimation of signals, their power
spectrum, and other statistics, when noisy observations on the ball
can only be made over partial fields-of-views.  For example, surveys
of the distribution of Galaxies in our Universe are observed only over
partial fields-of-view \cite{ahn:2012} on the ball, while their
statistical properties can be used to infer the physics of our
Universe, such as the nature of dark energy and dark matter.

\section*{Acknowledgements}
We thank Hiranya Peiris and Boris Leistedt for useful discussions, for
providing the processed Horizon simulation data-set, and for their
hospitality during the visit of Z.~Khalid to UCL.

\bibliography{bib_Zubair} 

\end{document}